\documentclass[11pt,twoside,leqno]{aomamlt2e} 
  \pageno{1187}
\received{December 4, 2000}
  \renewcommand{\thesection}{\arabic{section}}
\makeatletter
\renewcommand{\@seccntformat}[1]{\csname
the#1\endcsname.\hspace{0.5em}\setcounter{Subsec}{0}\setcounter{Subsubsec}{0}}\makeatother

 \newcommand{\ee}{{\hskip1pt\rm \'{\hskip-6.5pt \it e}}}
 
\newtheorem{theorem}{Theorem}[section]

\newtheorem{lemma}[theorem]{Lemma}
\newtheorem{proposition}[theorem]{Proposition}

\theorembodyfont{\upshape}
\newtheorem{definition}[theorem]{{\it Definition}}

\newcommand\lp{\left(}
\newcommand\rp{\right)}
\newcommand\lb{\left\{}
\newcommand\rb{\right\}}
\newcommand\lk{\left[}
\newcommand\rk{\right]}

\newcommand\C{{\mathbb C}}
\newcommand\D{{\mathbb D}}
\newcommand\N{{\mathbb N}}
\newcommand\T{{\mathbb T}}
\renewcommand\P{{\mathbb P}}
\newcommand\R{{\mathbb R}}
\newcommand\Z{{\mathbb Z}}
\newcommand\cA{{\mathcal A}}
\newcommand\cB{{\mathcal B}}
\newcommand\cC{{\mathcal C}}
\renewcommand\cD{{\mathcal D}}

\renewcommand\cH{{\mathcal H}}

\newcommand\cM{{\mathcal M}}
\newcommand\cN{{\mathcal N}}
\newcommand\cO{{\mathcal O}}

\renewcommand\cR{{\mathcal R}}
\newcommand\cS{{\mathcal S}}
\newcommand\cU{{\mathcal U}}
\newcommand\cV{{\mathcal V}}
\newcommand\cW{{\mathcal W}}
\newcommand\cZ{{\mathcal Z}}
\newcommand\rD{{\rm D}}

\newcommand\rR{{\rm R}}

\newcommand{\isdef}{\stackrel{\textrm{def}}{=}}
\newcommand\downto{\downarrow}

\newcommand\eps{\varepsilon}
\newcommand\vth{\vartheta}
\newcommand\vka{\mathcal{K}}
\newcommand\pfi{\varphi}
\newcommand\comp{\circ}
\newcommand\bD{{\pt\D}}
\newcommand\Om{\Omega}

\newcommand\OMD{{\Omega\setminus\Delta}}
\newcommand{\MG}{\cM_{\Gamma}}
\newcommand\pt{\partial}
\newcommand\tv{{\hbox{\kern1pt$\cap$\kern-5pt$|$\kern3pt}}}
\newcommand\cBl{{\cB^\ell}}
\newcommand\hBl{{\hat{\cB}^\ell}}
\newcommand\cspace[1]{\mathcal{H}({#1})}

\newcommand{\bg}{{\bar\gamma}}
\newcommand{\norm}{\mathbf{N}}

\newcommand{\dist}{\mathrm{dist}}
\newcommand\CS{curve shortening}
\newcommand\cross{\mathrm{Cross}}
\newcommand\Cross{\mathrm{Cross}}

\newcommand\utb{{T^1(M)}}       
\newcommand{\ptm}{\P TM}        

\newcommand\hg{{\hat\gamma}}
\newcommand\ptx{{\pt_{x}}}
\newcommand\pty{{\pt_{y}}}

\newcommand\Lin{\mathrm{L}}
\newcommand{\hB}{\hat{\cB}}
\renewcommand{\Re}{\mathfrak{Re}}

\begin{document}
\currannalsline{162}{2005} 

 \title{Curve shortening and the
  topology\\ of closed geodesics on surfaces}

 \acknowledgements{Supported by NSF through a grant from DMS, and by
  the NWO through grant NWO-600-61-410.}
 \author{Sigurd B.\   Angenent}

 \institution{University of Wisconsin-Madison, Madison, Wisconsin\\
\email{angenent@math.wisc.edu}}

  \shorttitle{Curve shortening and geodesics}
  

\centerline{\bf Abstract}
\vskip12pt
  We study ``flat knot types'' of geodesics on compact surfaces $M^2$.
  For every flat knot type and any Riemannian metric $g$ we introduce
  a Conley index associated with the curve shortening flow on the
  space of immersed curves on $M^{2}$.  We conclude existence of
  closed geodesics with prescribed flat knot types, provided the
  associated Conley index is nontrivial.

\section{Introduction}
If $M$ is a surface with a Riemannian metric $g$ then closed geodesics
on $(M, g)$ are critical points of the length functional $L(\gamma) =
\int |\gamma'(x)|dx$ defined on the space of unparametrized $C^{2}$
immersed curves with orientation, i.e.\ we consider closed geodesics
to be elements of the space
\[
\Omega = \mathrm{Imm}(S^{1}, M)/\mathrm{Diff}_{+}(S^{1}).
\]
Here \( \mathrm{Imm}(S^{1}, M) = \{ \gamma\in C^2(S^1, M) \mid
\gamma'(\xi)\neq 0 \text{ for all $\xi\in S^1$} \} \) and
$\mathrm{Diff}_{+}(S^1)$ is the group of $C^{2}$ orientation
preserving diffeomorphisms of $S^1=\R/\Z$.  (We will abuse notation
freely, and use the same symbol $\gamma$ to denote both a convenient
parametrization in $C^2(S^1; M)$, and its corresponding equivalence
class in $\Om$.)
 
The natural gradient flow of the length functional is given by curve shortening, i.e. by the evolution equation
\begin{equation}
    \frac{\pt\gamma}{\pt t} = \frac{\pt^2\gamma}{\pt s^2}= 
    \nabla_{T}(T),\qquad
    T\isdef \frac{\pt \gamma}{\pt s}.
    \label{eq:CS}
\end{equation}

In 1905 Poincar\'e \cite{poincare} pointed out that geodesics on
surfaces are immersed curves \emph{without self-tangencies. }
Similarly, different geodesics cannot be tangent -- all their
intersections must be transverse.  This allows one to classify closed
geodesics by their number of self-intersections, or their ``flat knot
type,'' and to ask how many closed geodesics of a given ``type'' exist
on a given surface $(M, g)$. Our main observation here is that the
curve shortening flow \eqref{eq:CS} is the right tool to deal with
this question.

We formalize these notions in the following definitions (which are a
special case of the theory described by Arnol'd in \cite{Arnold}.)

\demo{Flat knots} A curve $\gamma\in\Omega$ is a flat knot if it
has no self-tangencies. Two flat knots $\alpha$ and $\beta$ are
equivalent if there is a continuous family of flat knots
$\{\gamma_{\theta} \mid 0\leq \theta\leq 1\}$ with $\gamma_0=\alpha$
and $\gamma_1=\beta$.

\demo{Relative flat knots} For a given finite collection of
immersed curves,
\[
\Gamma=\{\gamma_1, \dots , \gamma_N \}\subset \Omega,
\]
we define a \emph{flat knot relative to $\Gamma$\/} to be any
$\gamma\in\Omega$ which has no self-tangencies, and which is
transverse to all $\gamma_j\in\Gamma$.  Two flat knots relative to
$\Gamma$ are equivalent if one can be deformed into the other through
a family of flat knots relative to~$\Gamma$.

\begin{figure}[htb]
  \centering \includegraphics{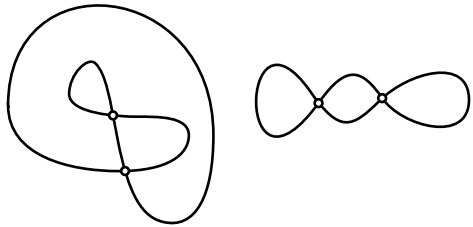}
    \caption{Two flat knots in $\R^2$ with two self-intersections}
    \label{fig:someflatknots}
\end{figure}

Clearly equivalent flat knots have the same number of self-intersections since this number cannot change
during a deformation through flat knots. The converse is not true: Flat knots with the same
number of self-intersections need not be equivalent. See Figure
\ref{fig:someflatknots}.
Similarly, two equivalent flat knots relative to $\Gamma=\{\gamma_{1}, \dots ,
\gamma_{N}\}$ have the same number of self-intersections, and the same
number of intersections with each $\gamma_j$.

In this terminology any closed geodesic on a surface is a flat knot,
and for given closed geodesics $\{\gamma_{1}, \dots , \gamma_{N}\}$
any other closed geodesic is a flat knot relative to $\{\gamma_{1},
\dots , \gamma_{N}\}$.

One can now ask the following question: Given a Riemannian metric $g$
on a surface $M$, closed geodesics $\gamma_{1}, \dots , \gamma_{N}$
for this metric, and a flat knot $\alpha$ relative to
$\Gamma=\{\gamma_{1}, \dots , \gamma_{N}\}$, \emph{how many closed
  geodesics on $( M, g)$ define flat knots relative to $\Gamma$ which
  are equivalent to $\alpha$\/}? In this paper we will use \CS\ to
obtain a lower bound for the number of such closed geodesics which
only depends on the relative flat knot $\alpha$, and the linearization
of the geodesic flow on $(T M, g)$ along the given closed geodesics
$\gamma_{j}$.

Our strategy for estimating the number of closed geodesics equivalent
to a given relative flat knot $\alpha$ is to consider the set
$\cB_\alpha\subset\Omega$ of all flat knots relative to $\Gamma$ which
are equivalent to $\alpha$.  This set turns out to be \emph{almost} an
isolating block in the sense of Conley \cite{Conley} for the \CS\ 
flow.  We then define a Conley index $h(\cB_{\alpha})$ of $\cB_\alpha$
and use standard variational arguments to conclude that nontriviality
of the Conley index of a relative flat knot implies existence of a
critical point for \CS\ in $\cB_\alpha$.

To do all this we have to overcome a few obstacles.

First, the \CS\ flow is not a globally defined flow or even semiflow.
Given any initial curve $\gamma(0)\in\Omega$ a solution $\gamma:[0,
T)\to\Omega$ to \CS\ exists for a short time $T=T(\gamma_{0})>0$, but
this solution often becomes singular in finite time.  What helps us
overcome this problem is that the set of initial curves
$\gamma(0)\in\cB_\alpha$ which are close to forming a singularity is
attracting.  Indeed, the existing analysis of the singularities of
\CS\ in \cite{Gray1989}, \cite{Ang2}, \cite{Ha:isoper}, \cite{Hu:isoper}, \cite{Oaks} shows that such
singularities essentially only form when ``a small loop in the curve
$\gamma(t)$ contracts as $t\nearrow T(\gamma(0))$.''  A calculation
involving the Gauss-Bonnet theorem shows that once a curve has a
sufficiently small loop the area enclosed by this loop must decrease
under \CS. This observation allows us to include the set of curves
$\gamma\in \cB_\alpha$ with a small loop in the exit set of the \CS\ 
flow.  With this modification we can proceed as if the \CS\ flow were
defined globally.

Second, $\cB_\alpha$ is not a closed subset of $\Omega$ and its
boundary may contain closed geodesics, i.e.\ critical points of \CS:
such critical points are always multiple covers of shorter geodesics.
To deal with this, one must analyze the \CS\ flow near multiple covers
of closed geodesics.  It turns out that all relevant information to
our problem is contained in Poincar\'e's rotation number of a closed
geodesic.  In the end our Conley index $h(\cB_{\alpha})$ depends not
only on the relative flat knot class $\cB_{\alpha}$, but also on the
rotation numbers of the given closed geodesics
$\{\gamma_1,\dots ,\gamma_N\}$.

Finally, the space $\cB_\alpha$ on which \CS\ is defined is not
locally compact so that Conley's theory does not apply without
modification.  It turns out that the regularizing effect of \CS\ 
provides an adequate substitute for the absence of local compactness
of $\cB_\alpha$.

After resolving these issues one merely has to compute the Conley
index of any relative flat knot type to estimate the number of closed
geodesics of that type.  To describe the results  we need to
discuss satellites and Poincar\'e's rotation number.

\Subsec{Satellites} Let $\alpha\in\Omega$ be given, and let
$\alpha:\R/\Z\to M$ also denote a constant speed parametrization of
$\alpha$.  Choose a unit normal $\norm$ along $\alpha$, and consider
the curve $\alpha_\epsilon:\R/\Z\to M$ given by
\[
\alpha_{\epsilon}(t) = \exp_{\alpha(qt)}\bigl(\epsilon \sin(2\pi p
t)\norm(qt)\bigr)
\]
where $\frac{p}{q}$ is a fraction in lowest terms.  When $\epsilon=0$,
$\alpha_{\epsilon}$ is a $q$-fold cover of $\alpha$.\break\vskip-11pt\noindent  For sufficiently
small $\epsilon\neq0$ the $\alpha_{\epsilon}$ are flat knots relative
to $\alpha$.  Any flat knot relative to $\alpha$ equivalent to
$\alpha_{\epsilon}$ is by definition a $(p, q)$-satellite of $\alpha$.

Poincar\'e \cite{poincare} observed that a $(p, q)$-satellite of a
simple closed curve $\alpha$ has $2p$ intersections with $\alpha$ and
$p(q-1)$ self-intersections. See also Lemma \ref{lemma:poincare}.

\Subsec{Poincar\ee \/{\rm '}\/s rotation number} Let $\gamma(s)$ be an
arc-length parametrization of a closed geodesic of length $L>0$ on $(M,
g)$.  Thus $\gamma(s+L)\equiv \gamma(s)$, and $T=\gamma'(s)$ satisfies
$\nabla_{T}T=0$.  Jacobi fields are solutions of the second order ODE
\begin{equation}
    \frac{d^2 y}{ds^2}+ K(\gamma(s))y(s) = 0,
    \label{eq:jacobi}
\end{equation}
where $K:M\to \R$ is the Gaussian curvature of $(M, g)$.

Let $y:\R\to\R$ be any Jacobi field, and label the zeroes of $y$ in
increasing order
\[
\ldots <s_{-2}<s_{-1}<s_{0}<s_{1}<s_2<\dots 
\]
with $(-1)^{n}y'(s_n)>0$.  Using the Sturm oscillation theorems one
can then show that the limit
\[
\omega(\gamma)= \lim_{n\to\infty} \frac{s_{2n}}{nL}
\]
exists and is independent of the chosen Jacobi field $y$.  We call
this number the Poincar\'e rotation number of the geodesic $\gamma$.
If there is a Jacobi field with only finitely many zeroes then the
oscillation theorems again imply that $y(s)$ has either one or no
zeroes $s\in\R$. In this case we say the rotation number is infinite.

For an alternative definition we observe that if $y(s)$ is a Jacobi
field then $y(s)$ and $y'(s)$ cannot vanish simultaneously. Thus one
can consider
\[
\rho(\gamma) = \lim_{s\to\infty} \frac{L}{2\pi s}\arg\{y(s)+iy'(s)\}.
\]
Again it turns out that this limit exists and is independent of the
particular choice of Jacobi field $y$. Moreover one has
\[
\rho=\frac{1}{\omega}.
\]
We call $\rho$ the inverse rotation number of $\gamma$.  See
\cite{JohnsonMoser} where the much more complicated case of quasi-periodic potentials is treated.  The
inverse rotation number $\rho$ is analogous to the ``amount of rotation'' of a periodic orbit of a twist
map introduced by Mather in \cite{mather-rotation}.

\Subsec{Allowable metrics for a given relative flat knot and the 
  nonresonance 
\pagebreak
condition}\label{sec:intro-resonance} Let
$\Gamma=\{\gamma_1, \dots , \gamma_N\}\subset \Omega$ be a collection
of curves with no mutual or self-tangencies, and denote by
$\cM_{\Gamma}$ the space of $C^{2, \mu}$ Riemannian metrics $g$ on $M$
for which the $\gamma_i\in\Gamma$ are geodesics (thus the metric has
continuous derivatives of second order which are H\"older continuous
of some exponent $\mu\in(0,1)$).  When written out in coordinates one
sees that this condition is quadratic in the components $g_{ij}$ and
$\pt_ig_{jk}$ of the metric and its derivatives. Thus $\cM_{\Gamma}$
is a closed subspace of the space of $C^{2, \mu}$ metrics on $\cM$.

If $\alpha\in\Omega$ is a flat knot rel $\Gamma$ then it may happen
that $\alpha$ is a $(p_1, q_1)$ satellite of, say, $\gamma_1$.  In
this case the rotation number of $\gamma_1$ will affect the number of
closed geodesics of flat knot type $\alpha$ rel $\Gamma$.  To see
this, consider a family of metrics $\{g_{\lambda} \mid
\lambda\in\R\}\subset\cM_{\gamma}$ for which the inverse rotation
number $\rho(\gamma; g_\lambda)$ is less than $p_1/q_1$ for negative
$\lambda$ and more than $p_1/q_1$ for positive $\lambda$. Then, as
$\lambda$ increases from negative to positive, a bifurcation takes
place in which generically two $(p_1, q_1)$ satellites of $\gamma_1$
are created. These bifurcations appear as resonances in the Birkhoff
normal form of the geodesic flow on the unit tangent bundle near the
lift of $\gamma$. This is described by Poincar\'e in \cite[\S6,
p.~261]{poincare}. See also \cite[Appendix 7D,F]{Arnold2}.

In studying the closed geodesics of flat knot type $\alpha$ rel
$\Gamma$ we will therefore exclude those metrics for which a
bifurcation can take place.  To be precise, given $\alpha$ we order
the $\gamma_i$ so that $\alpha$ is a $(p_i, q_i)$ satellite of
$\gamma_i$, if $1\leq i\leq m$, but not a satellite of $\gamma_i$ for
$m<i\leq N$.  We then impose the nonresonance condition
\begin{equation}
    \rho(\gamma_i) \neq  \frac{p_i}{q_i} \text{ for 
    $i\in\{1,\dots ,m\} $}.
    \label{eq:nonresonance}
\end{equation}
The metrics $g\in\MG$ which satisfy this condition can be separated
into $2^m$ distinct classes.  For any subset $I\subset\{1, \dots ,
m\}$ we define $\cM_{\Gamma}(\alpha; I)$ to be the set of all metrics
$g\in \cM_{\Gamma}$ such that the inverse rotation numbers
$\rho(\gamma_1)$, \dots , $\rho(\gamma_{m})$ satisfy
\begin{equation}
    \rho(\gamma_i) <\displaystyle{\frac{p_{i}}{q_i}}   \text{ if } 
    i\in I \text{ and }
    \rho(\gamma_i) >\displaystyle{\frac{p_{i}}{q_i}}   \text{ if } 
    i\not\in I . \label{eq:rotationnumbers}
\end{equation}
For each $I\subset \{1, \dots , m\}$ we define in Section
\ref{sec:theindex} a Conley index $h^I$.  This is done by choosing a
metric $g\in\MG(\alpha;I)$, suitably modifying the set
$\cB_{\alpha}\subset\Omega$ and its exit set for the \CS\ flow,
according to the choice of $I\subset\{1,\dots ,m\}$ and then finally
setting $h^I$ equal to the homotopy type of the modified
$\cB_{\alpha}$ with its exit set collapsed to a point.  Thus the index
we define is the homotopy type of a topological space with a
distinguished point.  We show that the resulting index $h^I$ does not
depend on the choice of metric $g\in\cM(\alpha; I)$, and also that the
index $h^I$ does not change if one replaces $\alpha$ by an equivalent
flat knot rel $\Gamma$.

Using rather standard variational methods we then show in
\S\ref{sec:existence}:
\begin{theorem}\label{thm:main1}
  If $g\in \cM_{\Gamma}(\alpha; I)$ and if the index $h^I$ is
  nontrivial{\rm ,} then the metric $g$ has at least one closed
  geodesic of flat knot type $\alpha$ {\rm rel} $\Gamma$.
\end{theorem}

Using more standard variational arguments one could then improve on
this and show that there are at least $n-1$ closed geodesics of type
$\alpha$ rel $\Gamma$, where $n$ is the Lyusternik-Schnirelman
category of the pointed topological space $h^I$. We do not use this
result here and omit the proof.

Computation of the index $h^I$ for an arbitrary flat knot $\alpha$ rel
$\Gamma$ may be difficult. It is simplified somewhat by the
independence of $h^I$ from the metric $g\in\MG(\alpha;I)$. In addition
we have a long exact sequence which relates the homologies of the
different indices one gets by varying $I$.
\begin{theorem}\label{thm:main2}
  Let $\varnothing\subset J\subset I\subset \{1, \dots , m\}$ with
  $J\neq I$.  Then there is a long exact sequence
    \begin{equation}
        \dots  H_{l+1}(h^I) \stackrel{\partial_*}{\longrightarrow}
        H_{l}(\cA^I_J)\longrightarrow
        H_{l}(h^J)\longrightarrow
        H_{l}(h^I)\stackrel{\partial_*}{\longrightarrow}
        H_{l-1}(\cA^I_J)\dots 
    \end{equation}
    where
    \[
    \cA^I_J = \bigvee_{k\in I\setminus J}\left\{ \frac{S^{1}\times
        S^{2p_k-1}} {S^1\times\{\textrm{pt}\} } \right\}.
    \]
\end{theorem}

This immediately implies
\begin{theorem}\label{thm:main3}
  If $J\subset I$ with $J\neq I$ then $h^I$ and $h^J$ cannot both be
  trivial.
\end{theorem}

One may regard this as a global bifurcation theorem. If for some
choice of rotation numbers $I$ and some choice of metric $g\in
\MG(\alpha;I)$ there are no closed geodesics of type $\alpha$ rel
$\Gamma$, then the index $h^I$ is trivial. By increasing one or more
of the rotation numbers (i.e. increasing $I$ to $J$), or by decreasing
some of the rotation numbers (i.e.\ decreasing $I$ to $J$) the index
$h^I$ becomes nontrivial, and a closed geodesic of type $\alpha$ rel
$\Gamma$ must exist \emph{for any metric} $g\in\MG(\alpha; J)$.

When applied to the case where $M=S^2$ and $\Gamma$ consists of one
simple closed curve $\gamma$ this gives us the following result.

\begin{theorem}\label{thm:main4}
  Let $g$ be a $C^{2,\mu}$ metric on $M$ with a simple closed geodesic
  $\gamma\in\Omega$.  Let $\rho=\rho(\gamma, g)$ be the inverse
  rotation number of $\gamma$.
    
  If $\rho>1$ then for each $\frac{p}{q}\in(1, \rho)$ there is a
  closed geodesic $\gamma_{p/q}$ on $( M, g)$ which is a $(p, q)$
  satellite of $\gamma$.
    
  Similarly{\rm ,} if $\rho<1$ then for each $\frac{p}{q}\in(\rho, 1)$ there
  is a closed geodesic $\gamma_{p/q}$ on $( M, g)$ which is a $(p, q)$
  satellite of $\gamma$.
    
  In both cases the geodesic $\gamma_{p/q}$ intersects the given
  simple closed geodesic $\gamma$ exactly $2p$ times{\rm ,} and
  $\gamma_{p/q}$ intersects itself exactly $p(q-1)$ times.
\end{theorem}

{\it Acknowledgements}. The work in this paper was inspired by a
question of Hofer (Oberwollfach, 1993) who asked me if one could apply
the Floer homology construction to curve shortening, and which results
could be obtained in this way.  This turned out to be a very fruitful question, even though in the end curve
shortening appears to be sufficiently well behaved to use the Conley index instead of Floer's approach.

The paper was finished during my sabattical at the University of
Leiden. It is a pleasure to thank Rob van der Vorst, Bert Peletier and
Sjoerd Verduyn Lunel for their hospitality.

 \vglue12pt\centerline{\bf Contents}
\def\sni#1{\vskip-1pt\noindent{#1}.\phantom{0}}
\sni{1} Introduction
\sni{2} Flat knots
\sni{3} Curve shortening
\sni{4} Curve shortening near a closed geodesic
\sni{5} Loops
\sni{6} Definition of the Conley index of a flat knot
\sni{7} Existence theorems for closed geodesics
\sni{8} Appendices
\vskip-1pt\noindent References 

\section{Flat knots}
\vglue-12pt
\Subsec{The space of immersed curves}\label{sec:charts}
The space of immersed curves $\Om=\mathrm{Imm}(S^1$, $ M)/
\mathrm{Diff}_+\lp S^{1}\rp$ is locally homeomorphic to $C^2(\R/\Z)$.
The homeomorphisms are given by the following charts.  Let
$\gamma\in\Om$ be a given immersed curve.  Choose a $C^{2}$
parametrization $\gamma:\R/\Z\to M$ of this curve and extend it to a
$C^{2}$ local diffeomorphism $\sigma:(\R/\Z)\times (-r, r) \to M$ for
some $r>0$.  Then for any $C^{1}$ small function $u\in C^{2}(\R/\Z)$
the curve
\begin{equation}
        \gamma_{u}(x) = \sigma(x, u(x))
        \label{eq:chart}
\end{equation}
is an immersed $C^2$ curve.  Let $\cU_r= \{u\in C^{2}(\R/\Z) : |u(x)|
< r \}$. For sufficiently small $r>0$ the map $\Phi: u\in\cU_{r}
\mapsto \gamma_{u}\in\Omega$ is a homeomorphism of $\cU_{r}$ onto a
small neighborhood $\Phi(\cU_{r})$ of $\gamma$. The open sets
$\Phi(\cU_{r})$ which one gets by varying the curve $\gamma$ cover
$\Omega$, and hence $\Omega$ is a topological Banach manifold with
model $C^{2}(\R/\Z)$.

A natural choice for the local diffeomorphism $\sigma$ would be
$$\sigma(x, u) = \exp_{\gamma(x)}(u\norm(x))$$ where $\norm$ is a unit
normal vector field for the curve $\gamma$. We avoid this choice of
$\sigma$ since it uses too many derivatives. For $\sigma$ to be $C^2$
one would want the normal to be $C^2$, so the curve would have to be
$C^3$; one would also want the exponential map to be $C^2$, which
requires the Christoffel symbols to have two derivatives, and so the
metric $g$ would have to be $C^3$.

For future reference we observe that if the curve $\gamma$ is
$C^{2,\mu}$ then one can also choose the diffeomorphism $\sigma$ to be
$C^{2,\mu}$.

\Subsec{Covers}
For any $\gamma\in\Om$ and any nonzero integer $q$ we define
$q\cdot\gamma$ to be the $q$-fold cover of $\gamma$, i.e.\ the curve
with parametrization
$$
(q\cdot\gamma)(t) = \gamma(qt),\;\; t\in\R/\Z,$$
where
$\gamma:\R/\Z\to M $ is a parametrization of $\gamma$. Thus
$(-1)\cdot\gamma$ is the curve $\gamma$ with its orientation reversed.

A curve $\gamma\in\Om$ will be called {\it primitive} if it is not a
multiple cover of some other curve, i.e.\ if there are no $q\ge2$ and
$\gamma_0\in\Om$ with $\gamma=q\cdot \gamma_0$.

\Subsec{Flat knots}
Let $\gamma_1$, $\dots $, $\gamma_N$ be a collection of primitive
immersed curves in $ M $.  Define
\begin{align}
  \Delta(\gamma_1, \dots , \gamma_N) & = \lb \gamma\in\Om
  \;\;\left|\;\; \parbox[c]{1.7in}{$\gamma$ has a self-tangency or a
      tangency with one of the $\gamma_i$}
  \right.\rb \\
  \intertext{and} \Delta &= \lb \gamma\in\Om \mid\hbox{$\gamma$ has a
    self-tangency}\rb.
\end{align}
Then $\Delta$ and $\Delta(\gamma_1, \dots , \gamma_N)$ are closed
subsets of $\Om$, and their complements $\Omega\setminus \Delta$ and
$\Omega\setminus \Delta(\gamma_1, \dots , \gamma_N)$ consist of flat
knots, and flat knots relative to $(\gamma_1, \dots , \gamma_N)$,
respectively.  Two such flat knots are equivalent if and only if they
lie in the same component of $\Omega\setminus\Delta$ or
$\Omega\setminus\Delta(\gamma_1, \dots , \gamma_N)$.

\Subsec{Flat knots as knots in the projective tangent bundle}
Let $\ptm$ be the {\it projective tangent bundle} of $ M $, i.e.\ $\P
T M$ is the bundle obtained from the unit tangent bundle
\[
\utb = \{(p, v)\in T( M ) \mid g(v, v)=1\}
\]
by identification of  all antipodal vectors $(x, v)$ and $(x, -v)$.  The
projective tangent bundle is a contact manifold.  If we denote the
bundle projection by $\pi:\ptm\to M$, then the contact plane $L_{(x,
  \pm v)}\subset T(\ptm)$ at a point $(x, \pm v)\in\ptm$ consists of
those vectors $\xi\in T(\ptm)$ for which $d\pi(\xi)$ is a multiple of
$v$.  Each contact plane $L_{(x, \pm v)}$ contains a nonzero vector
$\vartheta$ with $d\pi(\vartheta)=0$ ($\vartheta$ corresponds to
infinitesimal rotation of the unit vector $\pm v$ in the tangent space
$T_{x} M$, while   the base point $x$ remains fixed).

Any $\gamma\in\Om$ defines a $C^{1}$ immersed curve $\hg$ in the
projective tangent bundle $\ptm$ with parametrization $\hg(s)=
(\gamma(s), \pm\gamma'(s))$, where $\gamma(s)$ is an arc length
parametrization of $\gamma$.  We call $\hg$ the lift of $\gamma$.

An immersed curve $\tilde{\gamma}$ in $\ptm$ is the lift of some
$\gamma\in\Omega$ if and only if $\tilde{\gamma}$ is everywhere
tangent to the contact planes, and nowhere tangent to the special
direction $\vartheta$ in the contact planes.

Self-tangencies of $\gamma\in\Omega$ correspond to self-intersections of
its lift $\hat{\gamma}\subset\ptm$.  Thus an immersed curve
$\gamma\in\Omega$ is a flat knot exactly when its lift $\hat{\gamma}$
is a knot in the three manifold $\ptm$.  If two curves $\gamma_{1},
\gamma_{2}\in\Omega$ define equivalent flat knots then one can be
deformed into the other through flat knots.  By lifting the
deformation we see that $\hat{\gamma}_{1}$ and $\hat{\gamma}_{2}$ are
equivalent knots in $\ptm$.

\Subsec{Intersections} If $\alpha \in \Omega \setminus 
\Delta(\gamma_{1}$, $\dots $, $\gamma_{n})$ then $\alpha$ is
transverse to each of the $\gamma_{i}$.  Hence the number of
intersections in $\alpha\cap\gamma_{i}$ is well defined.  This only
depends on the flat knot type of $\alpha$ relative to $\gamma_{1}$,
$\dots $, $\gamma_{n}$.

\begin{figure}[ht]
  \centering \includegraphics{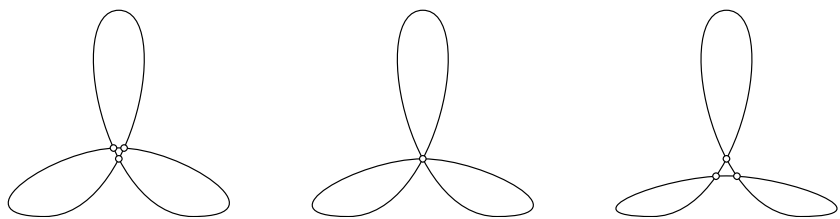}
    \caption{Equivalent flat knots with 3 self-intersections.}
    \label{fig:threeintersections} \vglue-12pt
\end{figure}

If $\alpha\in\Omega\setminus\Delta$ then $\alpha$ only has transverse
self-intersections, so their number is well defined by
$\#\alpha\cap\alpha = \#\{0\leq x<x'<1\mid \alpha(x)=\alpha(x')\}$.
From a drawing of $\alpha$ they are easily counted.  An
$\alpha\in\Omega\setminus\Delta$ can only have double points, triple
points, etc.\  (see Figure \ref{fig:threeintersections}).  If $\alpha$
only has double points (a generic property) then their number is the
number of self-intersections.  Otherwise one must count the number of
geometric self-intersections where a $k$-tuple point counts for
$\binom{k}{2}$ self-intersections.  Again this number only depends on
the flat knot type of $\alpha\in\OMD$.

\Subsec{Nontransverse crossings of curves}
\label{sec:nontransverse} If $\gamma_{1}, \gamma_{2}\in\Omega$ are
not necessarily transverse then we define the number of crossings of
$\gamma_{1}$ and $\gamma_{2}$ to be
\begin{equation}
        \cross(\gamma_{1}, \gamma_{2})
        = \sup_{\gamma_{i}\in\cU_{i}}
                \inf\left\{\#(\gamma_{1}' \cap \gamma_{2}')\; \left|
                \begin{array}{c}
                        \gamma_{1}'\in\cU_{1}, \: \gamma_{2}'\in\cU_{2}  \\
                        \gamma_{1}'\tv \gamma_{2}'
                \end{array}
                \right.\right\}     
        \label{eq:crossingsdef}
\end{equation}
where the supremum is taken over all pairs of open neighborhoods
$\cU_{i}\subset \Omega$ of $\gamma_{i}$.  Thus $\cross(\gamma_{1},
\gamma_{2})$ is the smallest number of intersections $\gamma_{1}$ and
$\gamma_{2}$ can have if one perturbs them slightly to become
transverse.

The number of self-crossings $\cross(\gamma, \gamma)$ is defined in a
similar way.

Clearly $\cross(\gamma_{1}, \gamma_{2})$ is a lower semicontinuous
function on $\Omega \times \Omega$.

\Subsec{Satellites}
\label{sec:BirkhoffOrbits}
We first describe the local model of a satellite of a primitive flat
knot $\gamma\in\Omega\setminus\Delta$ and then transplant the local
model to primitive flat knots on any surface.

Let $q\geq 1$ be an integer, and let $u\in C^{2}(\R/q\Z)$ be a
function for which
\begin{equation}\hbox{all zeroes of $u$ are simple}\label{eq:uzerosimple} 
\end{equation}
and
\begin{equation}
          \hbox{all zeroes of $v_{k}(x)\isdef u(x)-u(x-k)$ are simple for $k=1, 2, \cdots, q-1$.} 
                \label{eq:uselftransverse}
\end{equation}
Consider the curve $\alpha_{u}$ in the cylinder
$\Gamma=(\R/\Z)\times\R$, parametrized by
\begin{equation}
        \alpha_{u}:\R/q\Z \to \Gamma, \quad
        \alpha_{u}(x)= (x, u(x)).
        \label{eq:satellitegraph}
\end{equation}
The conditions \eqref{eq:uzerosimple} and \eqref{eq:uselftransverse}
imply that $\alpha_{u}$ is a flat knot relative to $\alpha_{0}$, where
$\alpha_{0}=(\R/\Z)\times\{0\}$ is the zero section (i.e., the curve
corresponding to $u(x)\equiv 0$).

Now consider a primitive flat knot $\gamma\in\Omega\setminus\Delta$.
Denote by $\gamma:\R/\Z\to M $ any parametrization, and choose a local
diffeomorphism $\sigma:\R/\Z\times (-r, r)\to M$ with $\gamma(x) =
\sigma(x, 0)$.  As in \S\ref{sec:charts} we then identify any curve
$\gamma_{u}$ which is $C^1$ close to $\gamma$ with a function $u\in
C^2(\R/\Z)$ via \eqref{eq:chart}.

If $u\in C^{2}(\R/q\Z)$ then the curve defined by
\begin{equation}
    \alpha_{\varepsilon,u}(x) = \sigma(x, \varepsilon u(x))
    \label{eq:satellite-def}
\end{equation}
is a flat knot relative to $\gamma$.  For given $u\in C^{2}(\R/q\Z)$
and small enough $\eps>0$ the $\alpha_{{\eps, u}}$ all define the same
relative flat knot.

By definition, a curve $\alpha\in\Omega\setminus \Delta(\gamma)$ is a
satellite of $\gamma\in\Omega\setminus\Delta$ if for some $u\in
C^{2}(\R/q\Z)$ it is isotopic relative to $\gamma$ to all
$\alpha_{\eps, u}$ with $\eps>0$ sufficiently small.

To complete this definition we should specify the orientation of the
satellite $\alpha_{\eps, u}$.  One can give $\alpha_{\eps, u}$ as
defined in (\ref{eq:satellite-def}) the same orientation as its base
curve $\gamma$, or the opposite orientation.  We will call both curves
satellites of $\gamma$.  In general the satellites $\alpha_{\eps, u}$
and $-\alpha_{\eps, u}$ can define different flat knots relative to
$\gamma$ or they can belong to the same relative flat knot class.

\demo{Example}\label{sssec:cb11example} Let $\gamma$ be the
equator on the standard two sphere $ M =S^{2}$.  Then any other great
circle is a satellite of $\gamma$.  Moreover, all these great circles
with either orientation define the same flat knot relative to the
equator.  For example, if $\alpha$ is a great circle in a plane through the
$x$-axis which makes an angle $\pfi\ll\pi/2$ with the $xy$-plane, then
one can reverse its orientation by first rotating it through
$\pi-2\pfi$ around the $x$-axis, and then rotating it through $\pi$
around the $z$-axis.  Throughout this motion the curve remains
transverse to the equator, so that $\alpha$ and $-\alpha$ indeed
belong to the same component of $\Delta\setminus \Omega(\gamma)$.
Below we will show that this example is exceptional.

As defined in the introduction, one obtains $(p,q)$ satellites by
setting
\begin{equation}
        u(x) = \sin(2\pi\frac{p}{q}x).
\end{equation}

Let $p\neq 0$, and let $\alpha$ be the $(p, q)$ satellite of $\gamma$
given by $u(x) = \epsilon\sin(2\pi\frac{p}{q}x)$.  Then we can
translate $\alpha$ along the base curve $\gamma$; i.e.\ we can
consider the $(p, q)$ satellites $\alpha_{\tau}$ given by $u_{\tau}(x)
= \epsilon \sin(2\pi\frac{p}{q}(x-\tau))$.  By translating from
$\tau=0$ to $\tau=\frac{q}{2p}$ one finds an isotopy from $\alpha$ to
the curve $\bar\alpha$ given by $\bar u(x) = \sin(2\pi\frac{-p}{q}x)$.
Hence one can turn any $(p, q)$ satellite into a $(-p, q)$ satellite,
and we may therefore always assume that $p$ is nonnegative.

We will denote the set of $(p, q)$-satellites of $\gamma\in\Omega$ by
$\cB_{p, q}(\gamma)$, always assuming that $p\geq 0$ and $q\geq 1$.

More precisely we will let $\cB^{+}_{p, q}(\gamma)$ be the set of $(p,
q)$-satellites of $\gamma$ which have the same orientation as
$\gamma$, and we let $\cB^{-}_{p, q}(\gamma)$ be those $(p, q)$
satellites with opposite orientation. With this notation we always
have
\[
\cB_{p, q}(\zeta) = \cB_{p, q}^{+}(\zeta) \cup \cB_{p, q}^{-}(\zeta).
\]

It is not \emph{a priori} clear that all these classes are disjoint,
but by counting the number of self-intersections of $(p, q)$ satellites
one can at least see that there are infinitely many disjoint $\cB_{p,
  q}$'s.
\begin{lemma}\label{lemma:poincare}
  Let $\gamma\in\Omega\setminus\Delta$ be a flat knot with $m$ self-intersections.  Then any
$\alpha\in\cB_{p, q}(\gamma)$ has exactly
  $2p+2mq$ intersections with $\zeta${\rm ,} and $p(q-1)+mq^{2}$
  self\/{\rm -}\/intersections.
\end{lemma}

This was   observed by Poincar\'e   \cite{poincare}. We include
a proof for completeness' sake.

\Proof Intersections of $\alpha$ and $\gamma$ are of two types.  Each zero
  of $u(x)$ corresponds to an intersection of $\alpha$ and $\gamma$.
  At each self-intersection of $\gamma$ the two intersecting strands
  of $\gamma$ are accompanied by $2q$ strands of $\alpha$ which
  intersect $\gamma$ in $2q$ points.  Since $u(x)$ has $2p$ zeroes and
  $\gamma$ has $m$ self-intersections we get $2mq+2p$ intersections of
  $\alpha$ and~$\gamma$.
  
  To count self-intersections one must count the intersections of the
  graph of $u(x)=\sin(2\pi\frac{p}{q}x)$ wrapped up on the cylinder
  $\Gamma=(\R/\Z)\times \R$, i.e.\ the intersections of the graphs of
  $u_{k}(x)= u(x-2k)$ ($k=0 , 1,  \dots  ,  q-1$) with $0 \leq x <
  2\pi$.  After some work one finds that these are arranged in $q-1$
  horizontal rows, each of which contains $p$ intersections.
  
  At each self-intersection of $\gamma$ two strands of $\gamma$ cross.
  If $\eps$ is small enough then $\alpha_{\eps, u}$ is locally almost
  parallel to $\gamma$, so that any pair of crossing strands of
  $\gamma$ is accompanied by a pair of $q$ nearly parallel strands of
  $\alpha$ which cross each other.  This way we get $q^{2}$ extra
  self-crossings of $\alpha$ and $2q$ extra crossings of $\gamma$ with
  $\alpha$ per self-crossing of $\gamma$.
\hfill\qed

\begin{lemma}\label{lem:pqisrs}
  If $\cB_{p, q}(\gamma)\cap\cB_{r, s}(\gamma)\neq\varnothing$ then
  $p=r$ and $q=s$.
\end{lemma}

\Proof   If $\alpha\in\cB_{p, q}$ has $2k$ intersections with $\gamma$
and $l$ self-intersections then
\[
p(q-1)+mq^{2} = l, \qquad p+mq=k.
\]
Substitute $p=k-mq$ in the first equation to get
\[
l=(k+m)q-k = mq+(q-1)k
\]
from which one finds $q=\frac{k+l}{k+m}$. In particular, the numbers
$k$, $l$ and $m$ determine $p$ and $q$.  \Endproof\vskip4pt

The proof also shows that most satellites are not $(p,q)$-satellites
for any $(p, q)$.  Indeed, given $\alpha\in\cB_{p, q}(\gamma)$ one can
modify it near one of its crossings with $\gamma$ so as to increase
the number $k$ of intersections with $\gamma$ arbitrarily without
changing the number of self-intersections $l$, or $m$.  Unless both
$l=0$ and $m=0$, then for large enough $k$ the fraction
$\frac{k+l}{k+m}$ will not be an integer, so the modified curve can no
longer be a $(p,q)$ satellite.  If both $l=m=0$ then both $\gamma$ and
its satellite $\alpha$ must be simple curves.

\Subsec{$(p,q)$ satellites along a simple closed curve on
  $S^{2}$}\label{sec:pq-satellites} In this section we consider the
case in which $ M =S^{2}$ and $\zeta\in\Omega$ is a simple closed
curve. We will show that for all $(p, q)$ except $p=q=1$ the classes
$\cB^{\pm}_{p,q}(\zeta)$ are different.

After applying a diffeomorphism we may assume that $ M $ is the unit
sphere in $\R^{3}$ and that $\zeta$ is the equator, given by $z=0$.

To study curves in $\Omega\setminus\Delta(\zeta)$ it is useful to
recall that one can identify the unit tangent bundle $T^{1}(S^{2})$ of
the 2-sphere with the group SO$(3, \R)$.  Indeed, by definition,
\[
T^{1}(S^{2})=\{(\vec x, \vec \xi) \in\R^{3}\times\R^{3} \mid |\vec
x|=|\vec \xi|=1, \vec x\perp\vec \xi \}
\]
so that any unit tangent vector $(\vec x, \vec \xi)\in T^{1}(S^{2})$
determines the first two columns of an orthogonal matrix. The third
column of this matrix is the cross product $\vec x\times \vec \xi$.
The map
\[
(\vec x, \vec \xi)\in T^{1}(S^{2}) \mapsto (\vec x, \vec \xi, \vec
x\times\vec \xi) \in\textrm{SO}(3, \R)
\]
is a diffeomorphism, and from here on we will simply identify
$T^{1}(S^{2})$ and SO$(3, \R)$.

\begin{figure}[ht]
 \centering \includegraphics{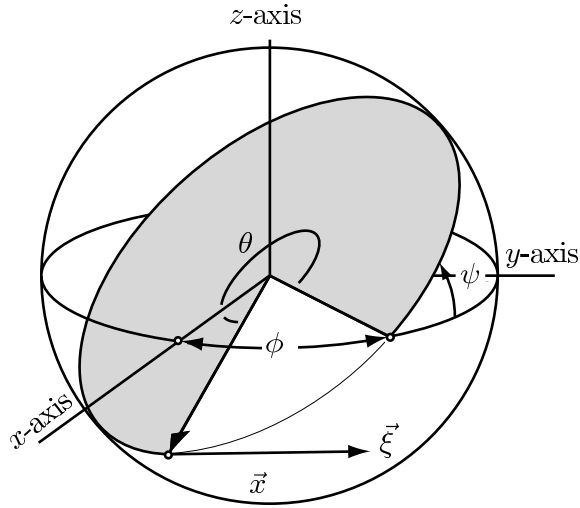}
  \caption{Euler angles  $\phi$, $\psi$ and $\theta$.}
  \label{figure:eulerangles} \vglue-12pt
\end{figure}

Let $\cU \subset T^{1}(S^{2})$ be the complement of the set of tangent
vectors to $\zeta$ and $-\zeta$.  One can describe $\cU$ very
conveniently using ``Euler Angles''.  For the definition of these
angles we refer to Figure \ref{figure:eulerangles}.  Any unit tangent
vector $(\vec x, \vec \xi)$ defines an oriented great circle,
parametrized by
\[
X(t) = (\cos t)\vec x+(\sin t)\vec \xi.
\]
Unless $(\vec x, \vec \xi)$ is a tangent vector of the equator
$\pm\zeta$, the great circle through $(\vec x, \vec \xi)$ intersects
the equator in two points.  In one of these intersections the great
circle crosses the equator from south to north.  Let $\theta$ be the
angle from the upward intersection to $x$, so that $X(-\theta)$ is the
upward intersection point.  We define $\psi$ to be the angle between
the plane through the great circle $\{X(t) \mid t\in\R\}$ and the
$xy$-plane (so that $0<\psi<\pi$).  Finally we let $\phi$ be the angle
along the equator $\zeta$ from the $x$-axis to the upward intersection
point $X(-\theta)$.

If we denote the matrix corresponding to a rotation by an angle
$\alpha$ around the $x$ axis by $\rR_{x}(\alpha)$, etc.\ then the
relation between the Euler angles $(\theta, \pfi, \psi)$ and the unit
tangent vector $(x, \xi)$ they represent is given by
\begin{equation}
        (\vec x, \vec \xi, \vec x\times\vec \xi) 
        = \rR_{z}(\phi)\cdot  \rR_{x}(\psi)\cdot \rR_{z}(\theta).
        \label{eq:Eulerangles}
\end{equation}
The map $(\vec x, \vec \xi) \mapsto (\theta, \psi, \phi)$ is a
diffeomorphism between $\cU$ and $(\R/2\pi\Z)\times (0, \pi) \times
(\R/2\pi\Z) \cong \T^{2}\times\R$.

Given this identification we can now define two numerical invariants
of flat knots $\alpha$ relative to the equator $\zeta$.  By the lift of  a
unit speed parametrization, any flat knot
$\alpha\in\Omega\setminus\Delta(\zeta)$ defines a closed curve
$\hat\alpha:S^{1}\to\cU$.  The numerical invariants are then the
increments of the Euler angles $\theta$ and $\phi$ along $\hat\alpha$,
which we will denote by $\Delta\theta(\alpha)$ and
$\Delta\phi(\alpha)$, respectively.  Both are integral multiples of
$2\pi$.
\begin{lemma}
  If $\alpha$ is a satellite of $\zeta$ given by
  \eqref{eq:satellite-def} then
  \begin{subequations}
    \begin{eqnarray}
      \pm\Delta\theta+\Delta\phi & = & 2q\pi,\label{eq:2qformula} \\
      \Delta\theta & = &  2p\pi
    \end{eqnarray}
  \end{subequations}
  where $2p$ is the number of zeroes of $u\in
  C^{2}(\R/2q\pi\Z)$.  In the first equation one must take the
  \/{\rm ``}\/$+$ sign\/{\rm ''}\/ if $\alpha$ has the same orientation as $\zeta${\rm ,}
  and the \/{\rm ``}\/$-$ sign\/{\rm ''}\/ otherwise.
\end{lemma}

Note that the number of zeroes of $u\in C^{2}(\R/2\pi\Z)$ must always
be even (assuming they are all simple zeroes, of course).
\Proof We project the sphere onto the cylinder $x^{2}+y^{2}=1$ and write
  $z$ and $\vth$ for the usual coordinates on this cylinder. We assume
  that $\alpha$ projects to the graph of $z=u(\vth)$ on the cylinder,
  and that $u$ is a 2$q\pi$ periodic function with simple zeroes only,
  and for which $|u(\vth)|+|u'(\vth)|$ is uniformly small.
  \begin{figure}[ht]
     \centering \includegraphics{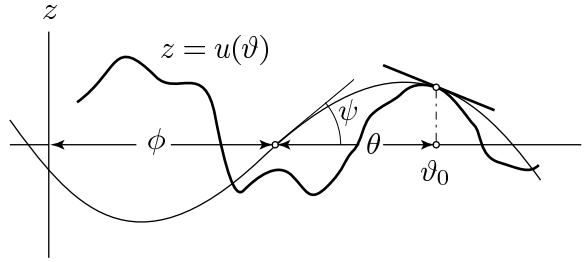}
    \caption{A great circle projected onto the cylinder.}
    \label{fig:Euler-angles-near-equator} \vglue-12pt
  \end{figure}
  Let $\alpha$ have the same orientation as the equator (from west to
  east).  We compute the Euler angles corresponding to the unit
  tangent vector to $\alpha$ at the point which projects to
  $(\vartheta_0, u(\vartheta_0))$ on the cylinder.  In Figure
  \ref{fig:Euler-angles-near-equator} we have sketched the great
  circle which passes through $(\vth_0, u(\vth_0))$ with slope
  $u'(\vartheta_0)$ as it appears in\break ($\vartheta$, $z$) coordinates on
  the cylinder.  Since great circles are intersections of planes
  through the origin with the sphere, they project to intersections of
  such planes with the cylinder, and are therefore graphs of
  $z=\psi\sin(\vartheta-\phi)$.
    
  From Figure \ref{fig:Euler-angles-near-equator} one finds
  \begin{equation} \label{eq:th+psi}
    \theta+\phi  =  \vth_{0} ,\qquad
    u(\vth_0)  =  \psi\sin\theta ,\qquad
    u'(\vth_0)  =  \psi\cos\theta, 
  \end{equation}
  so that
  \begin{equation}
    \theta = \arg(u'(\vth_{0})+iu(\vth_{0})).
    \label{eq:psi-is-arg}
  \end{equation}
  From \eqref{eq:th+psi} we see that $\theta+\phi$ increases by
  $2q\pi$ along the curve $\alpha$. To compute $\Delta\phi$ we use
  (\ref{eq:psi-is-arg}) and count the number of times the curve
  $u'(\vth_{0})+iu(\vth_{0})$ in the complex plane crosses the
  positive real axis. Every such crossing corresponds to a zero of
  $u$ with positive derivative, and hence there are $\frac{2p}{2}=p$
  of them. We conclude that $\Delta\theta=p\times 2\pi$, as claimed.
  
  Similar arguments also allow one to find $\Delta\phi$ and
  $\Delta\theta$ if one gives $\alpha$ the orientation opposite to
  that of the equator.
\Endproof\vskip4pt 

We have observed that $\cB_{1,1}^{+}(\zeta)$ and
$\cB_{1,1}^{-}(\zeta)$ coincide.  If $p/q$ is any fraction in lowest
terms then $\cB_{p,q}^{+}(\zeta) = \cB_{p,q}^{-}(\zeta)$ combined with
\eqref{eq:2qformula} implies $\Delta\theta=0$, and hence $p=q$. Since
$\mathrm{gcd}(p,q)=1$ we conclude

\begin{lemma} If $\zeta$ is a simple closed curve on $S^{2}${\rm ,} and  
  $\cB_{p,q}^{+}(\zeta) = \cB_{p,q}^{-}(\zeta)$ then $p=q=1$.
\end{lemma}

\section{Curve shortening}
\vglue-12pt
\Subsec{The gradient flow  of the length
functional}\label{sec:gradflow}
Let $g$ be a $C^{2,\mu}$ metric on the surface $M$.  Then for any
$C^1$ initial immersed curve $\gamma_0$ a maximal classical solution
to \CS\ exists on a time interval $0\leq t<T(\gamma_0)$.  We denote
this solution by $\{\gamma_t : 0\leq t<T(\gamma_0)\}$.  The solution
depends continuously on the initial data $\gamma_0\in\Om$, so that
\CS\ generates a continuous local~semiflow
\begin{align*}
  &\Phi : \cD\to\Om, \qquad \Phi^{t}(\gamma_0) \isdef \gamma_t, \\
  &\cD = \lb (\gamma, t)\in\Om\times[0, \infty) \mid 0\le
  t<T(\gamma)\rb.
\end{align*}
One can show that if $T(\gamma_{0})<\infty$ then the geodesic
curvature $\kappa_{\gamma_{t}}$ of $\gamma_{t}$ ``blows-up'' as
$t\nearrow T(\gamma_{0})$, i.e.
\[
\lim_{t\nearrow T(\gamma_{0})} \sup_{\gamma_{t}} |\kappa_{\gamma_{t}}|
= \infty.
\]
Since the geodesic curvature itself satisfies a parabolic equation
\begin{equation}
                \frac{\pt\kappa_{\gamma}}{\pt t} = 
                        \frac{\pt^{2}\kappa_{\gamma}}{\pt s^{2}}
                        + \lp K\circ{\gamma} + \kappa_{\gamma}^{2} \rp \kappa_{\gamma}
        \label{eq:kevolve}
\end{equation}
($K\circ{\gamma}$ is the Gauss curvature of the surface evaluated
along the curve) the maximum principle implies that one has the
following lower estimate for the lifetime of any solution.  If
$T(\gamma_{0})\leq 1$ then
\begin{equation}
        T(\gamma_{0}) \geq \frac{C}{\sqrt{\smash[b]{\sup_{\gamma_{t}}|\kappa|}}}
        \label{eq:lifetimebound}
\end{equation}
where $C$ is some constant  depending on $\sup_{M} |K|$ only.  See
\cite{GH1986} or \cite{Ang1}.

The \CS\ flow on $\Omega$ provides a gradient flow for the length
functional. Indeed, one has
\begin{equation}
        \frac{dL(\gamma_{t})}{dt} = -\int_{\gamma_{t}} 
        \lp\kappa_{\gamma_{t}}\rp^{2} ds
        \label{eq:Lderiv}
\end{equation}
where $ds$ represents arclength along $\gamma_{t}$. Thus solutions of
\CS\ do indeed always become shorter, unless $\gamma_{t}$ is a
geodesic, in which case the solution $\gamma_{t}\equiv\gamma_{0}$ is
time independent. From the above description of $T(\gamma_{0})$ one
easily derives the following (see \cite{Gray1987}, \cite{Gray1989}, also
\cite{Ang1}, \cite{Ang2}).

\begin{lemma}\label{lemma:omegalimit}
  If $T(\gamma_{0})=\infty$ then
        \[
        \lim_{t\to\infty}\sup_{\gamma_{t}} |\kappa_{\gamma_{t}}| =0.
        \]
        Moreover{\rm ,} any sequence $t_{i}\nearrow\infty$ has a subsequence
        $t_{i}'$ for which $\gamma_{t_{i}'}$ converges to some
        geodesic of $(M, g)$.
\end{lemma}

In other words, orbits of the \CS\ flow $\Phi$ which exist for all
$t\geq 0$ have (compact) omega-limit sets in the sense of dynamical
systems.  Such $\omega$-limit sets,
\[
\omega(\gamma_{0})\isdef \left\{\gamma_{*}\in\Omega \mid \exists
  t_{i}\uparrow\infty : \gamma_{t_{i}}\to\gamma_{*} \right\}
\]
are of course connected, and if the geodesics of $(M, g)$ are isolated
then any orbit of \CS\ either becomes singular or else converges to
one 
\break
geodesic.

The same is true for ``ancient orbits,'' i.e.\ orbits $\{\gamma_{t}\}$
which are defined for all $t\leq 0$ and for which $\sup_{t\leq
  0}L(\gamma_t)<\infty$. For such orbits one can define the $\alpha$
limit set
\[
\alpha(\gamma_{0})\isdef \left\{\gamma_{*}\in\Omega \mid \exists
  t_{i}\searrow-\infty : \gamma_{t_{i}}\to\gamma_{*} \right\},
\]
and this set consists of closed geodesics.

\Subsec{Parabolic estimates} Since \CS\ is a nonlinear heat 
equation solutions are generally smoother than their initial data.
This provides a compactness property which we will use later to
construct the Conley-index.  There are various well-known ways of
deriving the smoothing property of nonlinear heat equations.  Here we
show which estimate one can easily obtain assuming only that the
metric $g$ is $C^2$.

\begin{lemma}\label{lemma:parabolicestimate}
  If $\{\gamma_t \mid 0\leq t\leq t_0\}$ is a solution of \CS\ whose
  curvature is bounded by $|\kappa|\leq A$ at all times{\rm ,} then
\begin{equation}
        \int_{\gamma_t} \kappa_s^2 ds \leq \frac{C}{t}
        \label{eq:ks2estimate}
\end{equation}
where the constant $C$ only depends on $A${\rm ,} $t_0${\rm ,} the length $L$ of
$\gamma(0)$ and $\sup_{M}|K|$.
\end{lemma}

By adding a Nash-Moser iteration to the following arguments one could
improve the estimate \eqref{eq:ks2estimate} to an $L^\infty$ estimate
for $\kappa_s$ of the form $|\kappa_s|\leq C/\sqrt{t}$. However,
\eqref{eq:ks2estimate} will be good enough for us in this paper.

\Proof Let $\gamma:\R/\Z\times [0, T)\to M$ be a normal parametrization of
  a solution of \CS, i.e.\ one with $\pt_{t}\gamma\perp\pt_{s}
  \gamma$. Then the curvature $\kappa$ satisfies \eqref{eq:kevolve},
  and using the commutation relation $[\pt_t, \pt_s] = \kappa^2\pt_s$
  one obtains
    \begin{equation}
        \frac{\pt\kappa_{s}}{\pt t} = 
            \frac{\pt^2 \kappa_s}{\pt s^2}
            + \frac{\pt}{\pt s}\bigl((K\circ\gamma) \kappa +\kappa^3\bigr).
        \label{eq:ksevolve}
    \end{equation}
    The arclength ${d}s$ on the curve evolves by $\frac{\pt}{\pt
      t}{d}s = -\kappa^2{d}s$. Therefore we have
    \begin{align}
      \frac{d}{dt}\int_{\gamma_t} (\kappa_s)^2 ds &=
      \int_{\gamma_t}\lp 2\kappa_s\kappa_{st} - \kappa^2\kappa_s^2\rp
      ds
    \label{eq:ks2growth}  \\
      &= \int_{\gamma_t}\lp -2(\kappa_{ss})^2 + 5 \kappa^2\kappa_s^2 -
      2\kappa (K\circ\gamma) \kappa_{ss} \rp ds
      \nonumber\\
      &\leq C + C\int_{\gamma_t}\kappa_s^2ds -
      \int_{\gamma_t}(\kappa_{ss})^2 ds
    \nonumber    
    \end{align}
    where the constant $C$ only depends on $A$, $L$ and $\sup_{M}|K|$.
    
    By expanding $\kappa(\cdot, t)$ in a Fourier series in $s$ one
    finds that
    \[
    \lp\int_{\gamma_t} \kappa_s^2 ds\rp^2 \leq {\int_{\gamma_t}
      \kappa^2 ds}\; {\int_{\gamma_t} \kappa_{ss}^2 ds},
    \]
    which implies
    \[
    \int_{\gamma_t} \kappa_{ss}^2 ds \geq \frac{1}{C}
    \lp{\int_{\gamma_t} \kappa_{s}^2 ds}\rp^{2}
    \]
    where the constant $C$ only depends on $A=\sup|\kappa|$ and $L$.
    Combined with \eqref{eq:ks2growth} this leads to a differential
    inequality for $\int \kappa_s^{2}ds$,
    \[
    \frac{d}{dt}\int_{\gamma_t} (\kappa_s)^2 ds \leq C +
    C\int_{\gamma_t}\kappa_s^2ds - \frac{1}{C}
    \lp\smash[b]{\int_{\gamma_t} \kappa_{s}^2 ds}\rp^{2}.
    \]
    Integration of this inequality gives \eqref{eq:ks2estimate}.
\Endproof\vskip4pt

This lemma implies that for solutions with bounded curvature the
curvature becomes H\"older continuous with exponent $1/2$, since
\begin{align}
  |\kappa(P, t) - \kappa(Q, t)| &\leq \int_P^Q |\kappa_s|ds  \label{eq:kHolderestimate}\\
  & \leq \lp{\int_P^Q \kappa_s^2ds}\rp^{1/2}\; {\dist_{\gamma_t}(P,
    Q)}^{1/2}
  &\hbox{(Cauchy)}\nonumber\\
  &\leq \frac{C(L, A, \sup|K|)}{t}{\dist_{\gamma_t}(P, Q)}^{1/2},
       \nonumber
\end{align}
$\dist_{\gamma_t}(P, Q)^{1/2}$ being the distance from $P$ to $Q$
along the curve $\gamma_t$.

\Subsec{The nature of singularities in \CS}%
\label{sec:singularities} 
Consider a solution $\{\gamma(t) : 0\leq t<T\}$ of \CS\ with
$T=T(\gamma_0)<\infty$. Then, as $t\nearrow T$, the curve $\gamma_t$
converges to a piecewise smooth curve $\gamma_T$ which has finitely
many singular points $P_1, \dots , P_m$; i.e.\ $\gamma_T$ is the union
of finitely many immersed arcs whose endpoints belong to $\lb P_1,
\dots , P_m\rb$.

Either $\gamma_t$ shrinks to a point (in which case $m=1$, and
$\gamma_T$ consists only of the point $P_1$), or else any
neighborhood $\cU\subset M^2$ of any of the $P_i$ will contain a self-intersecting arc of $\gamma_t$ for
$t$ sufficiently close to $T$.  In other words, $\gamma_{t}\cap\cU$ is the union of a finite number of
arcs, at least one of which has a self-intersection (a parametrization
$x\in\R/\Z \mapsto \gamma_{t}(x)$ of the curve will enter $\cU$ and
self-intersect before leaving the neighborhood).

This description of the singularities which a solution of \CS\ may
develop follows from work of Grayson \cite{Gray1987}, \cite{Gray1989}; see
also \cite{Ang1}, \cite{Ang2}, \cite{Oaks} for a similar result applicable to more
general flows; an alternative proof of the above result can now be
given using the Hamilton-Huisken distinction between ``type 1 and type
2'' singularities (see \cite{ElEscorial} for a short survey), where we apply 
a monotonicity formula in the type 1 case, and either Hamilton's
\cite{Ha:isoper} or Huisken's isoperimetric ratios \cite{Hu:isoper} in
the type 2 case.

\Subsec{Intersections and Sturm\/{\rm '}\/s theorem}
We recall Sturm's theorem \cite{Sturm} which states that if $u(x, t)$
is a classical solution of a linear parabolic equation
\[
\frac{\pt u}{\pt t} = a(x, t) u_{xx} + b(x, t)u_{x}+c(x, t)u
\]
on a rectangular domain $[x_{0}, x_{1}]\times[t_{0}, t_{1}]$, with
boundary conditions
\[
u(x_{0}, t) \neq 0,\; u(x_{1}, t)\neq 0,\;\textrm{ for } t_{0}\leq
t\leq t_{1},
\]
then the number of zeroes of $u(\cdot, t)$
\[
z(u; t)\isdef \#\{ x\in [x_{0}, x_{1}] \mid u(x, t)=0\}
\]
is finite for any $t>t_{0}$, and does not increase as $t$ increases.
Moreover, at any moment $t_{*}$ at which $u(\cdot, t_{*})$ has a
multiple zero, $z(u, t)$ drops.  This theorem goes back to Sturm
\cite{Sturm} who gave a rigorous proof assuming the solutions and
coefficients are analytic functions, which has been rediscovered and
reproved under weaker hypotheses many times since then. See
\cite{Nickel}, \cite{Matano}, \cite{Ang-sturm}.

In \cite{Ang-sagarro} we argue that Sturm's theorem may be considered
as a ``degenerate version'' of the well-known principle that the local
mapping degree of an analytic function $f:\C\to\C$ near any of its
zeroes is always positive (so that one can count zeroes of $f$ by
computing winding numbers, etc.).

Using Sturm's theorem we proved the following in \cite{Ang1}, \cite{Ang2}.

\begin{lemma}\label{lemma:intersections}
  Any smooth solution $\lb \gamma_t\mid 0<t<T\rb$ of \CS\ which is not
  a multiple cover of another solution{\rm ,} always has finitely many self-intersections{\rm ,} all of which are transverse{\rm ,} except at a discrete set
  of times $\lb t_j\rb\subset(0, T)$. At each time $t_j$ the number of
  self-intersections of $\gamma_t$ decreases.
\end{lemma}

A similar statement applies to intersections of two different
solutions: if $\{\gamma_t^1\mid $ $0<t<T\}$ and $\{\gamma_t^2\mid $ $
0<t<T\}$ are solutions of \CS\ then they are transverse to each other,
except at a discrete set of times $\lb t_j\rb\subset(0, T)$, and at
each $t_j$ the number of intersections of $\gamma_t^1$ and
$\gamma_t^2$ decreases.

\vglue-8pt

\section{Curve shortening near a closed geodesic}
\vglue-12pt
\Subsec{Eigenfunctions as $(p,q)$ satellites}
\label{sec:Birkhoffeigenfunctions}
Let $\gamma\in\Omega$ be a primitive closed geodesic of length $L$ for
a given $C^{2, \mu}$ metric $g$.  We consider a $C^{1}$ neighborhood
$\cU\subset\Omega$ and parametrize it as in \S\ref{sec:charts}.  Since
the metric $g$ is $C^{2, \mu}$, geodesics of $g$ are $C^{3,\mu}$, and
the unit normal to a geodesic will be $C^{2, \mu}$.  We can therefore
choose the local diffeomorphism $\sigma:\R/L\Z\times(-\delta, +\delta)
\to M$ so that $x\mapsto \sigma(x, 0)$ is a unit speed parametrization
of $\gamma$ and such that $\sigma_{y}(x, 0)$ is a unit normal to
$\gamma$ at $\sigma(x, 0)$.

The pullback of the metric under $\sigma$ is
\[ 
\sigma^{*}(g) = E(x, u) (dx)^{2} +2F(x, u)\, dx\,du + G(x, u)(du)^{2},
\]
for certain $C^{2,\mu}$ functions $E$, $F$, $G$.\medbreak

One can map a $C^{1}$ neighborhood of $q\cdot\gamma$ in $\Omega$ onto
a neighborhood of the origin in $C^{2}(\R/qL\Z)$ via \eqref{eq:chart}:
\begin{equation}
        u\in C^{2}(\R/qL\Z)\; \mapsto \; \alpha_{u}\in\Omega, \qquad
        \alpha_{u}(x) = \sigma(x, u(x)).
    \label{eq:qchart}
\end{equation}
In this chart the length functional $L:\Omega\to\R$ is given by
\[
L(\alpha_{u}) = \int_0^{qL} \sqrt{E(x, u)+2F(x, u)u_{x}+G(x,
  u)u_{x}^{2}}\: dx.
\]
The curve $\alpha_{u}$ will be a geodesic if and only if $u$ satisfies
the Euler-Lagrange equations corresponding to $L$.  Since we assume
$\gamma$ is already a geodesic,\break $u(x)\equiv 0$ satisfies the
Euler-Lagrange equations.  As is well-known, the second variation of
$L$ at $u=0$ is then given by
\[
d^2 L(\gamma) \cdot(v, v) = \left.\frac{d^2 L(\eps v)}{
    d\eps^2}\right|_{\eps=0} = \int_0^{qL} \bigl( v'(x)^2 -
K(\gamma(x)) v(x)^2\bigr)\: dx
\]
where $K(\gamma(x))$ is the Gauss curvature of $(M, g)$ evaluated at
$\gamma(x)$.

Consider the associated Hill's equation
\begin{equation}
        \frac{d^2 \varphi}{ dx^2}+(Q(x)+\lambda)\varphi(x) = 0\qquad   \;(x\in\R)
        \label{eq:Hill}
\end{equation}
where $\lambda$ is an eigenvalue parameter, and where
$Q(x)=K(\gamma(x))$ (although in what follows $Q \in C^{0}(\R/L\Z)$
could be arbitrary).

Let $\varphi_i(x)$ be the solutions with initial conditions
\begin{equation}
                \varphi_0(0)= 1,\quad \varphi_0'(0)=0,\quad
                \varphi_1(0)=0,\quad \varphi_1'(1)=1,
        \label{eq:y01def}
\end{equation}
and define the solution matrix
\begin{equation}
        M(\lambda; x) = 
                \left(\begin{array}{cc}
                        \varphi_0(x) &\varphi_1(x)\\ \varphi_0'(x) &\varphi_1'(x)
                \end{array}\right),
        \label{eq:Mdef}
\end{equation}
which belongs to $\textrm{SL}(2, \R)$.

If we identify the set of rays $\{\binom{ta}{tb} \mid t\geq0,
a^{2}+b^{2}=1\}$ emanating from the origin in $\R^2$ with their
intersections with the unit circle, then the linear transformation
defined by $M(\lambda; x)$ also defines a homeomorphism of the unit
circle to itself.  This homeomorphism has a rotation number
$\rho(\lambda, x)$, which is determined up to its integer part (see
\cite[\S 17.2]{Coddington-Levinson}).  To fix the integer part of
$\rho(\lambda, x)$, we require that $\rho(\lambda, 0)=0$ for all
$\lambda\in\R$ and that $\rho(\lambda, x)$ vary continuously with
$\lambda$ and $x$. The inverse rotation number of the geodesic mentioned in
the introduction is precisely $\rho(\lambda=0,x=L)$.

Since the coefficient $Q(x)$ is an $L$ periodic function, one has
\begin{align}
  M(\lambda; qL) &= M(\lambda; L)^{q}     \label{eq:Mperiodic} \\
  \intertext{and hence} \rho(\lambda, qL) &= q\rho(\lambda, L).
  \label{eq:rhoperiodic}
\end{align}

The rotation number $\rho(\lambda, L)$ is a continuous nondecreasing
function of the eigenvalue parameter $\lambda$, and thus for each
fraction $p/q$ the set of $\lambda$ with $\rho(\lambda, L)=p/q$ is a
closed interval $[\lambda^{-}_{p/q}, \lambda^{+}_{p/q}]$.  Indeed, if
$2p/q$ is not an integer, then $\lambda^{-}_{p/q} =
\lambda^{+}_{p/q}$, and we just write $\lambda_{p/q}$.

The $\lambda^{\pm}_{p/q}$ depend on the potential $Q$, and depending
on the context we will either write $\lambda_{p/q}(Q)$ or
$\lambda_{p/q}(\gamma)$ if $Q=K\circ\gamma$ is the Gauss curvature
evaluated along $\gamma$, as above.

Both for $\lambda=\lambda^{-}_{p/q}$, and $\lambda=\lambda^{+}_{p/q}$,
Hill's equation \eqref{eq:Hill} has a $qL$ periodic solution which we
denote by $\varphi^{\pm}_{p/q}(x)$. When $\lambda^{-}_{p/q} =
\lambda^{+}_{p/q}$ both solutions $\varphi_{i}(\lambda; x)$ are $qL$
periodic, and we let $\varphi^{\pm}_{p/q}(x)$ be $\varphi_{0}$,
$\varphi_{1}$ respectively.

Let $E_{p/q}(Q)$ be the two dimensional subspace of $C^{2}(\R/qL\Z)$
defined by
\begin{equation}
        E_{p/q}(Q)\isdef
                \left\{\left.
                \textstyle{c_{+}\varphi_{p/q}^{+}(x)  + c_{-}\varphi_{p/q}^{-}(x) }
                                \;\right|\;
                                        c_{\pm}\in\R
                \right\}.
        \label{eq:Epqdef}
\end{equation}
This space is determined by $Q\in C^{0}(\R/L\Z)$, i.e.\ does not
require the geodesic $\gamma$ or the surface $M$ for its definition.
It is the spectral subspace corresponding to the eigenvalues
$\lambda_{p/q}^{\pm}$ of the unbounded operator
$-\frac{d^{2}}{dx^{2}}-Q(x)$ in $L^2(\R/qL\Z)$ and as such depends
continuously on the potential $Q\in C^{0}(\R/qL\Z)$.
\begin{lemma}
  Let $\alpha_{\eps}$ be the satellite of $\gamma$ given by
  $\alpha_{\eps u}(x) = \sigma(x, \varepsilon u(x))${\rm ,} with $u(x)\in
  E_{p/q}(K\comp\gamma)${\rm ,} $u\neq 0${\rm ,} and $\eps$ sufficiently small.
  Then $\alpha_{\eps u}$ is a $(p,q)$ satellite of $\gamma${\rm ,} i.e.\ 
  $\alpha_{\eps u}\in\cB_{p, q}(\gamma)$.
\end{lemma}

\Proof   The space $E_{p/q}(Q)\subset C^{2}(\R/qL\Z)$ depends
continuously on $Q\in C^{0}(\R/qL\Z)$. For $Q(x)\equiv 0$ one has
\[
E_{p/q}(0) = \{\textstyle{ A\cos 2\pi \frac{p}{q}\frac{x}{L} + B\sin
  2\pi \frac{p}{q}\frac{x}{L}} \mid A, B\in\R \}.
\]
Choose a continuous family of $\varphi_{\theta}\in E_{p/q}(\theta
K\comp\gamma)$, $\varphi_{\theta}\neq 0$ with $\varphi_{0}(x) = \cos
2\pi \frac{p}{q}\frac{x}{L}$.
\smallbreak

We must now show that for sufficiently small $\eps\neq0$ the
corresponding curves
\[
\alpha_{\eps, \theta}(x) = \sigma(x, \varepsilon \varphi_{\theta}(x))
\]
define flat knots relative to $\gamma$.  To prove this we will show
(i) that the graph of $\varphi_{\theta}(x)$ has no double zeroes
(which implies that $\alpha_{\theta, \eps}$ is never tangent to
$\gamma$), and (ii) that the graphs of $\varphi_{\theta}(x)$ and
$\varphi_{\theta}(x-kL)$  ($k=1$, 2, $\dots $, $q-1$) have no
tangencies (which implies that $\alpha_{\theta, \eps}$ has no
self-tangencies).

The following arguments are inspired by those in
\cite[\S2]{AngenentFiedler}.

If $\lambda_{p/q}^{-}(\theta K\comp\gamma)=\lambda_{p/q}^{+}(\theta
K\comp\gamma)$, then $\varphi_{\theta}$ is a solution of Hill's
equation \eqref{eq:Hill} and cannot have a double zero without
vanishing identically.

If $\lambda_{p/q}^{-}(\theta K\comp\gamma)\neq \lambda_{p/q}^{+}
(\theta K\comp\gamma)$ then
\[
\varphi_{\theta}(x) = c_{-}(\theta)\varphi_{p/q}^{-}(x) +
c_{+}(\theta)\varphi_{p/q}^{+}(x)
\]
for certain constants $c_{\pm}(\theta)$, at least one of which is
nonzero. If one of these constants vanishes then $\varphi_{\theta}$ is
again a solution of Hill's equation and therefore cannot have a double
zero.  If both coefficients $c_{\pm}$ are nonzero then we consider
\[
u(t, x) = c_{-}(\theta)e^{\lambda^{-}_{{p/q}}t}\varphi_{p/q}^{-}(x) +
c_{+}(\theta)e^{\lambda^{+}_{{p/q}}t}\varphi_{p/q}^{+}(x).
\]
This function is a solution of the heat equation corresponding to
Hill's equation, i.e.
\[
\frac{\pt u}{\pt t}=\frac{\pt^2 u}{\pt x^2}+\theta K\comp\gamma(x)u,
\]
and by Sturm's theorem the number of zeroes of $u(t, \cdot)$ must
decrease at any moment $t$ at which $u(t, \cdot)$ has a double zero.
For $t\to\pm\infty$, $u(t, \cdot)$ is asymptotic to
$c_{\pm}e^{\lambda^{\pm}t}\varphi_{p/q}^{\pm}(x)$, and since both
$\varphi_{p/q}^{\pm}(x)$ have $2p$ zeroes in the interval $[0, qL)$
none of the intermediate functions $u(t, \cdot)$ can have a double
zero. In particular $\varphi_{\theta}=u(0, \cdot)$ only has simple
zeroes.

To prove (ii) one applies exactly the same arguments to the difference
$\varphi_{\theta}(x)-\varphi_{\theta}(x-kL)$.  The conclusion then is
that this difference either only has simple zeroes (as desired), or
else must vanish identically.  To exclude the second possibility we
observe that $\varphi_{\theta}(x)\equiv \varphi_{\theta}(x-kL)$
implies that $\varphi_{\theta}$ is an $lL$ periodic function with
$1\leq l<q$ some divisor of gcd$(k, q)$.  The number of zeroes of
$\varphi_{\theta}$ then equals $\frac{q}{l}$ times the number of
zeroes $m$ of $\varphi_{\theta}$ in its minimal period interval $[0,
lL)$.  This number $m$ is even, so the number of zeroes of
$\varphi_{\theta}$ in the interval $[0, qL)$ is a multiple of $2q/l$.
However, this number is $2p$ and so $q/l$ must be a common divisor of
$p$ and $q$.  This contradicts the hypothesis $\gcd(p, q)=1$.  \hfill\qed

\Subsec{The linearized flow at a closed geodesic}
In the chart \eqref{eq:qchart} \CS\ is equivalent to the following
parabolic equation for $u(x, t)$ (see \cite{Ang1} and also \S
\ref{sec:CS-in-coords}):
\begin{equation}
        u_{t} = \frac{u_{xx}+P(x, u) + Q(x, u) u_{x} +R(x, u) (u_{x})^{2}
                                                        + S(x, u) (u_{x})^{3}}
                                {E(x, u) +2F(x, u) u_{x} + G(x, u)(u_{x})^{2}}.
        \label{eq:CSpde}
\end{equation}
The coefficients $P$, $Q$, $R$ and $S$ are $C^{1}$ functions of their
arguments, and they satisfy
\begin{equation}
        \left\{
        \begin{array}{rl}
                \displaystyle{P(x, 0)}\hskip-8pt &\displaystyle{= Q(x, 0) =0, } \\
                \displaystyle{P_{y}(x, 0)}\hskip-8pt &\displaystyle{= K(\sigma(x, 0))}
        \end{array}
        \right.
        \label{eq:CSpdecoeffs}
\end{equation}
in which $K$ is the Gauss curvature on the surface.

One can apply classical results on parabolic equations to deduce
short-time existence for \CS\ from \eqref{eq:CSpde}.  In this section
we shall use the 
\pagebreak
local form of \CS\ to prove
\begin{lemma}\label{lem:asymptoticbirkhoff}
  If $\{\gamma_{t} \mid t\geq 0\}$ is an orbit of \CS\ which converges
  to a closed geodesic $\alpha\in\Omega${\rm ,} then for $t$ sufficiently
  large $\gamma_{t}$ is a $(p,q)$ satellite of $\alpha${\rm ;} i.e.{\rm ,} 
  $\gamma_{t} \in\cB_{p, q}(\alpha)$ for some $p, q$.  Moreover{\rm ,}
\begin{equation}
                \lambda^{-}_{p/q}(\alpha)\leq 0.
\end{equation}
If $\{\gamma_{t} \mid t\leq 0\}$ is an \/{\rm ``}\/ancient orbit\/{\rm ''}\/ of \CS\ with
$\lim_{t\to-\infty}\gamma_{t}=\alpha$ for some closed geodesic
$\alpha\in\Omega${\rm ,} then for $-t$ sufficiently large $\gamma_{t}$ is a
$(p,q)$ satellite of $\alpha$ for some $p, q$.  In this case{\rm ,}
\begin{equation}
                \lambda^{+}_{p/q}(\alpha)\geq 0.
\end{equation}
\end{lemma}
\vskip9pt

\Proof   We only prove the first statement; the second can be shown in
the same way.

If $\gamma_{t}$ converges to $\alpha$ in $C^{1}$ then we can choose
coordinates as above, and for large $t$ the curves $\gamma_{t}$
correspond to a solution $u(x, t)$ of \eqref{eq:CSpde}.  This solution
is defined for, say, $t\geq t_{0}$, and $u(\cdot, t)\to 0$ in
$C^{1}(\R/\Z)$ as $t\to\infty$.  By parabolic estimates we also have
$u(\cdot, t)\to 0$ in $C^{2}(\R/\Z)$ as $t\to\infty$.

We can write \eqref{eq:CSpde} as
\[
u_{t} = a(x, u, u_{x})u_{xx} +b(x, u, u_{x})u_{x}+c(x, u, u_{x})u
\]
where, using \eqref{eq:CSpdecoeffs} and $E(x, 0)\equiv 1$,
\begin{align*}
  a(x, u, p) &= \lp E(x, u) +2F(x, u) p + G(x, u)p^{2}\rp^{-1}, \\
  b(x, 0, 0) &= 0, \\
  c(x, 0, 0) &= K(\sigma(x, 0)).
\end{align*}
Thus \eqref{eq:CSpde} can be written as a quasilinear equation
\[  u_{t} = \cA(u)u \]
in which $\cA(u)$ is the linear differential operator
\[
\cA(u) = a(x, u, u_{x})\frac{d^{2}}{dx^{2}} +b(x, u,
u_{x})\frac{d}{dx}+c(x, u, u_{x}).
\]
For $u=0$ this operator reduces to
\[
\cA(0) = \frac{d^{2}}{dx^{2}} + K(\alpha(x))
\]
whose spectrum we have just discussed.

Since $u$ tends to zero, $u$ asymptotically satisfies the equation
$u_{t}=\cA(0)u$, and thus for some $j\geq 0$ and some constant $C\neq
0$ one has
\begin{equation}
        \lim_{t\to\infty} \frac{u(x, t)}{\|u(\cdot, t)\|_{L^{2}}} = 
        C\pfi_{j}(x)
        \label{eq:ulimit}
\end{equation}
where $\pfi_{j}(x)$ is an eigenfunction of $\cA(0)$ with 2$j$ zeroes.
See Lemmas \ref{lemma:0or1} and~\ref{lemma:lowerexponential}.  For
large $t$ the curve $\gamma_{t}$ is therefore parametrized by
\[
x \mapsto \sigma\bigl(x, \eps(t) \{ C\pfi_{j}(x)+o(1) \}\bigr),
\]
where $\eps(t)\to 0$ as $t\to\infty$. This implies that $\gamma_{t}$
is a satellite of $\alpha$.

If both eigenvalues $\lambda_{\pm}(p/q, K\circ\alpha)$ were positive
then for large $t$ one would have
\begin{align*}
  \frac{d}{dt}\|u(\cdot, t)\|_{L^{2}}^{2}
  &= \lp u(t), \cA(u(t))u(t)\rp_{L^{2}} \\
  &= (\lambda_{\pm}(p/q, K\circ\alpha) + o(1))\|u(\cdot,
  t)\|_{L^{2}}^{2}\\
  &>0
\end{align*}
which would keep $u(\cdot, t)$ from converging to zero.\hfill\qed

\section{Loops}
\vglue-12pt
\Subsec{Loops{\rm ,} simple loops{\rm ,} and filled loops}
Let $\gamma\in\Omega\setminus\Delta$ be a flat knot, and choose a
parametrization $\gamma\in C^{2}(S^{1}, M)$, also denoted by $\gamma$.
By definition a \textit{loop} for $\gamma$ is a nonempty interval $(a,
b)\subset\R$ for which $\gamma(a)=\gamma(b)$ is a transverse
self-intersection.

If we identify $S^{1}$ with $\bD$, where $\D$ is the unit disc in the
complex plane, then $\gamma(a)=\gamma(b)$ implies that any simple loop
$(a, b)\subset\R$ for $\gamma$ defines a map $\bg:S^{1}\to M$ via
\[
\bg\lp e^{ 2\pi i\frac{t-a}{b-a}}\rp = \gamma(t), \quad \text{for~}
t\in(a, b).
\]
By definition we will say that one can \textit{fill in a loop} $(a,
b)$ if the map $\bg:\bD\to M$ can be extended to a local homeomorphism
$\pfi:\D\to M$. We will always assume that a filling is at least
$C^{1}$ on $\D\setminus\{1\}$, and that $\pfi$ is a local
diffeomorphism on $\D\setminus \{1\}$.

If $\bg:S^{1}\to M$ is contractible, and one-to-one, then by the
Jordan curve theorem one can fill $\bg$. We call such a loop an
\emph{embedded loop}.

Fillings come in two varieties which are distinguished by the way they
approach the corner at the intersection $\gamma(a)=\gamma(b)$.  The
arcs $\gamma((a-\eps, a+\eps))$ and $\gamma((b-\eps, b+\eps))$ divide
a small convex neighborhood of this intersection into four pieces
(``quadrants'').  The image $\pfi(D(1, \delta))$ of a small disk will
intersect either one or three of these quadrants. If $\pfi(D(1,
\delta))$ lies in one quadrant we call the corner convex, otherwise we
call the corner concave.
\begin{figure}[h]
  \centering \includegraphics{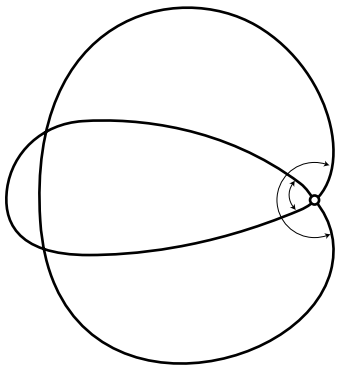}
        \caption{Convex and concave corners.}
        \label{:convex-concave}
\end{figure}

\Subsec{Continuation of loops and their fillings}
Let $\{\gamma_{\theta} \mid \theta\in[0, 1]\} \subset \Omega\setminus
\Delta$ be a smooth family of flat knots, and let $\gamma_{\theta}$
stand for smooth parametrizations of the corresponding curves.  If
$(a_{0}, b_{0})\subset\R$ is a loop for $\gamma_{\theta_{0}}$ then,
since all $\gamma_{\theta}$ 
\pagebreak
have transverse self-intersections, the
Implicit Function Theorem implies the existence and uniqueness of
smooth functions $a(\theta)$, $b(\theta)$ for which $(a(\theta),
b(\theta))$ is a loop for $\gamma(\theta)$, and such that
$a(\theta_{0})=a_{0}$ and $b(\theta_{0})=b_{0}$.  Thus \emph{any loop
  of a flat knot can be continued along homotopies of that flat knot.}

\begin{figure}[h]
  \centering \includegraphics{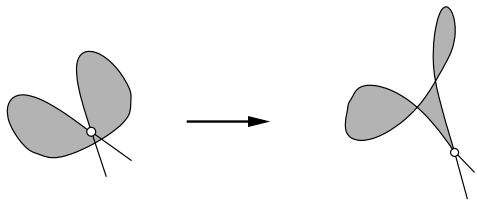}
        \caption{Inward corners may cut up fillings.}
        \label{figure:no-filling}
\end{figure}

\begin{figure}[h]
  \centering \includegraphics{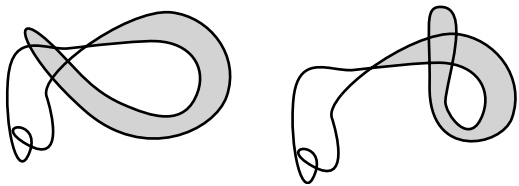}
        \caption{An embedded loop becomes nonembedded.}
        \label{figure:nonembedded}
\end{figure}

Now assume that the loop $(a_{0}, b_{0})\subset\R$ of
$\gamma_{\theta_{0}}$ has a filling: can one continue this filling in
the same way?  In general the answer is no, as the example in Figure
\ref{figure:no-filling} shows.  It is also not true that embedded
loops must remain embedded under continuation (see Figure
\ref{figure:nonembedded})


\begin{lemma}\label{lma:continuation}
  If the filling $\pfi_{0}:\D\to M$ of the loop $(a_{0}, b_{0})$ has a
  convex corner{\rm ,} then there exists a continuous family of fillings
  $\pfi_{\theta}:\D\to M$ for the loops $(a(\theta), b(\theta))$ for
  all $\theta\in[0, 1]$.
\end{lemma}

\Proof   We may assume, by changing the parametrizations if necessary,
that $a(\theta)$ and $b(\theta)$ are constant, so that $(a, b)$ is a
loop for each $\theta\in [0, 1]$.

If one has a filling of a loop for some parameter value $\theta_{0}$,
then by constructing a tubular neighborhood of the arc $\gamma:[a,
b]\to M$ one can adapt the given filling $\pfi_{0}$ to a filling
$\pfi_{\theta}$ of the loops $(a, b)$ for all $\theta$ in some
interval $(\theta_{0}-\eps_{0}, \theta_{0}+\eps_{0})$.  To obtain a
continuation from $\theta=0$ all the way to $\theta=1$ we must find a
fixed lower bound for the size of the tubular neighborhoods.  Such a
lower bound then implies a lower bound for the length $2\eps_{0}$ of
the intervals on which one can construct local continuations, so that
a finite number of such local continuations will take one from
$\theta=0$ to $\theta=1$.  We will therefore now describe the
construction of the tubular neighborhoods of the $\gamma_{\theta}$ and
the local continuations of the fillings in more detail.

Choose a suitable smooth metric $g$ on the surface $M$.  Then the
Gauss curvature of $(M, g)$ and geodesic curvatures of the
$\gamma_{\theta}$ are uniformly bounded, say by some constant $\vka$.
We can therefore choose a small $\sigma>0$ (much smaller than the
injectivity radius of $(M, g)$) such that the intersection of any disk
with radius $\sigma$ at any point $P\in M$ with any of the curves
$\gamma_{\theta}$ looks like a finite collection of straight line
segments.  More precisely, if we define the map $\phi_{P, \sigma}$
from the unit disc $\D_{P}\subset T_{P}M\cong\R^{2}$ to $M$, by
$\phi_{P, \sigma}(x) = \exp_{P}(\sigma x)$, then the preimage
$\phi_{P, \sigma}^{-1} (\gamma_{\theta})$ is a finite collection of
nearly straight arcs whose curvature is bounded by $C(\vka)\sigma$,
which can be made arbitrarily small by decreasing $\sigma$.

For each $\theta\in [0, 1]$ we construct a smooth vector field
$X_{\theta}$ along $\gamma_{\theta}$ (i.e.\ $X_{\theta}:S^{1}\to TM$
satisfies $X_{\theta}(t)\in T_{\gamma_{\theta}(t)}(M)$ for all $t\in
S^{1}$), which is nowhere tangent to $\gamma_{\theta}$, in particular
$\angle(X_{\theta}(t), \gamma_{\theta}'(t))\geq\delta$ for some
constant $\delta>0$.  This $\delta$ can be chosen independently of
$\theta$.  We can also choose the $X_{\theta}$ so that their
derivatives are uniformly bounded, i.e.\ $|\nabla^{j} X_{\theta}| \leq
C_{j}$ with $C_{j}$ independent of $\theta$.  (Note that we do not
assume that the $X_{\theta}$ vary continuously with $\theta$.)
Indeed, once one has constructed such a vector field for some value
$\theta_{1}$ of $\theta$ one can use the same vector field for all
$\theta$ in an interval containing $\theta_{1}$.  A finite number of
these intervals cover the interval $[0, 1]$, so that we really only
need a finite number of vector fields $X_{\theta}$.

Let some $\theta_{0}\in[0, 1]$ be given, and let $\pfi_{0}:\D\to M$ be
a filling for the loop $(a, b)$ of $\gamma_{\theta_{0}}$.  Since
$\gamma_{\theta_{0}}$ is the image $\pfi_{0}(\D)$ of the boundary of
the unit disc one can define an ``outward direction'' at each
$\gamma_{\theta_{0}}(t)$.  We will assume that our vector field $X$
along $\gamma_{\theta_{0}}$ is directed inward.

A tubular neighborhood is constructed from the mapping
\[
S(t, s) = \exp_{\gamma_{\theta_{0}}(t)}(s X(t)).
\]
This map is smooth from $S^{1}\times\R\to M$.  It is a local
diffeomorphism on some neighborhood $\cU= S^{1}\times [-\rho, \rho]$
of $S^{1}\times\{0\}$, where $\rho>0$ is independent of $\theta_{0}$.

If we choose $X$ so that $X(a)=X(b)$, then this map is a local
\emph{homeomorphism} from the annulus $I \times [-\rho, \rho]$ to $M$,
where $I=[a, b]/\{a, b\}$ (i.e.\ the interval $[a, b]$ with its
endpoints identified so that  $I\cong S^1$).

Now consider the curves $I\times\{s\}$ for $0\leq s<\rho$.  For
sufficiently small $s\geq 0$ there exist closed curves
$\Gamma_{s}\subset\D$ for which
\[
S(I\times\{s\}) = \pfi_{\theta_{0}}(\Gamma_{s}).
\]
Each $\Gamma_s$ is parametrized by $t\mapsto w(t, s)$, where $w$ is
the solution of
\begin{displaymath}
        S(t, s) = \pfi_{\theta_{0}}(w(t, s)).
\end{displaymath}
From $S(t, 0)= \gamma_{\theta_{0}}(t)$ it follows that $w = \exp \lp
2\pi i \frac{t-a}{b-a}\rp$ is a solution for $s=0$.  Fix $t$ and let
$s$ increase, starting at $s=0$; then, since $\pfi_{\theta_{0}} : \D
\to M$ is a local homeomorphism, one can continue the solution $w(t,
s)$ to a solution $w=w(t, s)\in\mathrm{int}(\D)$ for $0\leq
s<\sigma(t) \leq \rho$, where $\sigma(t)$ is a positive l.s.c.\ 
function of $t$.  In particular, $\sigma(t)$ is bounded from below by
some constant $\sigma>0$.  If for some $t\in I$ one has
$\sigma(t)<\rho$, then as $s\uparrow \sigma(t)$ the solution $w(t, s)$
must tend to the boundary $\pt\D$ (otherwise one could continue the
solution beyond $s=\sigma(t)$.)

It follows from $X(a)=X(b)$ that $w(a, s)\equiv w(b, s)$, and so
$t\in[a, b] \mapsto w(t, s)$ parametrizes a closed curve $\Gamma_{s}$.

\begin{proposition}
  There exists a $\sigma'>0${\rm ,} independent of $\theta$ such that all
  $\Gamma_{s}$ with $0<s<\sigma'$ are disjoint embedded curves in
  $\D$.
\end{proposition} 

\Proof   To begin, there is some $\sigma''>0$ such that none of the
smooth immersed curves $t\in\R/\Z\mapsto S(t, s)$ with $|s|\leq
\sigma''$ has a self-tangency.  This $\sigma''$ only depends on the
choice of the vector fields $X_{\theta}$, and we may thus assume that
it is independent of $\theta$.

The curves $\Gamma_{s}$ are smooth, except at $w(a, s)=w(b, s)$, where
they have a corner.  Since the derivatives of the vector fields
$X_{\theta}$ are bounded, we can find a $\sigma'''>0$ independent of
$\theta$ such that all curves $S(I\times\{s\})$ with
$|s|\leq\sigma'''$ have convex corners (in the sense that
$X_{\theta_{0}}$ points ``into the corner.'')  Hence the $\Gamma_{s}$
also have convex corners for all $0<s<\sigma'''$ for which they are
defined.

Let
\[
\sigma'=\min(\sigma, \sigma'', \sigma''').
\]

As $s$ increases from $0$ to $\sigma'$ the $\Gamma_{s}$ must remain
embedded, for the only way they can loose their embeddedness is by
first forming a self-tangency.  However, the smooth parts of the curves
$\Gamma_{s}$ are mapped to $S(I\times\{s\})$ which has no
self-tangency.  On the other hand the corner of $\Gamma_{s}$ is convex,
and
so it cannot take part in a first self-tangency.  Therefore the
$\Gamma_{s}$ remain embedded.

The $\Gamma_{s}$ are nested. Indeed, they move with velocity
\begin{displaymath}
        \frac{\pt w}{\pt s}  = \rD\pfi(w(t, s))^{-1}\frac{\pt S}{\pt s}
\end{displaymath}
which is never tangent to $\Gamma_{s}$. Thus the $\Gamma_{s}$ always
move in the same direction, which must be inward, since they start at
$\Gamma_{0}=\pt\D$.

Being nested, the $\Gamma_{s}$ can never reach the boundary $\pt\D$
again, and hence they exist for all $s\in (0, \sigma')$.  \hfill\qed

\demo{Conclusion of proof of Lemma {\rm \ref{lma:continuation}}}
By ``straightening'' the curves $\Gamma_{s}$, we see that the above construction
allows us to modify the filling $\pfi_{0}$ so that on the annulus
$e^{-\sigma'/2} \leq |w| \leq 1$ it is given by
\begin{equation}\label{eq:pfi-mod}
        \pfi_{0}(re^{i\phi}) = S(a+\frac{\phi}{2\pi}(b-a), -\ln r).
\end{equation}
For this $\varphi_{0}$ the curves $\Gamma_s$ are circles centered at
the origin.  Then we use this same expression \eqref{eq:pfi-mod} to
extend $\pfi_{0}$ to a local homeomorphism $\bar\pfi_{0} :
\D_{e^{\sigma'}} \to M$.

Since all $\gamma|_{[a,b]}$ with $\theta$ close to $\theta_0$ are
transverse to the vector field $X$, the preimage under $\bar\pfi_{0}$
of a nearby loop $\gamma_{\theta} \left|_{[a, b]} \right.$ appears as
a graph $r=r(\phi)$ in polar coordinates. One easily adapts the
filling $\pfi_{0}$ to a filling of $\gamma_{\theta} \left|_{[a, b]}
\right.$ by first mapping the unit disk $\D$ to the region enclosed by
the polar graph $r=r(\phi)$, and then composing with $\bar\pfi_{0}$.
The length of the interval of $\theta$'s for which one can do this is
bounded from below by some $\delta>0$ which is independent of
$\theta$, and hence a finite number of these local continuations will
allow one to fill $\gamma_{\theta}\left|_{[a, b]}\right.$ for all
$\theta\in[0, 1]$.  \hfill\qed

\Subsec{Loops and singularities in \CS}
In \S\ref{sec:singularities} we considered a solution $\{\gamma_{t}
\mid 0\leq t<T\}$ of \CS\ which becomes singular at time $t=T$ without
shrinking to a point.  In the notation of \S\ref{sec:singularities} we
recalled that Grayson's work implies that for every neighborhood
$\cU$ of a singular point $P_{j}$ there is a time $T_{\cU}\in(0, T)$
such that for $T_{\cU}<t<T$ the curve $\gamma_{t}$ has a loop $(a',
b')\subset[0, 1)$ with $\gamma_{t}([a', b'])$ contained in $\cU$.
Such a loop need not be simple, but one can easily extract a subloop
$(a, b)\subset(a', b')$ which is simple.  Since $\gamma_{t}|(a, b)$ is
simple it is also a fillable loop.  Still, the loop could have a
nonconvex corner, but if this is the case, and if the neighborhood
$\cU$ is homeomorphic to a disc, then we claim one can find another
loop, which is contained in $\cU$, which is simple, and whose filling
has a convex corner.

\begin{figure}[h]
   \centering \includegraphics{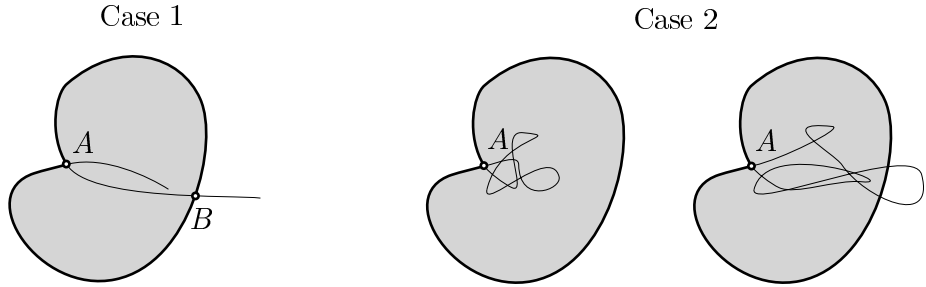}
        \caption{Finding a fillable loop with a convex corner.}
        \label{fig:goodloop}
\end{figure}

Indeed, let $\cR\subset M$ be the region enclosed by the loop, and let
$A$ be the (nonconvex) corner point of $\cR$.  Since $A$ is a
nonconvex corner point the two arcs of $\gamma\setminus\pt\cR$ enter
\emph{into} the region $\cR$ (see Figure \ref{fig:goodloop}).  There
are now two possibilities:

\demo{Case $1$} If one of these arcs exits $\cR$ again (say, at
$B\in\pt\cR$) without first forming a self-intersection, then the arc
$AB$ divides $\cR$ into two pieces, the boundary of one of which is a
simple loop with a convex corner $B$.

\demo{Case $2$} If both arcs starting at $A$ self-intersect before
leaving $\cR$, then each of these arcs contains a simple loop whose
area is strictly smaller than that of $\cR$.  If this smaller loop
still does not have a convex corner then we repeat the argument,
thereby obtaining a nested sequence of smaller simple loops.  Since
$\gamma$ only has finitely many loops this sequence must terminate
either with a simple loop with a convex corner, or with a loop as in
Case 1.

Thus we can refine the description of singularities in
\S\ref{sec:singularities} to the following:
\begin{lemma}\label{lma:filsingloop}
  If $\{\gamma_{t} \mid 0\leq t<T\}$ is a solution of \CS\ which
  becomes singular at $t=T${\rm ,} then for any $\eps>0$ there exists a
  $t_{\eps}\in (0, T)$ such that $\gamma_{t_{\eps}}$ has a convexly
  fillable loop with area no more than $\eps$.
\end{lemma}
\vglue-12pt

\begin{figure}[ht]
 \centering \includegraphics{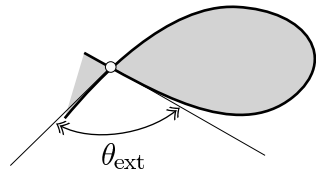}
        \caption{Definition of $\theta_{\textrm{ext}}$}
        \label{fig:thetaext}
\end{figure}

\Subsec{Decrease of area of small loops}
Let $\{\gamma_{t} \mid t_{0}\leq t<t_{1}\}$ be a solution of \CS\ with
$\gamma\in\Omega\setminus \Delta$ for all $t\in(t_{0}, t_{1})$.
Assume $\gamma_{t_{0}}$ has a fillable loop with a convex corner.
Then one can continue this loop for all values of $t\in(t_{0},
t_{1})$.  Let $\pfi_{t}:\D^{2}\to M$ be a filling of these loops.
Since $\pfi_{t}$ is a local diffeomorphism away from $1\in\D^{2}$, we
can pull the metric back from $M$ to $\D^{2}$ and define the area form
$dS_{t}$ and Gauss curvature $K_{t}$ of $\pfi_{t}^{*}(g)$, as well as
the geodesic curvature $\kappa_{t}$ and arc length $ds_{t}$ of the
boundary $\pt\D^{2}$.  The Gauss-Bonnet formula states that \[
\int_{\pt\D}\kappa_{t}ds_{t} +\iint_{\D}K_{t}dS_{t}
+\theta_{\text{ext}} =2\pi.
\]
Here $\theta_{\text{ext}}$ is the exterior angle at the corner of the
filling. See Figure \ref{fig:thetaext}.

Using this we find that the area $A(t)$ of the filling $\pfi_{t}$
satisfies
$$
  \frac{dA(t)}{dt}
   = \frac{d}{dt}\iint_{\D^{2}} \pfi_{t}^*(dS)   = -\int_{\pt\D^{2}} \kappa ds = -2\pi + \theta_{\text{ext}}
+ \int\int KdS.
$$
Since $0<\theta_{\text{ext}}<\pi$ this implies 
 $$ \frac{dA(t)}{dt}   < -\pi + (\sup_{M}K)A(t).
$$
Define
\begin{displaymath}
        \eps(g) = \frac{\pi}{2 \sup_{M}K}
\end{displaymath}
if $\smash{\sup_{M}}K>0$ and $\eps(g)=\infty$ otherwise.  We may then
conclude:
\begin{lemma}\label{lemma:Areadecrease}
  Let $\gamma_{0}\in\Omega\setminus \Delta$ have a convexly filled
  loop with area at most $\eps(g)${\rm ,} and consider the corresponding
  solution $\{\gamma_{t}\mid 0\leq t < T\}$ of \CS.  As long as the
  solution stays in $\Omega\setminus \Delta$ one can continue the
  loop{\rm ,} and its area satisfies
\begin{equation}
        A'(t)\leq\frac{\pi}{2}, \text{~and~}
        A(t)\leq \eps(g)-\frac{\pi}{2}t.
        \label{eq:Areadecrease}
\end{equation}
        
In particular the solution must either become singular or cross
$\Delta$ before $t_{*}= \frac{2\eps(g)}{\pi}$.
\end{lemma}

\section{Definition of the Conley index of a flat knot}
  \label{sec:theindex}
\vglue-12pt
\Subsec{The boundary of a relative flat knot type}
Let $\cB\subset\Omega\setminus\Delta(\Gamma)$ be a relative flat knot
type, for some $\Gamma=\{\gamma_1,\dots ,\gamma_N\}\subset\Omega$.
Throughout we will make the following assumption concerning multiple
covers
\begin{equation}\label{eq:boundarycondition}
\text{If $\alpha=m\cdot\beta\in\Omega$, $m\geq 2$, $\beta\in\Omega$,
 is tranverse to all $\gamma_i$ then $\alpha\not\in\bar{\cB}$.}
\end{equation}
We mention some examples.

\Subsubsec{$(p,q)$ satellites}
\label{sec:pq-satellites-example}
If $M$ is the sphere and $\zeta$ is the equator, then consider
$\cB=\cB_{p,q}(\zeta)$.  Let $\cU$ be the subset of the unit tangent
bundle which consists of all vectors not tangent to $\zeta$.  We have
seen in \S \ref{sec:pq-satellites} that $\cU$ has the homotopy
type of $\mathbb{T}^2$.  Any $\alpha \in\Omega$ which is transverse to
$\zeta$ lifts to a curve $\hat\alpha$ in $\cU$, and hence defines a
homotopy class $[\hat\alpha]$ in $\pi_1(\mathbb{T}^2)\cong\Z^2$.  The
homotopy class $[\hat\alpha] $ does not depend on $\alpha\in\cB$, and
hence on $\alpha\in\bar{\cB}$.  Since $\gcd(p,q)=1$ this homotopy
class is not a multiple of any other element of $\pi_1(\mathbb{T}^2)$,
and therefore $\alpha$ cannot be a multiple of another curve.  We
conclude that the relative flat knot types $\cB_{p,q}(\zeta)$ satisfy
condition \eqref{eq:boundarycondition}.

This example is easily generalized to any relative flat knot type
$\cB\subset\Omega\setminus\Delta(\Gamma)$.  Define $\cU$ to be the
unit tangent bundle of $M$ with the tangent vectors to the
$\pm\gamma_i$ removed, and assume that the homotopy class
$[\hat\alpha]\in\pi_1(\cU)$ is not a multiple of any other element of
$\pi_1(\cU)$.  Then $\bar\cB$ cannot contain multiple covers
transverse to the $\gamma_i \in\Gamma$.

\Subsubsec{Simple closed curves}
Let $\cS$ be the set of simple closed curves on $M=S^2$.  If
$\alpha=m\beta$ is a multiple cover, then any $\alpha'\in\Omega$ near
$\alpha$ must have at least one self-intersection.  Hence $\cS$
satisfies condition \eqref{eq:boundarycondition}.  However in this
case $\cU$ is the entire unit tangent bundle $T^1S^2\cong \R\P^3$ whose
fundamental group is $\Z_2$, in which $2\cdot 1=0$ and $3\cdot 1=1$,
i.e.\ in which all elements are nontrivial multiples.  So just like
$\cB_{p,q}(\zeta)$ the flat knot type $\cS$ satisfies
\eqref{eq:boundarycondition}, but it does so for different reasons.

\Subsubsec{Free satellites}
Let $p,q$ with $\gcd(p,q)=1$ be given and consider the set $\cB$ of
all $\alpha\in\Omega$ for which a simple closed curve $\zeta\in\Omega$
exists such that $\alpha$ is a $(p,q)$ satellite of $\zeta$.  Since
any two simple closed curves can be deformed into each other by
isotopy of $S^2$, the set $\cB$ is a connected component of
$\Omega\setminus\Gamma$, and hence the set of curves which are $(p,q)$
satellites is a flat knot type.  Note that, in contrast with the
example from \S \ref{sec:pq-satellites-example} the curve $\zeta$ here
is not fixed, and the set $\Gamma $ is empty.  Our current set $\cB$
is a flat knot type, while the set $\cB$ from
\S \ref{sec:pq-satellites-example} was only a relative flat knot type.

For any simple closed curve $\zeta$ the $q$ fold cover $q\cdot\zeta$
lies on the boundary $\partial\cB$ since one can approximate it by
$(p,q)$ satellites of $\zeta$.  Since there are no $\gamma_i$ in this
example, this flat knot type does not satisfy the condition
\eqref{eq:boundarycondition}.

  We consider the closure $\bar\cB$ of $\cB$ in $\Omega$ and
define
\begin{align*}
  \hB&= \bar\cB \setminus
  \{\pm m\cdot\gamma_i \mid m\geq2, i=1, \dots , N\} ,\\
  \partial\hB &= \hB\cap\partial\cB.
\end{align*}

\begin{lemma}
  Let $g\in\cM_{\Gamma}$. For any $\alpha\in\pt\hB$ a $t_\alpha>0$
  exists such that $\Phi^{(0, t_\alpha)}(\alpha) \subset \cB\cup
  \bigl(\Omega\setminus\bar\cB\bigr)$.
\end{lemma}

Recall that the curve shortening flow $\Phi^{t}$ was defined in
\ref{sec:gradflow}.
\Proof Since $\alpha\in\bar\cB$ the curve $\alpha$ has only finitely many
  crossings with any of the $\gamma_i$. Hence for some $t_1>0$ all
  $\Phi^t(\alpha)$ with $0<t<t_1$ are transverse to all $\gamma_i$.
  If $\alpha$ is not primitive, then condition
  \eqref{eq:boundarycondition} implies that $\Phi^t(\alpha) \in
  \Omega\setminus \bar\cB$.  If $\alpha$ is primitive, then we may
  assume that the $\Phi^t(\alpha)$ with $0<t<t_1$ also have transverse
  self-intersections. Hence $\Phi^t(\alpha)\in \Omega\setminus\Delta
  \subset \cB\cup \bigl(\Omega\setminus\bar\cB\bigr)$.
\Endproof\vskip4pt 

The following lemma states that orbit segments cannot touch the
boundary of a flat knot type $\cB$ ``from the inside''.
\begin{lemma}\label{lem:interiortangency}
  Let $g\in \cM_\Gamma$. If $\Phi^{[0,t]}(\alpha) \subset \hB$ and
  $t>0$ then $\Phi^s(\alpha) \in\cB$ for all $s\in (0, t)$.
\end{lemma}

\Proof Suppose for some $s\in(0, t)$ one has $\Phi^s(\alpha) \in\pt\hB$.
  Then $\Phi^s(\alpha)$ cannot be a multiple cover by condition
  \eqref{eq:boundarycondition}. By the Sturmian theorem
  $\Phi^{s'}(\alpha)$ is a flat knot rel $\Gamma$ for all $s'\neq s$
  close to $s$, and either the number of self-intersections or the
  number of intersections of $\Phi^{s'}(\alpha)$ with some $\gamma_i
  \in \Gamma$ must drop as $s'$ crosses $s$. This contradicts
  $\Phi^{[0,t]}(\alpha) \subset \hB$.
\Endproof\vskip4pt

We define the \emph{exit set} of $\hB$ to be the set $\cB^-$
consisting of those $\alpha\in\pt\hB$ for which $\Phi^{(0,
  t_\alpha)}(\alpha) \subset \Omega\setminus\bar\cB$. The complement
$\cB^+ = \pt\hB \setminus \cB^-$ is called the \emph{entry set}.

\begin{lemma}\label{lem:A}
  The sets $\cB^\pm$ do not depend on the metric $g\in\cM_\Gamma$
  chosen in their definition.
\end{lemma}

\begin{lemma}\label{lem:B}
  $\cB^-$ is a closed subset of $\hB$.
\end{lemma}

We prove these lemmas in reverse order.

\Subsubsec{Proof of Lemma {\rm \ref{lem:B}}}
We first show that $\cB^+$ is open in $\pt\hB$. Let $\alpha\in\cB^+$
be given. Then $\Phi^{(0, t_\alpha)} \subset \cB$ and in particular
$\Phi^{t_\alpha/2}(\alpha)\in \cB$.  By continuity of the local
semiflow $\Phi$ there is an open neighborhood $\cN \subset \Omega$
containing $\alpha$ such that $\Phi^{t_\alpha/2}(\cN)\subset \cB$.
Suppose some $\alpha'\in\cN$ belongs to $\cB^-$. Then there is a small
$t'\in(0, t_\alpha/2)$ such that $\Phi^{t'}(\alpha') \in
\Omega\setminus \bar\cB$. By continuity of $\Phi$ again, there is an
$\alpha''\in \cN\cap \cB$ with $\Phi^{t'} (\alpha'') \in \Omega
\setminus \bar\cB$. But then $\alpha''\in\cB$, $\Phi^{t'}(\alpha'')
\in \Omega\setminus \bar\cB$, and $\Phi^{t_\alpha/2}(\alpha'')
\in\cB$. This contradicts the Sturmian theorem. \phantom{overthere} \hfill\qed

\Subsubsec{Proof of Lemma {\rm \ref{lem:A}}}
We classify the possible curves $\alpha\in\pt\hB$ as follows.
\begin{enumerate}
\item $\alpha$ is primitive and transverse to all $\gamma_i$, but
  $\alpha$ has a self-tangency.
  
\item $\alpha$ is primitive and tangent to some $\gamma_i$ (but
  $\alpha\neq \pm \gamma_i$ by condition \eqref{eq:boundarycondition})
  
\item $\alpha=m\cdot \beta$ ($m\geq 2$, $\beta\in \Omega$) is a
  multiple cover. In this case $\alpha $ must be tangent to at least
  one of the $\gamma_i$, by condition \eqref{eq:boundarycondition}.
\end{enumerate}

The curves in Case 3 all belong to $\cB^-$, for under the \CS\ flow
they remain multiple covers, while they instantaneously become
transverse to the $\gamma_i$, so that condition
\eqref{eq:boundarycondition} forces them to leave $\bar\cB$.

The following proposition shows that we have in Case 1,
\[
\alpha\in\cB^- \Leftrightarrow \cross(\alpha, \alpha)<m_0,
\]
while in Case 2 we have
\[
\alpha\in\cB^- \Leftrightarrow \exists i : \cross(\alpha,
\gamma_i)<m_i.
\]
Thus we have a description of the exit set which is independent of the
chosen metric $g$.\hfill\qed

\begin{proposition} 
  Let $\alpha\in\hB$ be primitive with $\alpha\neq \pm \gamma_i$ for
  any $i$. If $\Cross(\alpha,\alpha)= m_0(\cB)$ and $\Cross(\alpha,
  \gamma_i) = m_i(\cB)$ for $i=1, \dots , N${\rm ,} then for some $\eps>0$
  one has $\Phi^{(0, \eps)}(\alpha) \subset \cB$.
\end{proposition}

\Proof Since $\Phi^{(0, \eps]}(\alpha) \subset \Omega\setminus
  \Delta(\gamma_1, \dots ,\gamma_N)$ we either have $\Phi^{(0,
    \eps]}(\alpha) \subset \cB$ or $\Phi^{(0, \eps]}(\alpha) \subset
  \Omega\setminus\bar\cB$.  We must show the latter cannot hold.
  Suppose it does hold. Then let $\alpha_n\in\cB$ be a sequence with
  $\alpha_n \to\alpha$. Since $\Phi^\eps(\alpha)\in\Omega \setminus
  \bar\cB$ one also has $\Phi^\eps(\alpha_n)\in\Omega \setminus
  \bar\cB$ for large enough $n\in\N$. Thus the orbit
  $\Phi^t(\alpha_n)$ crosses $\Delta(\gamma_1, \dots , \gamma_N)$ for
  some $t\in [0, \eps]$. By the Sturmian theorem one then has
\[
\#\bigl(\Phi^\eps(\alpha_n) \cap \Phi^\eps(\alpha_n)\bigr) <m_0 \text{
  or } \exists i : \#\bigl(\Phi^\eps(\alpha_n) \cap \gamma_i\bigr)
<m_i.
\]
On the other hand $\Phi^\eps(\alpha) \in \Omega \setminus
\Delta(\gamma_1, \dots , \gamma_N)$ so that for sufficiently large
$n\in\N$ one has
\begin{align*}
  \text{either}\quad&\#\bigl(\Phi^\eps(\alpha_n) \cap \Phi^\eps(\alpha_n)\bigr)
  =\#\bigl(\Phi^\eps(\alpha) \cap \Phi^\eps(\alpha)\bigr) < m_0 \\
  \text{  or}\quad & \exists i : \#\bigl(\Phi^\eps(\alpha_n) \cap \gamma_i\bigr)
  =\#\bigl(\Phi^\eps(\alpha) \cap \gamma_i\bigr) <m_i.
\end{align*}
If we now let $\eps\downto 0$ then we get
\[
\Cross(\alpha, \alpha) < m_0 \text{ or } \exists i : \Cross(\alpha,
\gamma_i) <m_i.
\]
Thus we have a contradiction, and the proposition is proved.
\Endproof\vskip4pt

Given a metric $g$ we now define
\begin{align}
  \cBl(g, \eps) &= \lb\alpha\in\cB \;\left|\; \vcenter{\hbox{$\alpha$
        has a filled loop with}
                 \vskip0.5ex%
                 \hbox{a convex corner, and area $\le\eps$}}
             \right. \rb,\\[1ex]
             \hBl(g, \eps) &= \text{the closure of}\;\cBl(g, \eps)\;
             \text{in}\; \Omega\setminus\left\{\pm q \gamma_i \mid
               q\in\N, 1\leq i\leq N\right\},
\end{align}
and we call
\[
h(\cB) = \left[^{\displaystyle \hB} \left/ _{\displaystyle \hBl(g,
      \eps)\cup \cB^-}^{}\right.  \right],
\]
the Conley-index of the component $\cB$.  Here for any closed subset
$A$ of a topological space $X$, $[X/A]$ stands for the homotopy type
of the pointed quotient space $X/A$. See \cite{Conley}.

\begin{lemma}\label{lemma:forwardinvariance}
  The set $\cB^-\cup\hBl(g, \varepsilon)$ is positively invariant
  relative to $\hB$.
\end{lemma}

\Proof Let $\Phi^{[0, t]}(\alpha) \subset \hB$ with $t>0$ and $\alpha\in
  \cB^-\cup\hBl(g, \varepsilon)$ be given.
  
  By Lemma \ref{lem:interiortangency} we have $\Phi^s(\alpha) \in \cB$
  for $0<s<t$.
  
  Fix an $s\in (0, t)$ and choose a sequence $\alpha_n\to\alpha$ with
  $\alpha_n\in \cBl(g, \eps)$. Since $\Phi^s(\alpha) \in \cB$
  continuity of the semiflow implies $\Phi^{s}(\alpha_n)\in\cB$ for
  large enough $n$. Foward invariance of $\cBl(g, \eps)$ in $\cB$ then
  implies $\Phi^{s}(\alpha_n)\in\cBl(g, \eps)$. Taking the limit as 
  $n\to\infty$ one finds $\Phi^s(\alpha)\in\hBl(g, \eps)$ for any
  $s\in (0, t)$. Taking another limit $s\to t$ one finds that
  $\Phi^t(\alpha)\in\hBl(g, \eps)$.
\hfill\qed

\begin{lemma}\label{lemma:invariance}
  The Conley index $h(\cB)$ does not depend on the metric
  $g\in\cM_\Gamma$ or the choice of $\eps>0${\rm ,} as long as
  $\eps<\eps(g)$.
\end{lemma}

This lemma justifies the absence of $g$ and $\eps$ in our notation
``$h(\cB)$'' for the index.

The proof of this lemma is essentially found in \cite{Conley}.
Our observation here is that although we do not have the desired local
compactness\footnote{K. Rybakowski has developed a version of Conley's
  theory for local semiflows on complete metric spaces, but we were
  unable to verify his ``admissibility condition,'' mainly because
  $\hB$ can contain arbitrarily long curves, and, possibly, geodesics of
  arbitrary length.} we are only trying to prove that the index is
independent of the ``index pair'' for a small class of index pairs.

We split the proof of Lemma \ref{lemma:invariance} into two pieces.
It will be convenient to write
\[
\cspace{g, \eps} = {}^{\displaystyle \hB} \left/
  _{\displaystyle{\cB^{-}\cup\hBl({g, \eps})}} \right.
\]
so that we have defined the Conley index of $\cB$ to be the homotopy
type of the quotient $\cspace{g, \eps}$, and we must now show that
this homotopy type does not depend on $g\in\cM_\Gamma$ or $\eps\in(0,
\eps(g))$.

\Subsubsec{$h(\cB)$ does not depend on $\eps$} Let 
$0<\eps_{1}<\eps_{2}<\eps(g)$ be given. Then trivially we have the
inclusion $\cBl(g, \eps_{1}) \subset\cBl(g, \eps_{2})$ which leads to
a natural mapping
\[
\cspace{g, \eps_{1}}\stackrel{f}{\longrightarrow}\cspace{g, \eps_{2}}.
\]
(Whenever $A_{1}\subset A_{2}\subset X$ are closed subsets there is a
natural mapping $X/A_{1}\to X/A_{2}$.)

We will show that this mapping is a homotopy equivalence.  For every
$\gamma\in\hB$ we define
\begin{equation}
    T_{*}(\gamma) = \inf\{ t\geq 0 \mid \Phi_{t}(\gamma)\in 
    \cB^{-}\cup\hBl(g, \eps_{1})\}
    \label{eq:Texitdef}
\end{equation}
with the understanding that $T_{*}(\gamma) =\infty$ if
$\Phi_{t}(\gamma)$ never reaches the exit set or $\hBl(g, \eps_{1})$.

\begin{proposition}\label{prop:Tcontinuous}
  The function $T_{*}:\hB \to [0, \infty]$ is continuous.
\end{proposition}
\Proof We check that both conditions $T_{*}(\gamma)<M$ and 
  $T_{*}(\gamma)>M$ define open subsets of $\hB$.
  
  If $T_{*}(\gamma)<M$ for some $0<M<\infty$, then
  $\Phi^{T_{*}(\gamma)}(\gamma)$ belongs to $\cB^{-}$ or $\hBl(g,
  \eps_{1})$. In the first case the orbit immediately leaves $\cB$, and  so
  there exists a $t_{0}\in(T_{*}(\gamma), M)$ with
  $\Phi^{t_{0}}(\gamma) \in\Omega\setminus \bar\cB$. By continuity of
  the semiflow $\Phi$ the same is then true for all $\gamma'$ near
  $\gamma$, so that $T_{*}(\gamma')<t_{0}<M$ holds on a neighborhood
  of $\gamma$.
  
  Consider the second case, in which $\Phi^{T_{*}(\gamma)}(\gamma) \in
  \hBl(g, \eps_{1})$.  If $T_*(\gamma)=0$ 
  \pagebreak
  then it is possible that
  $\gamma=\Phi^{T_{*}(\gamma)}(\gamma)$ lies on $\cB^+$. When this
  happens $\Phi^t(\gamma)$ must immediately enter $\cB$ and
  hence$\cBl(g, \eps_1)$, by forward invariance of $\hBl(g, \eps)$
  relative to $\hB$.  If on the other hand $T_*(\gamma)>0$ then
  $\Phi^{T_{*}(\gamma)}(\gamma)$ cannot lie on $\cB^+$. By assumption
  it does not lie on $\cB^-$ either, and thus it lies in $\cBl(g,
  \eps_1)$.  It follows from Lemma \ref{lemma:Areadecrease} that at
  $t=T_{*}(\gamma)$ the orbit $\Phi^{t}(\gamma)$ develops a convexly
  filled loop with area $\leq\eps_{1}$, and that for $t>T_{*}(\gamma)$
  the loop has area $\leq\eps_{1}-\frac{\pi}{2}(t-T_{*}(\gamma))$
  which is strictly less than $\eps_{1}$. Invoking continuity of the
  semiflow we conclude again that this condition also holds for
  $\gamma'$ near $\gamma$.
  
  Conversely, if $T_{*}(\gamma)>M$, then the (compact) orbit segment
  $\{\gamma_{t} \mid 0\leq t\leq M\}$ is contained in $\hB\setminus
  \hBl(g, \eps_{1})$ which is open relative to $\hB$.  Once more
  continuity of the semiflow guarantees that this is also the case for
  $\gamma'$ close to~$\gamma$. \phantom{moverover}
\Endproof\vskip4pt 

It follows from Lemma \ref{lemma:Areadecrease} that for all
$\gamma\in\cBl(g, \eps_{2})$ one has $T_{*}(\gamma)\leq
\frac{2}{\pi}(\eps_{2}-\eps_{1})$. By continuity this also holds for
all $\gamma\in \hBl(g, \eps_{2})$.  Now define
\[
T_{0}(\gamma) = \min\left(\frac{2}{\pi}(\eps_{2}-\eps_{1}),
  T_{*}(\gamma)\right) 
\]
and consider the following homotopy ($0\leq\lambda\leq 1$):
\[
G_{\lambda}: \cB^{-}\cup \hB \to \cB^{-}\cup \hB, \qquad \gamma
\mapsto \Phi^{\lambda T_{0}(\gamma)}(\gamma).
\]
Then
\begin{itemize}
\item $G_{0}$ is the identity map on $ \cB^{-}\cup \hB $,
        
\item $G_{\lambda}$ maps $\cB^{-}\cup\hBl(g, \eps)$ to itself for
  every $\eps \in (0, \eps(g))$ (by forward invariance of $\cB^- \cup
  \hBl$, Lemma \ref{lemma:forwardinvariance}),
        
\item $G_{1}$ maps $\cB^{-}\cup\hBl(g, \eps_{1})$ into
  $\cB^{-}\cup\hBl(g, \eps_{2})$
\end{itemize}
and it is easily verified from these facts that $G_{1}$ is a homotopy
inverse of $f$.

\Subsubsec{$h(\cB)$ does not depend on the metric}\label{sssec:metric-independence} Let $g_{1}, g_{2}\in \cM_\Gamma$
be two given metrics.  Then, since the surface $M$ is compact there
exists a constant $A>0$ such that one has $g_{1}\leq A g_{2}$ and
$g_{2}\leq Ag_{1}$ pointwise on $M$. In particular the area form of
either metric is bounded by $A^{2}$ times the area form of the other.
We therefore have the following inclusions
\[
\cB^{-}\cup\hBl(g_{1}, \eps) \subset \cB^{-}\cup\hBl(g_{2}, A^{2}\eps)
\subset \cB^{-}\cup\hBl(g_{1}, A^{4}\eps) \subset
\cB^{-}\cup\hBl(g_{2}, A^{6}\eps)
\]
and corresponding natural maps
\[
\cspace{g_{1}, \eps} \stackrel{f}{\longrightarrow} \cspace{g_{2},
  A^{2}\eps} \stackrel{g}{\longrightarrow} \cspace{g_{1}, A^{4}\eps}
\stackrel{h}{\longrightarrow} \cspace{g_{2}, A^{6}\eps} .
\]
Now it follows from the previous section that for sufficiently small
$\eps>0$ the compositions $g\circ f$ and $h\circ g$ are homotopy
equivalences, so that $g$ has a left and right homotopy inverse. Hence
$g$ is a homotopy equivalence.  \Endproof\vskip4pt   

\Subsec{Virtual satellites and the modified Conley index of a
  relative flat knot}\label{sec:virtualsatellites} Let
$\cB\subset\Omega\setminus\Delta(\Gamma)$ be a relative flat knot for
some $\Gamma=\{\gamma_1, \dots , \gamma_N\}\subset \Omega$.  Define
$\cM_\Gamma$ as in \S\ref{sec:intro-resonance}, and, as in
\S\ref{sec:intro-resonance},  order the $\gamma_i$ so that for
$i=1, \dots , m$ there exists $p_i/q_i$ with
\begin{equation}
    q_i\gamma_i\in \pt\cB \text{ and } \cB\subset\cB_{p_i,q_i}(\gamma_i)
    \label{eq:hypo-order}
\end{equation}
while no such $p_i/q_i$ exist for $i=m+1, \dots , n$.  By Lemma
\ref{lem:pqisrs} the $p_i/q_i$ are uniquely determined.

We impose the nonresonance condition \eqref{eq:nonresonance} and for
any $I\subset\{1, \dots ,m\}$ we define $\cM(\alpha; I)$ for
$\alpha\in\cB$ as in \S\ref{sec:intro-resonance}. Since $\cM(\alpha;
I)$ does not depend on $\alpha\in\cB$ we will write $\cM(\cB; I)$ for
$\cM(\alpha; I)$. Our discussion of the rotation number in
\S\ref{sec:Birkhoffeigenfunctions} shows that condition
\eqref{eq:rotationnumbers} is equivalent to
\begin{equation}
    \lambda_{p_{i}/q_{i}}^+(\gamma_i) < 0 \text{~for $i\in I$, and~}
    \lambda_{p_{i}/q_{i}}^-(\gamma_i) > 0 \text{~for $i\in I^c$,}
    \label{eq:5A}
\end{equation}
where $I^c=\{1, \dots , m\}\setminus I$.

For the moment write $q=q_{i}$ and $\gamma=\gamma_{i}$.  Let
$\cU\subset\bar\cB$ be a closed neighborhood in $\bar\cB$ of $q\gamma$
which is small enough for $q\gamma$ to be the only geodesic in $\cU$,
and for $\cU\cap\hBl(g, \varepsilon)$ to be empty (for some
$\varepsilon\in(0, \varepsilon(g))$ which we keep fixed throughout
this section).

Define
\begin{align}
  \cU_{\textrm{transient}} &= \{\alpha\in\cU\cap\bar\cB \mid
  \exists t>0 : \Phi^{t}(\alpha)\in\cB\setminus\cU \} \\
  \cU^{\#} &= \cU\setminus \lp \cU_{\textrm{transient}}\cup\{q\gamma\}
  \rp.
\end{align}
The conditions $\alpha\in\bar\cB$ and $\Phi^{t}(\alpha)\in\cB$ imply
that $\Phi^{(0, t]}\subset\cB$, since orbits cannot leave and then
enter $\bar\cB$ again. Thus $\cU^{\#}$ consists of those
$\alpha\in\cU$ which do not leave $\cU$ before leaving $\cB$.

Clearly $\cU_{\textrm{transient}}$ is open, so that $\cU^{\#}$ is closed in
$\hB$ and $\cU^{\#}\cup\{q\gamma\}$ is closed in $\bar\cB$.

By construction $\cU^{\#}$ is positively invariant relative to $\hB$;
$\cU\setminus\{q\cdot\gamma\}$ is positively invariant relative to
$\hB$ if and only if $\cU=\cU^{\#}\cup\{q\cdot\gamma\}$, or,
equivalently,
\[
\cU^{\#} = \cU\setminus\{q\cdot\gamma\}.
\]

\begin{lemma}
  $\cU^{\#}\cup\{q\gamma\}$ is a neighborhood in $\bar \cB$ of
  $q\cdot\gamma$.
\end{lemma}
\Proof If $\cU^{\#}\cup\{q\gamma\}$ is not a neighborhood of $q\gamma$ then
  a sequence $\alpha_{n}\in\bar\cB\setminus\cU^{\#}$ with
  $\lim_{n\to\infty} \alpha_{n} = q\gamma$ must exist. Since $\cU$ is
  assumed to be a neighborhood we may assume that all
  $\alpha_{n}\in\cU$, and thus
  $\alpha_{n}\in\cU_{\textrm{transient}}$.  Then $t_{n}'>0$ exist with
  $\Phi^{t_{n}'}(\alpha_{n})\in\cB\setminus \cU$. Choose $t_{n}$ to be
  the largest $t\in (0, t_{n}')$ with $\Phi^{[0, t_{n}]}(\alpha_{n})
  \subset \cU$. In particular $\Phi^{t_{n}}(\alpha_{n})\in\pt\cU$.
  
  Since $q\gamma$ is a fixed point for \CS, we have
  $\lim_{n\to\infty}t_{n}=\infty$.  By parabolic estimates
  \ref{lemma:parabolicestimate} we can extract a convergent
  subsequence of the sequence of solutions $\{\beta_{n}(t) \isdef
  \Phi^{t_{n}+t}(\alpha_{n}) \mid -t_{n}<t\leq0\}$.  The limit is an
  ``ancient orbit'' $\{\beta(t) \mid -\infty<t\leq 0\}$ which remains
  in $\cU$, and which reaches $\pt\cU$ at $t=0$.  The $\alpha$-limit
  of such an orbit must be $q\gamma$ (being the only closed geodesic
  in $\cU$) but this contradicts $\lambda^{+}_{p/q}(\gamma)<0$ and
  Lemma \ref{lem:asymptoticbirkhoff}.
\Endproof\vskip4pt 

For any neighborhood $\cU\subset\bar\cB$ of $q\gamma$ and $T>0$ we
define
\[
\cU^{\#T} = \{\Phi^{T}(\alpha) \mid \Phi^{[0,
  T]}(\alpha)\subset\cU^{\#}\}.
\]
As $T$ increases the set $\cU^{\#T}$ shrinks.  In general $\cU^{\#T}$
is not a neighborhood of $q\cdot\gamma$; in fact, due to the
regularizing effect of the heat flow, $\cU^{\#T}$ will have empty
interior.
\begin{lemma}
  Let $\cU, \cV\subset\bar\cB$ be neighborhoods of $q\gamma$.  Then
  for sufficiently large $T>0$ one has $\cU^{\#T}\subset\cV^{\#}$ and
  $\cV^{\#T}\subset\cU^{\#}$.
\end{lemma}
\Proof   We need only prove the first inclusion, and we may of course
assume that the neighborhoods $\cU=\cU^{\#}\cup\{q\gamma\}$,
$\cV=\cV^{\#}\cup\{q\gamma\}$ are positively invariant relative to
$\cB$.

Arguing by contradiction we assume that there exists a sequence
$\alpha_{k}\in\cU$ with $\Phi^{[0, k]}(\alpha_{k})\subset
\cU\setminus\cV$. Parabolic estimates yield an \emph{a priori} bound
for $\frac{\pt k}{\pt s}$ on the curves $\Phi^{1}(\alpha_{k})$, and
thus we can extract a convergent subsequence from the solutions
$\beta_{k}(t)=\Phi^{t+1}(\alpha_{k})$ of \CS. The limit would then be
an orbit of \CS\ which stays in $\overline{\cU\setminus\cV}$, in
particular its $\omega$-limit would be a closed geodesic other than
$q\gamma$ in $\cU$, which by assumption does not exist.\Endproof\vskip4pt

Let $I\subset\{1, \dots , m\}$ and $g\in\cM_\Gamma(\cB;I)$ be given.
For each $i\in I$ we choose a sufficiently small neighborhood
$\cU_{i}$ of $q_i\gamma_i$ and we set
\[
\cU^I \isdef \bigcup_{i\in I} \cU_i^{\#}.
\]
We will assume that the $\cU_{i}\setminus\{q_i\gamma_i\}$ are forward
invariant relative to $\hB$, i.e.\ $\cU_{i}=\cU_{i}^{\#}$.

\begin{definition}
  The \emph{modified Conley index} of the relative flat knot type
  $\cB$ is
\[ h^I(\cB) =
\left[ ^{\displaystyle{\hB}} \left/ _{\displaystyle{\cB^-\cup \hBl(g,
        \varepsilon) \cup \cU^I} }\right.  \right].
\]
\end{definition}

Our previously defined Conley index $h(\cB)$ is contained in this
definition as the special case in which $I\subset\{1, \dots , m\}$ is
empty.

\begin{lemma}\label{lemma:1}
  For sufficiently small $\cU_{i}$ and $\varepsilon>0$ the index
  $h^I(\cB)$ does 
\pagebreak  
not depend on either $\varepsilon${\rm ,} the metric
  $g\in\cM_\Gamma(\cB;I)$ or the neighborhoods $\cU_{i}$.
\end{lemma}

\Proof   We may assume that $\cU\supset \cV$ for otherwise we choose a
smaller neighborhood $\cW\subset \cU\cap\cV$ and compare the indices
$h^{I}(\cB)$ obtained by using $\cU$ and $\cV$ with the index obtained
by using $\cW$.

Choose a sufficiently large $T>0$ so that $\cU^{\#T}\subset \cV$ and
as before in \eqref{eq:Texitdef}   define
\[
T_{*}(\alpha) = \inf \{t\geq 0 \mid \Phi^{t}(\alpha)\in
\cB^{-}\cup\hBl \}.
\]
We showed in Proposition \ref{prop:Tcontinuous} that the exit time
$T_{*}(\alpha)$ is a continuous function with values in $[0, \infty]$.
The family of maps
\begin{displaymath}
    F_{\theta}(\alpha) \isdef 
    \Phi^{\theta\min(T_{*}(\alpha), T)}(\alpha)
\end{displaymath}
with $\theta\in [0, 1]$ is therefore a continuous homotopy $F_\theta:
\mathrm{id}\cong F_1$ of maps of the pairs $(\hB, \cB^- \cup \hBl \cup
\cU)$ and $(\hB, \cB^- \cup \hBl \cup \cU)$.  From $\cU^{\#T}\subset
\cV$ we conclude that $F_1$ maps the quotient $\hB/(\cB^-\cup \hBl
\cup \cU)$ to $\hB/(\cB^-\cup \hBl \cup \cV)$, and is a homotopy
inverse for the inclusion induced map from $\hB/(\cB^-\cup \hBl \cup
\cV)$ to $\hB/(\cB^-\cup \hBl \cup \cU)$.  \phantom{moveover} \hfill\qed

\Subsec{The Conley index of a virtual satellite}

\begin{lemma}
  The homotopy type of $\cU_{i}^{\#}/(\cU_{i}^{\#}\cap\cB^{-})$ is
  that of $S^{1}\times S^{2p_{i}-1}/S^1\break\times\{\mathrm{pt}\}$.
          Consequently the homotopy type of $\cU^{I}/(\cU^{I}\cap\cB^{-})$ is
  given by
        \[
        \bigvee_{i=1}^{m}[S^1\times
        S^{2p_{i}-1}/S^1\times\{\mathrm{pt}\}].
        \]
\end{lemma}

We will call the homotopy type of
$\cU_{i}^{\#}/(\cU_{i}^{\#}\cap\cB^{-})$ the Conley index of the
virtual satellite of $q_{i}\gamma_{i}$ in $\cB$.

In the following proof we omit the subscript $i$ and write $\cU$
instead of $\cU_{i}$, etc.

The same arguments as in Corollary \ref{lemma:1} show that
$\cU^{\#}/(\cU^{\#}\cap\cB^{-})$ is independent of both the metric
$g$, provided $\gamma$ is a geodesic with
$\lambda_{p/q}^{+}(\gamma)\break<0$, and the neighborhood $\cU$, provided it
is sufficiently small, meaning that it should not contain any other
closed geodesics besides $q\gamma$ and be disjoint from~$\hBl$.  Thus
we may choose our metric so that a neighborhood of $\gamma$ in the
surface $M$ is isometric to a part of the surface of revolution whose
metric is given by
\[
(ds)^2=e^{y^2/2}\{(dx)^2+(dy)^2\}, \quad (x, y)\in (\R/\Z)\times\R,
\]
where the ``waist'' $y=0$ corresponds to $\gamma$.  Curves
$\alpha\in\Omega$ which are $C^1$ close to $q\gamma$ are then given by
graphs of functions $u\in C^{2}(\R/q\Z)$, $\alpha(x) = (x, u(x))$
(such a graph wraps itself $q$ times around the waist $\{y=0\}$).  In
this section we will identify a small neighborhood of $q\cdot \gamma
\in \Omega$ with an open neighborhood of $u\equiv 0$ in
$C^{2}(\R/q\Z)$ without explicitely mentioning the identification
again.

Curve shortening for such graphs is equivalent to the PDE
\begin{equation}
    \frac{\pt u}{\pt t} = \frac{u_{xx}}{1+u_x^{2}}-u.
    \label{eq:csnearwaist}
\end{equation}
We choose
\begin{equation}
  \cU  = \cN_{\sigma}\cap\bar\cB 
\end{equation}
where 
\begin{equation} \cN_\sigma  = \{u\in C^{2}(\R/q\Z) \mid
  \sup_{x}|u(x)|\leq \sigma, \sup_{x}|u'(x)|\leq\sigma \},
\end{equation}
and $\sigma$ is sufficiently small.
\begin{lemma}
  $\cN_\sigma$ is invariant for the \CS\ flow{\rm ,} so that $\cU$ contains
  no transient part\/{\rm ,}\/ i.e.\ $\cU=\cU^{\#}\cup\{q\gamma\}$.
\end{lemma}
\Proof   The maximum principle implies that any solution of
\eqref{eq:csnearwaist} with $|u(x, 0)|\leq \sigma$ satisfies $|u(x,
t)|\leq \sigma e^{-t}$, since $\pm\sigma e^{-t}$ are sub- and
supersolutions for \eqref{eq:csnearwaist}.

By differentiating \eqref{eq:csnearwaist} one finds that $v=u_x$
satisfies
\[
v_t= \frac{v_{xx}}{1+u_x^2}-\frac{2u_{x}}{(1+u_x^2)^{2}}v_{x}-v
\]
so that the maximum principle again implies that
$\sup_{x}|v(x,0)|\leq\sigma$ leads to $\sup_{x}|v(x,t)|\leq\sigma
e^{-t}$. \Endproof\vskip4pt

We will identify $\Delta\subset\Omega$ with those $u\in C^{2}(\R/q\Z)$
which correspond to a curve $\alpha_u\in\Delta$.
\begin{lemma}\label{lemma:cones}
  Both $\cB^{-}$ and $\Delta$ are cones in $C^{2}(\R/q\Z)$.
\end{lemma}
\Proof   A function $u\in C^{2}(\R/q\Z)$ belongs to $\Delta$ if it
either has a multiple zero or if for some $k=1$, $\dots $, $q-1$ the
function $u(x)-u(x-k)$ has a multiple zero. This clearly holds for $u$ if and only if it holds for $\lambda u$,
for any $\lambda\neq 0$. Thus $\Delta$ is a cone.

Near $q\gamma$ the set $\cB^{-}$ consists of those $u\in\Delta$ which
have fewer self-inter\-sections, or fewer intersections with $u=0$ than a
general $u\in\cB$ has. This condition also holds for both $u$ and
$\lambda u$ or for neither.  \Endproof\vskip4pt   

Any $u\in C^2(\R/q\Z)$ has a Fourier series of the form
\begin{equation}
    u(x) = \sum_{k=0}^{\infty} \Re\lp u_k e^{2k\pi ix/q}\rp
    \label{eq:fourier-u}
\end{equation}
with $u_0\in\R$ and $u_k\in\C$ for $k\geq 2$.  The embedding
$C^2\hookrightarrow W^{2, 2}$ implies that
\begin{equation}
    \int_{0}^{q} \{u(x)^2+\frac{q^4}{16\pi^4}u''(x)^2\}dx
     = \sum_{k=0}^{\infty} (1+k^{4})|u_k|^2<\infty.
    \label{eq:w22norm}
\end{equation}
We now define
\[
\cS_{\varepsilon}=\lb u\in C^2(\R/q\Z) \;\left|\; \sum_{k=0}^{\infty}
  (1+k^{4})|u_k|^2 =\varepsilon^2 \right.  \rb;
\]
i.e., $\cS_{\varepsilon}$ is the intersection with $C^2(\R/q\Z)$ of
the sphere of radius $ \varepsilon$ in $W^{2,2}(\R/q\Z)$ with norm
given by \eqref{eq:w22norm}.

Since $W^{2,2}\hookrightarrow C^1$ one has $\cS_{\varepsilon} \subset
\cN_{\sigma}$ for small enough $\varepsilon>0$.

Lemma \ref{lemma:cones} implies that $(\cS_{\varepsilon}\cap\hB,
\cS_{\varepsilon}\cap\cB^{-})$ is a deformation retract of
$(\cU^{\#},\cU^{\#}\cap\cB^{-})$, the deformation going along rays
through the origin in $ C^2(\R/q\Z)$.  We therefore have a homotopy
equivalence
\[
\cU^{\#}/(\cU^{\#}\cap \cB^-) \cong (\cS_{\varepsilon} \cap \cU^{\#})
/ (\cS_{\varepsilon} \cap \cU^{\#} \cap \cB^{-}).
\]
For small enough $\varepsilon>0$ one has $\cS_{\varepsilon}\subset
\cN_\sigma$, so that $\cS_{\varepsilon}\cap\cU^{\#} =
\cS_{\varepsilon} \cap \hB$, and $\cS_{\varepsilon}\cap\cU^{\#}\cap
\cB^{-} = \cS_{\varepsilon} \cap \cB^{-}$. Hence we have a further
homotopy equivalence
\[
\cU^{\#}/(\cU^{\#} \cap \cB^-) \cong (\cS_{\varepsilon} \cap \hB) /
(\cS_{\varepsilon} \cap \cB^-) .
\]

The linear heat equation induces a continuous semiflow on
$\cS_{\varepsilon}$: for any $u\in\cS_{\varepsilon}$ let
\begin{equation}
    u(t,x) =  \sum_{k=0}^{\infty} \Re\lp u_k e^{2k\pi ix/q-4k^2\pi^2/q^2t}\rp
    \label{eq:fourierheatflow}
\end{equation}
be the solution of $u_{t}=u_{xx}$ starting from $u$, and define
$\Psi^t(u)$ to be the radial projection of $u(t, \cdot)$ onto
$\cS_{\varepsilon}$, so that 
\[
\bigl(\Psi^tu\bigr)(x) = \varepsilon \frac{u(t, x)\qquad}{\|u(t,
  \cdot)\|_{W^{2,2}}}.
\]

We will refer to $\Psi^{t}$ as the projected heat flow.

The essential insight which allows us to determine the homotopy type
of $(\cS_{\varepsilon} \cap \hB) / (\cS_{\varepsilon} \cap \cB^-)$ and
hence of $\cU^{\#}/(\cU^{\#}\cap\cB^{-})$ is that $(\cS_{\varepsilon}
\cap \hB, \cS_{\varepsilon} \cap \cB^-)$ turns out to be an index pair
for the projected heat flow $\Psi^{t} : \cS_{\varepsilon} \to
\cS_{\varepsilon}$ which isolates the invariant set
\[
\cC \isdef \left\{\left.  \Re \left( u_p e^{2\pi i px/q} \right)
    \;\right|\; u_p\in\mathbb{C},
  |u_{p}|=\frac{\varepsilon}{\sqrt{1+p^4}} \right\}.
\]
This invariant set is a normally hyperbolic circle whose unstable
manifold is $2p$ dimensional, so one expects its Conley index to be
$\left[S^1\times S^{2p-1} /S^{1}\times\mathrm{pt}\right]$.  Since we
do not have the required compactness hypothesis of \cite{Conley}, we
must prove these statements by hand, essentially verifying that
Conley's arguments still go through in our setting.

For $u\in\cS_{\varepsilon}$ we define the quantities
\begin{align*}
  w_{+}(u)&= \varepsilon^{-2}\sum_{k>p} \lp 1+k^{4}\rp |u_k|^2, \\
  w_{-}(u)&= \varepsilon^{-2}\sum_{0\leq k<p} \lp 1+k^{4}\rp |u_k|^2, \\
  w_{p}(u)&= \varepsilon^{-2} \lp 1+p^{4}\rp |u_p|^2.
\end{align*}
By definition we have
\[
w_{-}(u) + w_{p}(u) + w_{+}(u) = 1
\]
for all $u\in \cS_{\varepsilon}$.

\begin{lemma}
  Along any orbit $\Psi^t(u)$ of the projected heat flow one has
\begin{align}
  \frac{dw_+(\Psi^tu)}{dt} &\leq -Cw_+(1-w_+)<0 ,\label{eq:w+growth}\\
  \frac{dw_-(\Psi^tu)}{dt} &\geq -Cw_-(1-w_-)>0. \label{eq:w-growth}
\end{align}
\end{lemma}
\vglue5pt

\Proof   Let $u(t, x)$ be the solution to the linear heat equation
starting at $u\in\cS_{\varepsilon}$ given by
\eqref{eq:fourierheatflow}, so that $u(t,x)=\sum_{k\geq 0}\Re
(u_k(t)e^{2\pi ikx/q})$ with $u_k(t)=e^{-4\pi^2k^2/q^2t}u_k(0)$. Then
one has
\begin{align*}
  \frac{dw_+(u(t,\cdot))}{dt}
  &=\sum_{k>p}\lp 1+k^4\rp \Re(2u_k'(t)u_k(t))  \\
  &\leq-\frac{8\pi^2(p+1)^{2}}{q^2} w_+(u(t,\cdot)) \intertext{and}
  \frac{dw_-(u(t,\cdot))}{dt}
  &\geq -\frac{8\pi^2(p-1)^{2}}{q^2}  w_-(u(t,\cdot)),  \\
  \frac{dw_p(u(t,\cdot))}{dt} &= -\frac{8 \pi^2 p^{2}}{ q^2}
  w_p(u(t,\cdot)).
\end{align*}
Using
\[
w_{\pm}(\Psi^{t}(u))=\frac{w_\pm(u(t,\cdot))}{w_+(u(t,\cdot))+w_p(u(t,\cdot))+w_-(u(t,\cdot))}
\]
one then arrives at \eqref{eq:w+growth} and \eqref{eq:w-growth} with
$C=8\pi^2(2p\pm 1)/q^{2}$.  \Endproof\vskip4pt   

Consider the sets
\begin{align*}
  \cV_{\rho} &= \{u\in\cS_{\varepsilon}\cap\hB \mid w_{+}(u)\leq \rho\} ,\\
  \cV_\rho^- &= \{u\in\cV_{\rho} \mid w_{-}(u)\geq \rho\} .
\end{align*}
The differential inequalities \eqref{eq:w+growth}, \eqref{eq:w-growth}
imply that $(\cV_{\rho}, \cV_{\rho}^-)$ is an index pair.  It isolates
the same invariant set $\cC$ as $(\cU^{\#}, \cU^{\#} \cap \cB^-)$, so
one expects $\cV_{\rho}/\cV_{\rho}^-$ and $\cU^{\#}/(\cU^{\#} \cap
\cB^-)$ to have the same homotopy type.  To prove this we exhibit
homotopy equivalences obtained by ``flowing along'' with $\Psi^{t}$
\begin{align}
  \cU^{\#}/(\cU^{\#}\cap \cB_-)
  & \longrightarrow \cV_{\rho}/(\cV_{\rho} \cap \cB^-), \label{eq:homotopy1}\\
  \intertext{and} \cV_{\rho}/(\cV_{\rho} \cap \cB^-) & \longrightarrow
  \cV_{\rho}/\cV_{\rho}^-. \label{eq:homotopy2}
\end{align}

\Subsubsec{Construction of the homotopy
  equivalence \eqref{eq:homotopy1}} We define
\begin{displaymath}
    t_{*}(u) = \inf \{ t\geq 0 \mid w_{+}(\Psi^t(u))\leq\rho \}.
\end{displaymath}
\begin{proposition}\label{prop:A}
  The function $t_{*}$ is continuous and finite on
  $\cS_{\varepsilon}\cap\hB$.
\end{proposition}
\Proof   We first observe that one has $w_{+}(u)<1$ for every $u\in
\cU^{\#}$.  Indeed, if $w_{+}(u)=1$ then $w_{-}(u)=w_{p}(u)=0$ and so the
Fourier series \eqref{eq:fourier-u} only contains terms with $k\geq
p+1$.  Then $u$ has at least $2(p+1)$ sign changes and cannot belong to
$\hB$ or $\cU^{\#}\subset \hB$.

The differential inequality \eqref{eq:w+growth} implies that
$w_{+}(\Psi^t(u))$ will decrease to $\rho$ in finite time so that  $t_{*}$
is finite. Moreover $\frac{d}{dt}w_{+}<0$ implies that the time at
which $w_{+}(\Psi^t(u))$ becomes equal to $\rho$ depends continuously
on $u$. \Endproof\vskip4pt   

The family of maps $G_{\theta}(u) = \Psi^{\theta t_*(u)}(u)$ is a
continuous homotopy $G_{\theta}: \mathrm{id} \cong G_1$. The final map
$G_1$ sends $\cU^{\#}/(\cU^{\#}\cap \cB^-)$ to
$\cV_\rho/(\cV_{\rho}\cap \cB^-)$ and is a homotopy inverse for the
inclusion induced map $\cV_\rho / (\cV_{\rho}\cap \cB^-) \to
\cU^{\#}/(\cU^{\#}\cap \cB^-)$.

\Subsubsec{Construction of the homotopy equivalence
  \eqref{eq:homotopy2}} To construct a map from left to right in
\eqref{eq:homotopy2} we observe
\begin{proposition}\label{prop:vrho}
  If $\rho>0$ is small enough then $\cV_{\rho}\cap \cB^-
  \subset\cV^-$.
\end{proposition}
\Proof   We consider the sets
\begin{align}
  \cW_{\rho} &= \{u\in\cS_{\varepsilon} \mid
  w_+(u)\leq\rho, w_-(u)\leq\rho \}, \label{eq:cWdef1}\\
  \cW_{\rho}^- &= \{u\in\cS_{\varepsilon} \mid w_+(u)\leq\rho,
  w_-(u)=\rho \} .\label{eq:cWdef2}
\end{align}
By definition $\cW_{\rho}$ is a $W^{2,2}$ neighborhood of $\Gamma$
which can be made as small in $W^{2,2}$ as desired by decreasing
$\rho$.  Since $\Gamma$ is compact, and since $\Delta$ is closed in
$C^1$ and thus also in $W^{2,2}$, we conclude that, for sufficiently
small $\rho>0$, $\cW_{\rho}$ and $\Delta$ are disjoint.  Since
$\cV_\rho\setminus\cV_\rho^- \subset \cW_\rho$ and
$\cB^-\subset\Delta$ the proposition follows.\Endproof\vskip4pt

For small enough $\rho$ the proposition guarantees that we have an
inclusion induced map $\cV_{\rho}/(\cV_{\rho} \cap \cB^-)
\longrightarrow \cV_{\rho}/\cV_{\rho}^-$.  A homotopy inverse for this
map can again be found by following the flow. Define an ``exit time''
\begin{displaymath}
    t_*(u) \isdef \inf\{t\geq 0 \mid \Psi^t(u)\in\cB^-\}.
\end{displaymath}
If we allow $t_{*}(u)=\infty$ in case the orbit $\Psi^{t}(u)$ never
hits $\cB^-$ then $t_{*}(u)$ depends continuously on $u\in\cU$, again
because orbits cross $\cB^{-}$ in a topologically transverse way (this
is the same argument as in Proposition \ref{prop:Tcontinuous}).
\begin{proposition}
  If $w_{-}(u)>0$ then $t_{*}(u)<\infty$.
\end{proposition}

\Proof   Let $k_{0}$ be the smallest integer with $u_{k_0}\neq 0$. Then
$\Psi^{t}(u) = \varepsilon \frac{u(t, \cdot)}{\|u(t, \cdot)\|_{2,2}}$
with $u(t, x)$ given by \eqref{eq:fourierheatflow}. For large $t$ the
dominant term in \eqref{eq:fourierheatflow} is the term with
$k=k_{0}$, so that
\begin{displaymath}
    \lim_{t\to\infty} \Psi^{t}(u) = 
       \mathrm{Const}\cdot \Re\lp u_{k_0} e^{2\pi i k_0 x/q} \rp.
\end{displaymath}
Since $w_{-}(u)>0$ we have $k_{0}<p$, and hence for large $t$,
$\Psi^{t}(u)$ has less than $2p$ sign changes so that $\Psi^t(u)$
cannot lie in $\bar{\cB}$ anymore. The only way $\Psi^t(u)$ can leave
$\cB$ is by crossing $\cB^-$ first. \Endproof\vskip4pt   

We can now define the following family of maps,
\begin{displaymath}
    G_{\theta}(u) = \Psi^{\theta\eta(w_-(u)) t_*(u)}(u)
\end{displaymath}
in which $\eta:\R\to \R$ is a continuous nondecreasing function with
$\eta(w)\equiv 0$ for $w\leq\rho/2$ and $\eta(w)\equiv 1$ for $w\geq
\rho$.  Thus $\eta(w_-(u))$ vanishes in the region $w_-(u)\leq \rho/2$
while $t_*(u)$ is continuous for $w_-(u)>0$ so that the product
$\eta(w_-(u)) t_*(u)$ is continuous everywhere.

The $G_{\theta}$ are maps of the pairs $(\cV_{\rho}, \cV_{\rho}^-)$
and $(\cV_{\rho}, \cV_{\rho} \cap \cB^-)$ respectively, and the final
map $G_1$ sends $(\cV_{\rho}, \cV_{\rho}^-)$ to $(\cV_{\rho},
\cV_{\rho} \cap \cB^-)$.  It therefore provides a homotopy inverse for
the inclusion induced map $\cV_{\rho}/( \cV_{\rho} \cap \cB^-) \to
\cV_{\rho}/\cV_{\rho}^-$.

\Subsubsec{Computation of the homotopy type of
  $\cV_{\rho}/\cV_{\rho}^-$} Define $\cW_\rho$ and $\cW_\rho^-$ as
above in \eqref{eq:cWdef1}, \eqref{eq:cWdef2}.
\begin{proposition}
  For small enough $\rho>0$ one has
    \begin{displaymath}
        \cW_\rho = \overline{\cV_{\rho}\setminus\cV_\rho^{-}}, \qquad
        \cW_\rho^- = \cW_\rho \cap \cV_{\rho}^{-}.
    \end{displaymath}
    Consequently{\rm ,} for small $\rho>0$ one has $[\cV_\rho/\cV_\rho^-] =
    [\cW_\rho/\cW_\rho^-]$.
\end{proposition}

\Proof   This follows directly from the proof of Proposition
\ref{prop:vrho}.  \Endproof\vskip4pt

\begin{proposition}
  The pair $(\cW_{\rho}, \cW_{\rho}^-)$ contains $(\cZ_{\rho},
  \cZ_{\rho}^-)$ with
    \begin{displaymath}
        \cZ_\rho \isdef \{ u\in\cW_\rho \mid w_+(u)=0\},\quad 
        \cZ_\rho^- \isdef \cZ_\rho \cap \cW_{\rho}^{-}
    \end{displaymath}
    as a deformation retract.
\end{proposition}

\Proof   One can write $u\in \cW_\rho$ as $u= u^{-}+u^{p}+u^{+}$ and can
homotope it to $G_{\theta}(u) =
\mu(\theta)u^{-}+\nu(\theta)u^{p}+\theta u^{+}$, where $\mu(\theta)$,
$\nu(\theta)\in\R^{+}$ are chosen so as to keep $G_{\theta}(u)$ on
$\cS_{\varepsilon}$. \Endproof\vskip4pt   

\begin{proposition}
  $\cZ_\rho/\cZ_\rho^-$ is homeomorphic with $S^{1}\times S^{2p-1}
  /(S^{1}\times \mathrm{pt})$.
\end{proposition}
\Proof   We can write any $u\in\cZ_\rho$ as
\begin{equation}
    u(x) = \Re\;\textstyle{\sum_{k\leq p} u_k e^{2\pi i kx/q}}, 
    \label{eq:uincZ}
\end{equation}
with
\begin{equation}
    \textstyle{\sum_{k\leq p} (1+k^4)|u_k|^2} = \varepsilon^2.
    \label{eq:uonSeps}
\end{equation}
The condition $w_{-}(u)\leq\rho$ is equivalent to
\begin{equation}
    \textstyle{\sum_{k< p} (1+k^4)|u_k|^2} \leq \rho\varepsilon^2
    \label{eq:wminlessnrho}
\end{equation}
so that
\[
(1+p^{4})|u_{p}|^{2} \geq (1-\rho)\varepsilon^2.
\]
In particular, $u_p\neq 0$ if $\rho<1$, and we can write
\begin{equation}
    u_p = e^{i\theta}\sqrt{\frac{\varepsilon^2-\sum_{k< p} (1+k^4)|u_k|^2}
                                {1+p^4}}
    \label{eq:uppolar}
\end{equation}
with $\theta=\arg u_p$.  We have defined a map $f$ from $\cZ_{\rho}$
to $\C \times \R \times \C^{2p-2}$, given by
\[
f : u \mapsto (e^{i\theta}, u_0, u_1, \dots , u_{p-1}).
\]
(Recall that $u_0$ is real, while the other coefficients are complex.)

This map is one-to-one and hence a homeomorphism onto its image. The
image is clearly $S^1\times B^{2p-1}$, where $S^1$ is the unit circle
in $\C$ and $B^{2p-1}$ is the convex ball in $\R \times \C^{p-1}$
given by \eqref{eq:wminlessnrho}.

The subspace $\cZ_\rho^{-}$ consists of those $u\in\cZ_\rho$ for which
one has equality in \eqref{eq:wminlessnrho}, and therefore $f$ maps
$\cZ_\rho^-$ onto $S^1 \times \partial B^{2p-1}$.  We conclude that
$\cZ_{\rho}/\cZ_\rho^-$ is homeomorphic with $(S^1 \times B^{2p-1}) /
(S^1 \times \partial B^{2p-1})$ which in turn is homeomorphic with
$(S^1\times S^{2p-1})/(S^1 \times \{\mathrm{pt}\})$.
\hfill\qed

\Subsec{A long exact sequence relating the $h^{I}(\cB)$}
Let $\emptyset\subset J\subset I\subset \{1, \dots , m\}$ with $J\neq
I$ be given, and set $K=I\setminus J$.

Choose a metric $g\in \cM_\Gamma(\cB;J)$.  This metric can be modified
to a new metric $\tilde g\in \cM_\Gamma(\cB;I)$ so that $g$ and
$\tilde{g}$ coincide on an open neighborhood of the geodesics
$\gamma_{j}$, for all $j\in J$.

We can then construct punctured neighborhoods $\cU_{i}\subset\hB$ of
$q_{i}\gamma_i$ which isolate the $q_{i}\gamma_i$ for all $i\in I$,
and which are so small that \CS\ for $g$ and for $\tilde{g}$ coincide
on a neighborhood of $q_{i}\gamma_i$ in $\Omega$ for all $i\in J$.

The indices $h^I(\cB)$ and $h^J(\cB)$ are then defined to be the
homotopy types of the pointed spaces
\begin{displaymath}
    \cH^I(\cB) \isdef \frac{\hB}{\cU^I\cup\cB^-\cup\hBl(g, \varepsilon)}, \quad
    \cH^J(\cB) \isdef \frac{\hB}{\cU^J\cup\cB^-\cup\hBl(g, \varepsilon)},
\end{displaymath}
and since $\cU^J\subset\cU^I$ we have a natural map $\cH^J(\cB)\to
\cH^I(\cB)$ which collapses the set
\begin{equation}\label{eq:getscollapsed}
    \cA^I_J \isdef \frac{\cU^I\cup\cB^-\cup\hBl(g, \varepsilon)}{\cU^J\cup\cB^-\cup\hBl(g, \varepsilon)}
\end{equation}
to the base point in $H^{I}(\cB)$.  Since
$\cU^{I}=\cU^{J}\sqcup\cU^{K}$ is a disjoint union, the space
$\cA^I_J$ in \eqref{eq:getscollapsed} is homeomorphic to
\begin{align*}
  \frac{\cU^I\cup\cB^-\cup\hBl(g, \varepsilon)}
  {\cU^J\cup\cB^-\cup\hBl(g, \varepsilon)}
  &=\frac{\cU^J\cup\cU^K\cup\cB^-\cup\hBl(g, \varepsilon)}
  {\cU^J\cup\cB^-\cup\hBl(g, \varepsilon)}\\
  &=\frac{\cU^K} {\cU^K\cap\left(\cU^J\cup\cB^-\cup\hBl(g,
      \varepsilon)\right)}\\
  &=\frac{\cU^{K}}{\cU^{K}\cap\cB^{-}}\\
  &\cong \bigvee_{k\in K}\left\{ \frac{S^{1}\times S^{2p_k-1}}
    {S^1\times\{\textrm{pt}\} } \right\}.
\end{align*}

\begin{proposition}
  The subset $\cA^{I}_{J}$ of $\cH^{J}(\cB)$ is collared.
\end{proposition}
\Proof   We should have started with neighborhoods $\cV_i$, and then
chosen $\cU_{i}\subset\textrm{int}\cV_{i}$.  The curve shortening flow
then retracts the $\cV_{i}$ into the $\cU_{i}$.  \Endproof\vskip4pt   

This proposition implies an isomorphism
\[
H_{l}(\cH^{J}(\cB), \cA^{I}_{J}) \cong
H_{l}(\cH^{J}(\cB)/\cA^{I}_{J})=H_{l}(\cH^{I}(\cB))
\]
of relative singular homology groups.

The long exact sequence on homology for the pair $\left(\cH^{J}(\cB),
  \cA^{I}_{J} \right)$ then gives us the long exact sequence
\begin{equation}
    \dots  H_{l+1}(h^I(\cB)) \stackrel{\partial_*}{\longrightarrow}
    H_{l}(\cA^I_J)\longrightarrow
    H_{l}(h^J(\cB))\longrightarrow
    H_{l}(h^I(\cB))\stackrel{\partial_*}{\longrightarrow}
    H_{l-1}(\cA^I_J)\dots 
\label{eq:exactsequence}
\end{equation}
from Theorem \ref{thm:main2}.

\Subsec{Proof of Theorem {\rm \ref{thm:main3}}}
We know that not all homology groups of $\cA^I_J$ are trivial; so, if
$h^{I}(\cB)$ is the homotopy type of a point, then the exact sequence
implies that $\cA^I_J$ and $h^J(\cB)$ have the same homology groups.
Similarly, if $h^J(\cB)$ happens to be trivial, then
$H_l(\cA^I_J)\cong H_{l+1}(h^I(\cB))$ for all $l$, so that $h^I(\cB)$
cannot be trivial.

\section{Existence theorems for closed geodesics}
\label{sec:existence}
\vglue-12pt
\Subsec{Proof of Theorem {\rm \ref{thm:main1}}}
Fix $\Gamma=\{\gamma_1, \dots , \gamma_N\} \subset\Omega$, a relative
flat knot type $\cB\subset \Omega\setminus\Delta(\Gamma)$, an
$I\subset\{1, \dots , m\}$, and a metric $g\in\cM_\Gamma(\cB;I)$.
Assuming that there are no closed geodesics in $\cB$ for the metric
$g$ we will show that $h^I(\cB)$ is trivial.

Define $T_*:\hB\to [0, \infty]$ by
\[
T_*(\alpha) = \inf\{t\geq0 \mid \Phi^t(\alpha)\in \cB^-\cup \hBl(g,
\eps) \}.
\]
It was shown in Proposition \ref{prop:Tcontinuous} that $T_*$ is
continuous. Thus the set
\[
W= \{\alpha\in\hB \mid T_*(\alpha)=\infty\}
\]
is closed in $\hB$.

Choose neighborhoods $\cU_i\ni q_i\gamma_i$ with $\cU_i=\cU^\#_i$, as
in \S\ref{sec:virtualsatellites}. Let $\cU=\boldsymbol{\cup}_{i\in I} \cU_i$.

For each $\alpha\in W$ there is a $t_\alpha\in[0, \infty)$ such that
$\Phi^{t_\alpha}(\alpha) \in \mathrm{int}\,\cU$, where $\mathrm{int}\,
\cU$ is the interior of $\cU$ with respect to $\hB$. Indeed, if
$\alpha\in W$ then the entire orbit $\Phi^{[0, \infty)}(\alpha)$ is
contained in $\cB$. This orbit must converge to some closed geodesic,
and by assumption such a geodesic must lie on $\pt\bar\cB$. That is, the
orbit must converge to one of the $q_i\gamma_i$ with $i\in I$.

Continuity of the semiflow implies that some neighborhood $\cO_\alpha
\ni\alpha$ also gets mapped into $\mathrm{int}\,\cU$ under
$\Phi^{t_\alpha}$. Choose a sequence $\alpha_n$ so that the $\cO_n =
\cO_{\alpha_n}$ form a locally finite covering of $W$. Next let
$\cO=\boldsymbol{\cup}_n \cO_n$ and construct a continuous function $t_0:\cO \to
[0, \infty)$ with $t_0(\beta)\geq t_{\alpha_n}$ for all
$\beta\in\cO_n$. We may assume that $\lim_{\gamma\to\pt\cO}t_0(\gamma)
= \infty$ (add $\mathrm{dist}(\gamma,\pt\cO)^{-1}$ to $t_0(\gamma)$ if
necessary).

One has $\Phi^{t_0(\beta)}(\beta)\in \mathrm{int}\,\cU$ for all
$\beta\in \cO$. Moreover,
\[
T_\# (\alpha) = \min\bigl(t_0(\alpha), T_*(\alpha)\bigr)
\]
defines a continuous everywhere finite function on $\hB$ which
satisfies
\[
\Phi^{T_\#(\alpha)}(\alpha) \in \cB^- \cup \hBl(g, \eps) \cup \cU.
\]
The family of maps
\[
F^\theta(\alpha) = \Phi^{\theta T_\#(\alpha)}(\alpha)
\]
with $0\leq \theta\leq 1$ defines a deformation retraction of $(\hB,
\cB^- \cup \hBl \cup\cU)$ into $(\cB^- \cup \hBl \cup\cU, \cB^- \cup
\hBl \cup\cU)$. Thus the index $h^I(\cB)$ is trivial.

\Subsec{Proof of Theorem {\rm \ref{thm:main4}}}
Let $\gamma$ be a  simple closed geodesic on the sphere~$S^2$. After
applying a diffeomorphism we may assume that $\gamma$ is the equator.
We consider $p,q$ satellites of the equator. Thus in the notation we
have used so far, we have $\Gamma=\{\gamma_1, \dots , \gamma_N\} =
\{\zeta\}$, $N=1$. The unique curve 
$\zeta\in\Gamma$ belongs to the
boundary of $ \cB_{p,q}(\zeta)$, and thus $m=1$. There are two
\pagebreak
modified Conley indices to be considered, namely $h^\emptyset(\cB)$
and $h^{\{1\}}(\cB)$.

To compute the Conley indices $h^I(\cB_{p,q}(\zeta))$ for arbitrary
$p/q\neq 1$ we use the fact that the indices do not depend on the
metric, and consider the standard metric on the usual unit sphere
$S^2\subset \R^3$. For this metric the equator is indeed a closed
geodesic, while all geodesics are great circles. In particular, no
closed geodesic on the standard sphere is a $p,q$ satellite of the
equator.  Moreover, the rotation number of the equator is exactly
$\rho(\zeta)=1$.

For $p/q>1$ we therefore conclude that
\[
h^{\{1\}}(\cB_{p,q}(\zeta)) = [\text{\textup{point}}]
\]
while for $p/q<1$ we get
\[
h^{\emptyset}(\cB_{p,q}(\zeta)) = [\text{\textup{point}}].
\]
By the long exact sequence from Theorem \ref{thm:main3} we then find
that for $p/q>1$ the index $h^\emptyset(\cB_{p,q}(\zeta))$ is
nontrivial, while for $p/q<1$ the index $h^{\{1\}}(\cB_{p,q}(\zeta))$
is nontrivial.

If one now has another metric $g$ for which the simple closed curve
$\zeta$ is a geodesic with rotation number $\rho(\zeta, g)>p/q>1$,
then the nontriviality of $h^{\emptyset}(\cB_{p,q}(\zeta))$ implies
existence of at least one closed geodesic of $g$ which is a $(p,q)$
satellite of $\zeta$.  Similarly, if one has $1>p/q>\rho(\zeta, g)$,
then nontriviality of $h^{\{1\}}(\cB_{p,q}(\zeta))$ again leads to the
same conclusion.

\section{Appendices}
\vglue-12pt

\Subsec{Curve shortening\ in local coordinates}
\label{sec:CS-in-coords}
Assume $g$ is an $h^{2, \mu}$ metric on $M$ and let $\gamma\in\Omega$
be an $h^{2,\mu}$ curve of length $L$.  Then there exists an
$h^{2,\mu}$ local diffeomorphism $\sigma : \T\times(-r, r)\to M$ with
$\T=\R/L\Z$ such that $x\mapsto \sigma(x, 0)$ is an arclength
parametrization of $\gamma$.  If $\gamma$ is a $q$ fold cover, then we
may assume that $\sigma(x+L/q, y)\equiv \sigma(x, y)$.

In the local coordinates $\{x, y\}$ the metric $g$ is given by
\begin{equation}
        \sigma^* g=E(x, y) (dx)^{2} + 2F(x, y)\, dx\,dy + G(x, y)(dy)^{2}
        \label{eq:ap1}
\end{equation}
where $E, F, G$ are $h^{2, \mu}$ functions on $\T\times (-r, r)$.

We now compute the geodesic curvature of the graph of $y=u(x)$ and
determine the PDE which is equivalent to \CS\ in the coordinates $\{x,
y\}$.

The unit tangent to the graph $\{(x, u(x)) \mid x\in\T\}$ is
\begin{align*}
  T&= \frac{\ptx+u_{x}\pty}{|\ptx+u_{x}\pty|}= \frac{1}{\lambda}(\ptx+u_{x}\pty)\\
\end{align*}
where $\lambda=\sqrt{E+ 2F\, u_{x} + G \,u_{x}^{2}}$.

If we write $X\wedge Y$ for $\Omega_{g}(X, Y)$ where
$\Omega_{g}=(EG-F^{2})dx\wedge dy$ is the area form of the metric $g$,
then the geodesic curvature is $\kappa= T\wedge\nabla_{T}(T)$; i.e.,
\begin{align*}
  \kappa
  &= \lambda^{-3} (\ptx\wedge\pty) \lk u_{xx} + P(x, u) +Q(x, u)u_{x}+
  R(x, u)u_{x}^{2} + S(x, u) u_{x}^{3} \rk
\end{align*}
where
\begin{alignat*}{4}
  &P = \frac{\ptx\wedge\nabla_{\ptx}(\ptx)}{\ptx\wedge\pty}, &&\qquad Q =
  \frac{2\ptx\wedge\nabla_{\pty}(\ptx) +\pty\wedge\nabla_{\ptx}(\ptx)}
  {\ptx\wedge\pty}, \\
  &R = \frac{\ptx\wedge\nabla_{\pty}(\pty)
    +2\pty\wedge\nabla_{\pty}(\ptx)} {\ptx\wedge\pty}, &&\qquad S =
  \frac{\pty\wedge\nabla_{\pty}(\pty)}{\ptx\wedge\pty}. 
\end{alignat*}
If we now consider a moving family of graphs $y=u(x, t)$, then the
normal velocity of this family of curves is given by
\begin{align*}
  V & = T\wedge (u_{t}\pty) = \lambda^{-1} (\ptx\wedge\pty) u_{t}
\end{align*}
so that \CS, i.e.\ $V=\kappa$, is equivalent to
\begin{equation}
        u_{t} = \frac{u_{xx} + P(x, u) +Q(x, u)u_{x}+ R(x, u)u_{x}^{2}
                          + S(x, u) u_{x}^{3}}
                          {E(x, u)+ 2F(x, u)\, u_{x} + G(x, u) \,u_{x}^{2}}.
        \label{eq:CSloc}
\end{equation}

\Subsec{Interpretation of the coefficients $P${\rm ,} $Q${\rm ,} $R${\rm ,} $S$}
The coefficient $S(x, u)$ is the geodesic curvature of the vertical
lines $y=$ constant.  If the diffeomorphism $\sigma$ were obtained by
exponentiating normal vectors to the $x$-axis, as in
\eqref{eq:satellite-def} \S\ref{sec:BirkhoffOrbits}, then $S$ would
vanish. (However, \eqref{eq:satellite-def} contains the unit normal
vector $\norm$ which is only $h^{1,\mu}$, and so the resulting map
$\sigma$ is also only $h^{1,\mu}$ instead of $h^{2,\mu}$.)

The coefficient $P$ is proportional to the geodesic curvature of the
curves $u=$ constant.  In particular $P(x, 0)$ is proportional to the
geodesic curvature of the $x$-axis.  The $x$-axis is a geodesic only
if $P(x, 0)\equiv 0$.

If the $x$ axis is a geodesic then we can assume after an $h^{2,\mu}$
change of coordinates that $x\mapsto(x, 0)$ is a unit speed
parametrization of the $x$-axis, and that on the $x$-axis the vector
$\pty$ is a unit normal to the $x$-axis.  In other words, we assume
that
\[
E(x, 0) =1, F(x, 0) = 0,\quad  G(x, 0) =1
\]
for all $x$. (Use a Whitney type exension theorem, as in \cite[\S
VI.2.3, Th.~4]{Stein}.)

The derivative $P_{y}(x, 0)$ is then given by
\begin{align*}
  P_{y}(x, 0) &= \frac{\pt}{\pt y}
  \lp\frac{\ptx\wedge\nabla_{\ptx}(\ptx)} {\ptx\wedge\pty}
  \rp_{y=0} \\
  &= \ptx\wedge\nabla_{\pty}\nabla_{\ptx}(\ptx)\hskip48pt \textrm{(use~}
  \nabla_{\ptx}\ptx=0, \ptx\wedge\pty=1 \textrm{~for~} y=0)\\
  &= \ptx\wedge\lb\nabla_{\ptx}\nabla_{\pty}(\ptx) + \cR(\pty,
  \ptx)\ptx\rb\qquad \textrm{(definition of the  }
\\[-4pt]
&\hskip2.85in \textrm{Riemann tensor)}\\
  &= K(x, u)
\end{align*}
where $K$ is the Gauss curvature. (This last calculation is the
standard derivation of the equation for Jacobi fields.)

On the $x$-axis we have $\nabla_{\pty}(\ptx) = \nabla_{\ptx}(\pty)=0$
since $\pty$ is a unit normal to a geodesic. We also have
$\nabla_{\ptx}(\ptx)=0$ since the $x$ axis is a geodesic.  Thus
\[Q(x, 0)=0.\]
Both $R$ and $S$ are $h^{1, \mu}$ functions of their arguments.

Thus the linearization of \eqref{eq:CSloc} at $u=0$ is
\begin{equation}
        u_{t}=u_{xx}+K(x) u,
        \label{eq:CSloclin}
\end{equation}
$K(x)=K\circ\gamma(x)$ being the Gauss curvature along the $x$-axis.

\Subsec{Short time existence and the $C^{1}$ local semiflow property}
Equation~\eqref{eq:CSloc} is of the form
\begin{equation}
        u_{t} =  F(x, u, u_{x}, u_{xx})
        \label{eq:Nonlinpde}
\end{equation}
where $F$ is a $C^{1, \mu}$ function of its arguments with
\[
\lambda^{-1}\leq (1+p^2)\frac{\pt F(x, u, p, q)}{\pt q}\leq \lambda
\]
for some constant $\lambda$.  It is well-known (perhaps under higher
differentiablity assumptions on $F$) that solutions with initial data
$u(\cdot, 0)\in C^{2,\mu}(\T)$ exist on a short time interval.  We now
show that \eqref{eq:Nonlinpde} generates a $C^1$ local semiflow on an
open subset of $h^{2,\mu}(\T)$.

We may assume in our setting that $F(x, u, p, q)$ is defined for all
$(x, u, p, q)\in\T\times\R^3$ with $|u|\leq r$ for some $r>0$. Let
$V\subset C^1(\T)$ be defined by
\[
V=\{ u\in C^1(\T) \mid |u|<r \}.
\]
We write $V^{k, \lambda}$ for $V\cap h^{k, \lambda}(\T)$.

The PDE \eqref{eq:Nonlinpde} is actually quasilinear; i.e., $F$ has the
form $F(x, u, p, q) = a(x, u, p)q+b(x, u, p)$ where $a$ and $b$ are
$C^{1,\mu}$ in $x$ and $u$ and analytic in $p\in\R$.  This implies
that the substitution operator $u\mapsto F(x, u, u_x, u_{xx})$ is
continuously Fr\'echet differentiable from $V^{2, \mu}$ to $h^{0,
  \mu}(\T)$.  Since the Fr\'echet derivative of $F$ is the generator
of an analytic semigroup in $h^{0, \nu}(\T)$ for any $\nu \in(0,
\mu)$, we can apply \cite[Cor.~2.9]{Ang-semiflows} and conclude
that \eqref{eq:Nonlinpde} generates a $C^1$ local semiflow on $V^{2,
  \nu}$ for every $\nu\in(0, \mu)$.  This, by definition, means the
following:

\demo{Continuous local semiflow} The map $\Phi$ which maps the
initial data $u_{0}$ and time $t$ to the solution $u(t)$ at time $t$
is defined on an open subset $\cD\subset V^{2, \nu} \times [0,
\infty)$ containing $V^{2, \nu}\times\{0\}$ and satisfies
\begin{enumerate}
\item $F$ is continuous,
\item $F(u_{0}, 0)=u_0$ for all $u_0\in V^{2,\nu}$,
\item If $(u_0, t)\in\cD$ then $\{u_0\}\times[0, t]\subset \cD$,
\item If $(u_0, t)\in\cD$ and $(F(u_0, t), s)\in\cD$ then $(u_0,
  t+s)\in\cD$ and $F(u_0, t+s) = F(F(u_0, t), s)$.
\end{enumerate}

\emph{Differentiable local semiflow.\footnote{We repeat these
    definitions here because there seems to be no consensus on what a
    differentiable local semiflow should be.  In particular Amann
    \cite{amann1}, \cite{amann2} does not include or prove strong
    continuity of $d\Phi_t$ at $t=0$.} }%
For each $t\geq 0$ define $\cD_t = \{u\in V^{2, \nu} \mid (u,
t)\in\cD\}$ and write $\Phi_t(u) = \Phi(u, t)$.  Then the map
$u\mapsto \Phi_t(u)$ is continuously differentiable from $\cD_t$ to
$V^{2,\nu}$.  Moreover, the Fr\'echet derivative $d\Phi_t(u)$ is a
strongly continuous function of both variables $(u, t)\in\cD$, i.e.\ 
for any $v_{0}\in h^{2,\nu}(\T)$ the map $(u, t)\mapsto
d\Phi_t(u)v_{0}$ is continuous from $\cD$ to $h^{2, \nu}(\T)$.  One
obtains $d\Phi_t(u)v_{0}$ by formally linearizing
\eqref{eq:Nonlinpde}; i.e., $v(t) = d\Phi_t(u)v_{0}$ is the solution
of
\begin{eqnarray}\label{eq:linpde}
            v_t &= &F_q(x, u, u_x, u_{xx})v_{xx} + F_p(x, u, u_x,
          u_{xx})v_x +
          F_u(x, u, u_x, u_{xx})v ,\\
          v(\cdot, 0) &=& v_0(\cdot).
        \nonumber
\end{eqnarray}

\Subsec{Linearization at a closed geodesic}
If the curve $\gamma$ (the $x$-axis) is a geodesic so that $u\equiv 0$
is a solution to \eqref{eq:CSloc}, then $u=0$ is a fixed point of the
local semiflow $\Phi_t$ on $V^{2,\nu}$. The semigroup property
$\Phi_t\circ\Phi_s = \Phi_{t+s}$ and the chain-rule imply that the
linear operators $\{d\Phi_t(0) \mid t\geq 0\}$ form a $(C_0)$
semigroup on $h^{2, \nu}(\T)$. Since the linearized equation
\eqref{eq:linpde} coincides with \eqref{eq:CSloclin} the semigroup
$\{d\Phi_t(0)\}$ is generated by $\cA= \lp\frac{d}{dx}\rp^2 -
K\circ\gamma(x)$. ($\cA$ is an unbounded operator on $h^\nu$ with
domain $h^{2,\nu}$ and hence generates a semigroup on $h^\nu$,
$h^{2,\nu}$ and any of their interpolation spaces.)

Let the spectrum of $\cA$ be
\[
\lambda_0>\lambda_1\geq
\lambda_2>\cdots>\lambda_{2i-1}\geq\lambda_{2i}>\cdots
\]
with corresponding eigenfunctions $\{\varphi_k\}$.  For $j\in\N$ we
write $E_j $ for $\mathrm{span}\, \{\varphi_0, \dots\break\dots , \varphi_{2j}\}$
and $E_j^c$ for the closure in $h^{2,\nu}(\T)$ of the span of
$\{\varphi_{2j+1}, \varphi_{2j+2}, \dots \}$. Then $E_j$ and $E_j^c$
are spectral subspaces of the operator $\cA$ with $h^{2,\nu}(\T) =
E_j\oplus E_j^c$.  We let $\pi_j$ denote the projection of
$h^{2,\nu}(\T)$ onto $E_j$ along $E_j^c$.

\begin{lemma}\label{lemma:0or1} \hskip-4pt
  Let $\{u(t)\!  \mid\!  t\!\geq\!0\}\subset V^{2,\nu}$ be an orbit of $\Phi_t$
  with \hbox{$\lim_{t\to\infty} u(t)\!=\!0$} in the $h^{2,\nu}(\T)$ norm. Then
  for any $j\in\N${\rm ,}
\begin{equation}
        \lim_{t\to\infty}\frac{\|\pi_j u(t)\|}{\|u(t)\|} = 0\textrm{~or~}1.
        \label{eq:0or1}
\end{equation}
Here all norms are $h^{2,\nu}(\T)$ norms.

The same statement is true for \/{\rm ``}\/ancient orbits\/{\rm ''}\/ $\{u(t) \mid
-\infty<t\leq0\}$ provided all limits are taken for $t\to-\infty$.
\end{lemma}
\Proof   Once one gets away from the PDE and considers $\{u(n)\mid n=1,
2, \cdots\}$ as an orbit of the time-one map $\Phi_1$ the proof is
completely standard.

The map $\Phi_{1}:\cD_1\to V^{2,\nu}$ is $C^1$, and its Fr\'echet
derivative is given by $d\Phi_1(0) = e^{\mathcal{A}}$, a compact
operator with spectrum $e^{\lambda_{i}}$, $i\in\N_0$.  One can find
equivalent norms $|||\cdot|||$ on $E_j$ and $E_j^c$ so that
\begin{alignat*}{4}
  |||e^{\mathcal{A}} v|||&\geq e^{\lambda_{2j}}|||v|||, &&\qquad \forall v\in E_j ,\\
  |||e^{\mathcal{A}} v|||&\leq e^{\lambda_{2j+1}}|||v|||,&& \qquad
  \forall v\in E_j^c.
\end{alignat*}
Suppose now that
\[
\limsup_{n\to\infty}\frac{|||\pi_j u(n)|||}{|||u(n)|||}>0.
\]
Then for some $\varepsilon>0$ and any $r>0$ there exists a large $n_*$
such that $\frac{|||\pi_j u(n_*)|||}{|||u(n_*)|||}\break>\varepsilon$ and
$|||u(n)|||<r$ for all $n\geq n_*$.

We can now write $\Phi_{1}(u) = \cM(u)u$ where
\[ \cM(u) = \int_0^1 d\Phi(\theta u)d\theta.
\]
Since $\Phi_1$ is $C^1$ we have $|||\cM(u)-e^A|||<\sigma(|||u|||)$
where $\sigma(r)\searrow 0$ as $r\searrow 0$.

If one splits $u=v\oplus v^c \in E_j\oplus E_j^c$, as well as $\bar u
= M(u)u = \bar v\oplus \bar v^c$, and if one assumes $|||u|||\leq r$,
then
\begin{align*}
  |||\bar v ||| &= |||\pi_{j}\mathcal{M}(u)(v\oplus v^c)||| \\
  &\geq |||\pi_{j}e^{\mathcal{A}}(v\oplus v^c)||| - \sigma(r) |||u||| \\
  &\geq e^{\lambda_{2j}} |||v||| - \sigma(r) |||u|||\\
  &\geq \lp e^{\lambda_{2j}}-\sigma(r)\rp|||v|||-\sigma(r) |||v^c|||.
  \intertext{Similarly one finds} |||\bar v^c||| &\leq \lp
  e^{\lambda_{2j+1}}+\sigma(r)\rp |||v^c||| + \sigma(r) |||v||| .
\end{align*}
If one also assumes that $|||v|||\geq \varepsilon||| v^c|||$ then
\begin{equation}
    \frac{|||\bar v |||}{|||\bar v^c|||}
        \geq \frac{\lp e^{\lambda_{2j}}-\sigma(r)\rp \varepsilon-\sigma(r)}
                  {e^{\lambda_{2j+1}}+\sigma(r)+\varepsilon\sigma(r)} .
    \label{eq:vvcgrowth}
\end{equation}
Since $e^{\lambda_{2j}} > e^{\lambda_{2j+1}}$ one can choose
$1<\vartheta<e^{\lambda_{2j}-\lambda_{2j+1}}$. For sufficently small
$r>0$ one concludes from \eqref{eq:vvcgrowth} that
$$
\frac{|||\bar v |||}{|||\bar v^c|||}\geq \vartheta \varepsilon .
$$
Inductive application of this estimate shows that for $u(k) =
v(k)\oplus v^c(k)$ one has
$$
\frac{|||v(n_*+i)|||}{|||v^c(n_*+i)|||} \geq \varepsilon(n_*+i) =
\vartheta^i \varepsilon
$$
as long as $1+\varepsilon(n_*+i) < \frac{\delta}{2\sigma(r)}$. \medbreak

Thus if $\limsup_{n\to\infty}\|\pi_j u(n)\|/\|u(n)\|>0$, then
$$\liminf_{n\to\infty}\|\pi_j v(n)\|/\|v^c(n)\|\geq \delta/2\sigma(r),$$
with $r>0$ arbitrarily small. Hence $\|\pi_j u(n)\|/\|u(n)\|\to 1$ as
$n\to\infty$.

Having established the alternative \eqref{eq:0or1} along a sequence
$n\nearrow\infty$ we now assume that
$$\limsup_{t\nearrow\infty}\|\pi_j u(t)\|/\|u(t)\| >\varepsilon>0.$$
Then for any $r>0$ there is a $t_0>0$ such that $\|u(t)\|\leq r$ for
$t\geq t_0$ and $\|\pi_j u(t_0)\|/\|u(t_0)\| >\varepsilon$.

The previous arguments imply that $\lim_{n\nearrow\infty}\|\pi_j
u(t_n)\|/\|u(t_n)\|=1$, where $t_n=t_0+n$.

Splitting $u(t)=v(t)\oplus v^c(t)$ as before we have $\|v(t_n)\|=
o\bigl(\|v^c(t_n)\|\bigr)$ for $n\nearrow\infty$.  To estimate $v(t)$
and $v^c(t)$ for $t\in (t_n, t_{n+1})$ we write
\[
u(t_n+\theta) = d\Phi_\theta(u(t_n)) = \cM_\theta(u(t_n))u(t_n),
\]
where $\cM_\theta(u) = \int_0^1 d\Phi_\theta(su) ds$.

Since $\Phi$ is a differentiable semiflow the map $(\theta, u)\mapsto
\cM_\theta(u)$ is strongly continuous.  Hence, for small enough $r>0$
the operators $\{\cM_\theta(u) \mid 0\leq \theta\leq 1, \|u\|\leq r\}$
are uniformly bounded.  Since $E_{2j}$ is finite dimensional, the map
$(\theta, u)\mapsto \cM_\theta(u)|_{E_{2j}}$ is norm continuous. In
particular, there is a $\tau(r)>0$ with $\tau(r)\searrow 0$ for
$r\searrow0$, such that
\[
\|\cM_\theta(u)|_{E_{2j}} - e^{\theta \cA}|_{E_{2j}}\|_{\Lin(E_{2j},
  h^{2,\nu})} \leq \tau(r)
\]
if $\|u\|\leq r$ and $\theta\in [0, 1]$.

We have the following estimates:
\begin{align*}
  \|v(t_n+\theta)\| &\geq \|\pi_j\cM_\theta(u(t_n))v(t_n)\| -
  \|\pi_j\cM_\theta(u(t_n))v^c(t_n)\| \\
  &\geq \|\pi_j e^{\theta\cA}v(t_n)\| -
  \|\pi_j(\cM_\theta(u(t_n))-e^{\theta\cA})v(t_n)\| \\
&\quad -
  \|\pi_j\cM_\theta(u(t_n))v^c(t_n)\| \\
  &\geq e^{\theta\lambda_{2j}}\|v(t_{n})\| - \tau(r)\|v(t_{n})\| -o(1)
  \|v(t_{n})\|
  & \mbox{} \\
  \intertext{(use $\|v^{c}(t_{n})\|=o(\|v(t_{n})\|)$). Also}
  \|v^c(t_{n}+\theta)\| &\leq \|\pi_j^c\cM_\theta(u(t_n))v(t_n)\| -
  \|\pi_j^c\cM_\theta(u(t_n))v^c(t_n)\| \\
  &\leq \tau(r) \|v(t_{n})\| + o(1) \|v(t_n)\| \\
  &= o(1)\cdot \|v(t_{n})\|
\end{align*}
which together imply $\|v^c(t_{n}+\theta)\|=o(\|v(t_{n}+\theta)\|)$ as
$n\nearrow\infty$, uniformly in $\theta\in [0, 1]$.  \hfill\qed

\begin{lemma}[Notation as in Lemma \ref{lemma:0or1}] \label{lemma:lowerexponential}
    For any solution $\{u(t)
  \mid t\geq0\}$ of \CS\ which converges to $u=0$ there exists a
  $j\in\N$ such that
$$
\lim_{t\to\infty}\frac{\|\pi_j u(t)\|}{\|u(t)\|} = 1.
$$
In particular one has $\|u(t)\|\geq C e^{-\zeta t}$ for some
$\zeta<\infty$.
\end{lemma}

\Proof   If the limit were 0 for all $j$ then the solution $u(t)$ would
approach $u=0$ faster than any exponential, and  so we must prove the lower
bound $\|u(t)\|\geq C e^{-\zeta t}$.

There is a standard approach for proving exponential \emph{lower}
bounds on decay in heat equations due to Agmon (see \cite[\S2.18,
p.181]{Friedman}) which is used to prove backward uniqueness results.
This approach would work here, but it would require us to
differentiate the functions $a(x, u, p)$ in the PDE twice, thereby
forcing us to consider metrics $g$ on $M$ with at least three
derivatives.  In order not to use more than just $g\in h^{2,\mu}$ we
follow the less standard approach from the appendix in
\cite{Ang-MorseSmale} which applies to semilinear equations.

To rewrite \CS\ as a semilinear equation we study the evolution of the
curvature as a function of renormalized arclength.  Let
$\gamma:\R/\Z\times[0, t_{*})\to M$ be a normal parametrization (i.e.\ 
$\pt_{t}\gamma\perp\pt_x\gamma$) of a solution of \CS. Write $L(t)$
for length at time $t$, let $P_t$ be the point $\gamma(0, t)$ (so that 
$P_{t}$ moves with velocity perpendicular to the curve always) and
define the normalized arclength coordinate $\varsigma$ of any point
$Q=\gamma(x, t)$ on $\gamma_{t}$ by
\begin{equation}
        \varsigma(x, t)
                = \frac{1}{L(t)} \int_0^{x} |\pt_x\gamma(\xi, t)|d\xi
                = \frac{1}{L(t)} \int_{P_{t}}^{Q}ds.
        \label{eq:varsigdef}
\end{equation}
We also introduce a new time variable related to $t$ via
$$
\tau = \int_0^t \frac{dt}{L(t)^2}.
$$
\begin{proposition}
  The curvature $\kappa,$ as a function of $\tau$ and $\varsigma,$
  satisfies
\begin{equation}
        \kappa_{\tau}
                =\kappa_{\varsigma\varsigma}
                        + L(\tau)^{2} \lb J[\kappa]\kappa_\varsigma + (K\circ\gamma)\kappa
                        + \kappa^3 \rb
        \label{eq:ksemilin}
\end{equation}
where
$$
J[\kappa] = \int_0^\varsigma \kappa^2d\varsigma - \varsigma \int_0^1
\kappa^2d\varsigma.
$$
\end{proposition}

\Proof   A straightforward calculation begins with differentiating
\eqref{eq:varsigdef} with respect to $t$ to get
$$
\frac{\pt\varsigma}{\pt t} = \varsigma\int_{0}^{1}\kappa^{2}d\varsigma
- \int_0^\varsigma \kappa^{2}d\sigma.
$$
Then the chain rule
$$
\lp\frac{\pt \kappa}{\pt t}\rp_{x=\textrm{const}} = \lp\frac{\pt
  \kappa}{\pt t}\rp_{\varsigma=\textrm{const}} + \frac{\pt \kappa}{\pt
  \varsigma} \frac{\pt \varsigma}{\pt t}
$$
after some simplification leads to \eqref{eq:ksemilin}.  \Endproof\vskip4pt   

Since the limiting geodesic $u=0$ has positive length the new and old
time variables $t$ and $\tau$ are roughly proportional, so it suffices
to establish an exponential lower bound for the solution in the $\tau$
variable.

The equation \eqref{eq:ksemilin} is semilinear, and can be written as
\[
\kappa_{\tau} = \cA \kappa + \cR(\tau)\kappa
\]
where $\cA = (\pt_{\varsigma})^{2} - L_0^2 K_0(x)$, with $K_0(x) =
K(\sigma(x, 0))$, is the Gauss curvature on the $x$ axis, and $L_0$ is
the length of the $x$-axis. The ``remainder'' operator $\cR(\tau)$ is
\[
\cR(\tau) = L(\tau)^2J[\kappa(\tau)] \frac{\pt}{\pt \varsigma} +
\kappa^{2} + \bigl( K\circ\gamma_{\tau} - K\circ\gamma_{\infty}\bigr).
\]
This operator is bounded from the Sobolev space $W^{1, 2}(\R/\Z)$ to
$L^2(\R/\Z)$.  If we assume that $\|\pi_j u(t)\|/\|u(t)\|\to0$ for all
$j$ then the coefficients in $\cR$ decay faster than any exponential
$e^{-\zeta\tau}$ and thus the operator norm of $\cR(\tau)$ from $W^{1,
  2}(\R/\Z)$ to $L^2(\R/\Z)$ also tends to zero.

The eigenvalues of the self-adjoint operator $\cA$ on $L^2$ grow like
$n^2$, so the $n^{\mathrm{th}}$ gap in the spectrum of $\cA$ has
length proportional to $n$. This is exactly enough for the argument in
\cite[Appendix]{Ang-MorseSmale} and we can conclude that no solution
of \CS\ can approach a geodesic at a faster than exponential rate. \phantom{overthere}
\hfill\qed


\references {888}

\bibitem[1]{amann1} \name{H.\ Amann}, Quasilinear evolution equations and
    parabolic systems, {\it Trans.\ Amer.\ Math.\  Soc\/}.\  \textbf{293}
  (1986), 191--227.
  
\bibitem[2]{amann2} \bibline, Dynamic theory of quasilinear
    parabolic equations.\ I.\ Abstract evolution equations, {\it Nonlinear
  Anal\/}.\ \textbf{12} (1988), 895--919.
  
\bibitem[3]{Ang-MorseSmale} \name{S.\ B.\ Angenent}, The Morse-Smale
    property for a semi-linear parabolic equation, {\it J.\ 
  Differential Equations\/} \textbf{62} (1986),  427--442.
   
  
\bibitem[4]{Ang-twist} \bibline, The periodic orbits of an
    area preserving twist map, {\it Comm.\ Math.\ Phys\/}.\ \textbf{115}
  (1988),  353--374.
  
\bibitem[5]{Ang-semiflows} \bibline, Nonlinear
    analytic semiflows, {\it Proc.\ Royal Society of 
    Edinburgh\/} {\bf 115A} (1990), 91--107.
  
\bibitem[6]{Ang1} \bibline, Parabolic equations for
    curves on surfaces.\ I.\  Curves with $p$-integrable curvature,
   {\it Ann.\ of Math\/}.\ {\bf 132} (1990),  451--483.
  
\bibitem[7]{Ang2} \bibline,  Parabolic equations for
    curves on surfaces.\ II.\   Intersections, blow-up, and
    generalized solutions,  {\it Ann.\ of Math\/}.\ {\bf 133} (1991),
    171--215.
  
\bibitem[8]{Ang1991} \bibline, On the formation of
    singularities in the curve shortening flow,  {\it J.\ Differential
  Geometry\/} {\bf 33} (1991),  601--633.

\pagebreak  
\bibitem[9]{ElEscorial} \name{S.\ B.\ Angenent}, Recent results in
    mean curvature flow, in {\it Recent Advances in Partial
    Differential Equations\/} (El Escorial, 1992), 1--18,
    {\it RAM Res.\ Appl.\ Math\/}.\ {\bf 30}, Masson, Paris, 1994.
  
\bibitem[10]{Ang-sagarro} \bibline, Inflection points,
    extatic points and curve shortening,  in {\it Hamiltonian Systems
    Systems with Three or More Degrees of Freedom\/} (S'Agar\'o, 1995), 
    3--10, {\it NATO Adv.\ Sci.\ Inst.\ Ser.\ C Math.\ Phys.\ Sci\/}.\ 
    {\bf 533}, Kluwer Acad.\ Publ., Dordrecht, 1999.
  
\bibitem[11]{Ang-sturm} \bibline,  The zero set of a
    solution of a parabolic equation, {\it J.\  Reine Angew.\ Math\/}.\
   \textbf{390} (1988), 79--96.
   
\bibitem[12]{AngenentFiedler} \name{S.\ B.\ Angenent}
   and \name{B.\ Fiedler},
     The dynamics of rotating waves in reaction diffusion
         equations, {\it Trans.\ Amer.\ Math.\ Soc\/}.\ {\bf 307}
	 (1988),
	     545--568.

\bibitem[13]{Arnold} \name{V.\ I.\ Arnol'd}, {\it Topological Invariants of Plane
    Curves and Caustics,\/} {\it  University Lecture Series\/} 
  \textbf{5},  A.M.S., Providence, RI, 1994.
  
\bibitem[14]{Arnold2} \bibline,  {\it Mathematical Methods of
    Classical Mechanics\/},  {\it Grad.\ Texts in Math\/}.\ {\bf 60},
    Springer-Verlag, New York, 1978.

  
\bibitem[15]{birkhoff1} \name{G.\ D.\ Birkhoff}, {\it Dynamical
    Systems\/}, {\it A.M.S.\ Colloq.\ Publ\/}.\ {\bf IX}, A.M.S.,
    Providence, RI, 1966.
  
\bibitem[16]{birkhoff2} \bibline, Surface
    transformations and their dynamical applications, {\it Acta
  Math\/}.\ \textbf{43} (1922), 1--199.
  
\bibitem[17]{Conley} \name{C.\ Conley},  {\it Isolated Invariant Sets
    and The Morse Index\/}, {\it CBMS Regional Conf.\ Series in
    Math\/}.\ 
   \textbf{38} (1978), A.M.S.\ Providence, RI.
  
\bibitem[18]{Coddington-Levinson} \name{E.\ A.\ Coddington} and \name{N.\
Levinson}, 
  {\it Theory of Ordinary Differential Equations},  McGraw-Hill,
  New York, 1955.
  
  
\bibitem[19]{Friedman} \name{A.\ Friedman}, {\it Partial Differential
    Equations\/}, Holt, Rinehart and Winston, Inc., New York,  1969.

  
\bibitem[20]{Gage1983} \name{M.\ Gage}, An isoperimetric inequality with
    applications to curve shortening,  {\it Duke Math.\ J\/}.\  {\bf
    50} (1983), 1225--1229.
  
\bibitem[21]{Gage1984} \bibline, Curve shortening makes convex
    curves circular, {\it Invent.\ Math\/}.\ {\bf 76} (1984), 357--364.
  
\bibitem[22]{GH1986} \name{M.\ Gage} and \name{R.\ S.\ Hamilton}, The heat equation
    shrinking convex plane curves,  {\it J. Differential Geom\/}.\ {\bf
    23} (1986), 69--96.

  
\bibitem[23]{Gray1987} \name{M.\ Grayson}, \textit{The heat equation shrinks
    embedded plane curves to round points, } {\it J. Differential Geom\/}.\
  {\bf 26} (1987), 285--314.
  
\bibitem[24]{Gray1989} \bibline, Shortening embedded curves, 
  {\it Ann.\  of Math\/.}\ {\bf 129} (1989),  71--111.
  
\bibitem[25]{Ha:isoper} \name{R. S.\ Hamilton}, Isoperimetric estimates for
    the curve shrinking flow in the plane, in {\it Modern Methods in
  Complex Analysis\/} (Princeton, NJ, 1992), 201--222, {\it Ann. of Math.
  Studies\/} \textbf{137}, Princeton Univ.\ Press, Princeton, NJ, 1995.
  
\bibitem[26]{Hu:isoper} \name{G. Huisken}, A distance comparison
    principle for evolving curves, {\it   Asian J. Math\/}.\ \textbf{2} (1998),
   127--133
  
\bibitem[27]{JohnsonMoser} \name{R.\ Johnson} and \name{J.\ Moser}, The rotation number
    for almost periodic potentials, {\it Comm.\ Math.\ Phys\/}.\ 
  \textbf{84} (1982), 403--438.
  
\bibitem[28]{Klingenberg} \name{W.\ Klingenberg}, Lectures on closed
    geodesics,   {\it Grundlehren der mathematischen
  Wissenschaften\/}, \textbf{230}   (1978), Springer-Verlag, New York.
  
\bibitem[29]{Matano} \name{H.\ Matano}, Nonincrease of the lap-number of a
    solution for a one-dimensional semilinear parabolic equation,  {\it J.\
  Fac.\ Sci.\ Univ.\ Tokyo Math}.\ {\bf 29} (1982), 401--441.
  
\bibitem[30]{mather-rotation} \name{J.\  N. Mather}, Amount of rotation
    about a point and the Morse index, {\it Comm.\ Math.\ Phys\/}.\
  \textbf{94} (1984),  141--153.
  
\bibitem[31]{Nickel} \name{K.\ Nickel},  Gestaltaussagen \"uber
    L\"osungen parabolischer Differentialgleichungen, {\it J.\ Reine
  Angew.\ Math\/}.\ \textbf{211} (1962), 78--94.

  
\bibitem[32]{Oaks} \name{J.\ A.\ Oaks},  Singularities and
    self-intersections of curves evolving on surfaces, {\it  Indiana
  Univ.\ Math. J\/}.\ \textbf{43} (1994),  959--981.
  
\bibitem[33]{poincare} \name{H.\ Poincar\'e}, Sur les lignes
    geodesiques des surfaces convexes, {\it  Trans.\ Amer.\ Math.\ Soc\/}.\
  \textbf{6} (1905), 237--274.

 
\bibitem[34]{Stein} \name{E.\ M.\  Stein}, {\it Singular Integrals and
    Differentiability Properties of Functions\/}, {\it Princeton Math.\
    Series\/}  {\bf 30},  Princeton Univ.\
  Press, Princeton, NJ, 1970.

  
\bibitem[35]{Sturm} \name{C.\ Sturm},  M\'emoire sur une classe
    d'\'equations \`a diff\'erences partielles,  {\it J.\ de
  Math{\hskip1pt\rm \'{\hskip-5.5pt\it e}}-\break matiques Pures et 
Appliqu{\hskip1pt\rm \'{\hskip-5.5pt\it e}}es\/}  {\bf 1} (1836), 373--444.


\Endrefs
\end{document}